\documentclass{amsart}

\usepackage{amssymb,amsfonts,amscd,tikz,enumerate,mathtools,tikz-cd,amsmath}
\usepackage[utf8]{inputenc}
\usepackage[colorlinks,linktocpage]{hyperref}
\hypersetup{linkcolor=[rgb]{0,0,0.715}}
\hypersetup{citecolor=[rgb]{0,0.715,0}}
\usepackage[T1]{fontenc}
\usepackage{esint}

\newtheorem{theorem}{Theorem}[section]
\newtheorem{corollary}[theorem]{Corollary}

\newtheorem{lemma}[theorem]{Lemma}

\newtheorem*{thm1}{Theorem A}
\newtheorem*{thm2}{Theorem B}
\newtheorem*{thm3}{Theorem C}

\theoremstyle{definition}
\newtheorem{defn}[theorem]{Definition}

\theoremstyle{remark}
\newtheorem{remark}[theorem]{Remark}

\renewcommand{\subset}{\subseteq}

\makeatletter
\let\c@equation\c@theorem
\makeatother
\numberwithin{equation}{section}

\makeatletter
\renewcommand{\l@section}{%
  \@tocline{1}{4pt}{0pt}{}{\bfseries}}
\renewcommand{\l@subsection}{%
  \@tocline{2}{0pt}{1.5em}{}{}}
\makeatother

\title[]{Fibrations and nilpotent collapsing structure}
\author{Jikang Wang}
\address{Institute of Geometry and Physics, University of Science and Technology of China, Hefei, Anhui, China, 230026}

\email{wangjikang@ustc.edu.cn}
\thanks{}
\thanks{}

\thanks{}

\begin{document}
\date{}

\begin{abstract}
We construct fiber bundles, singular fibrations, and nilpotent
collapsing structures for Riemannian manifolds with lower Ricci
curvature bounds under additional geometric assumptions.

For every $n$ and $D>0$, there exists $\epsilon(n,D)>0$
such that a closed Riemannian $n$-manifold $M$ satisfying
\[
\operatorname{diam}(M)\le D,\qquad
\operatorname{sec}_M\ge-1,\qquad
\operatorname{Ric}_M\ge-\epsilon(n,D)
\] admits a finite cover of degree at most
$C(n)$ that is a smooth fiber bundle over a nilmanifold.
The fibers have finite fundamental group, and the base dimension
equals the Hirsch rank of $\pi_1(M)$. If $\pi_1(M)$ is torsion-free, then $M$
itself fibers smoothly over an infranilmanifold with simply
connected fibers.

We next consider a Gromov--Hausdorff convergent sequence $M_i\to X$
of closed Riemannian $n$-manifolds with uniformly bounded diameter,
$\operatorname{Ric}_{M_i}\ge -(n-1)$, and the $(\rho,\delta)$ local
rewinding Reifenberg condition, where
$\rho>0$ is fixed and $\delta$ is sufficiently small.
For all sufficiently large $i$, we construct singular
fibrations $\bar{F}_i:M_i\to X$ that are $\epsilon_i$-Gromov--Hausdorff
approximations, with $\epsilon_i\to0$, and whose fibers are homeomorphic
to infranilmanifolds. Under the same local hypotheses, complete, possibly noncompact manifolds admit nilpotent collapsing structures on their $\epsilon(n,\rho)$-collapsed regions.

The proofs combine the large-scale structure of approximate
groups with equivariant gluing of almost splitting maps.
We also establish corresponding results in the $\mathrm{RCD}$ setting.

Finally, we give an independent proof of a uniform positive lower
bound on the injectivity radius at the basepoints of local universal
covers of complete manifolds with $|\operatorname{sec}|\le1$. This implies the local rewinding Reifenberg condition. Our methods therefore also provide a new construction
of the singular fibrations and nilpotent collapsing atlases of
classical collapsing theory, with explicit descriptions of their
infranilmanifold fibers.
\end{abstract}
\maketitle

\clearpage
\tableofcontents

\clearpage
\section{Introduction}
We study the topology of manifolds with a lower Ricci curvature bound and, more generally, of $\mathrm{RCD}(-(n-1),n)$ spaces $(X,p,d,\mathfrak{m})$.

For a closed Riemannian manifold $M$ with nonnegative Ricci curvature, the Cheeger--Gromoll splitting theorem \cite{CheegerGromoll1971} implies that $\pi_1(M)$ is virtually abelian and that a finite cover of $M$ is a smooth fiber bundle over a torus with simply connected fibers.

For closed manifolds with bounded diameter and almost nonnegative Ricci curvature, Gromov conjectured that the fundamental group is virtually nilpotent. This conjecture follows from the generalized Margulis lemma of Kapovitch and Wilking \cite{KapovitchWilking2011}; see also \cite{DSZZ2023} for the $\mathrm{RCD}$ setting. Fukaya and Yamaguchi \cite{FukayaYamaguchi1992} proposed
a corresponding fibration statement: a finite cover of $M$ is a fiber bundle over a nilmanifold with simply connected fibers. Anderson, however, constructed a counterexample to this statement \cite{Anderson1992}.

Under almost nonnegative sectional curvature,
Kapovitch--Petrunin--Tuschmann proved that a finite cover
fibers over a nilmanifold with simply connected fibers
\cite[Theorem~C]{KPT2010}. When $\pi_1(M)$ is torsion-free and nilpotent, they further showed that $M$ itself fibers over a nilmanifold with simply connected fibers \cite[Theorem~5.1.1]{KPT2010}. Our first result replaces the almost nonnegative sectional curvature assumption by a fixed lower sectional curvature bound together with almost nonnegative Ricci curvature.

\begin{theorem}\label{theorem:smooth bundle over nil}
For each $n$, there exists a constant $C(n)>0$ such that, for every $D>0$, there exists $\epsilon(n,D)>0$ with the following property. Let $M$ be a closed Riemannian $n$-manifold satisfying
\[
\operatorname{diam}(M)\le D,\qquad
\operatorname{sec}_M\ge -1,\qquad
\operatorname{Ric}_M\ge -\epsilon(n,D).
\]
Then $M$ admits a finite cover of degree at most $C(n)$ that is a smooth fiber bundle over a nilmanifold, with fibers having finite fundamental group. Moreover, if $\pi_1(M)$ is torsion-free, then $M$ itself is a smooth fiber bundle over an infranilmanifold with simply connected fibers.
\end{theorem}

\begin{remark}
The finite cover constructed in the proof of
Theorem~\ref{theorem:smooth bundle over nil} has nilpotent
fundamental group, which has a torsion-free subgroup
of finite index.
Passing to the corresponding cover and the induced cover
of the nilmanifold base gives a smooth fiber bundle with
simply connected fibers. The degree of such a cover,
however, admits no uniform bound.
\end{remark}

By \cite{HHWZ}, the manifold $M$ in Theorem~\ref{theorem:smooth bundle over nil} also admits a fiber bundle over the torus $T^{b_1(M)}$, without any torsion-free assumption on $\pi_1(M)$. The base in our construction has dimension equal to the Hirsch rank of $\pi_1(M)$. Throughout, the rank of a virtually nilpotent discrete group means its Hirsch rank.

We prove Theorem~\ref{theorem:smooth bundle over nil} as a consequence of the following more general result in the $\mathrm{RCD}$ setting.
\begin{thm1}
For every $N\in[1,\infty)$, there exist constants $\delta(N),C(N)>0$ such that, for every $D>0$, there exists $\epsilon(N,D)>0$ with the following property. Let $(X,p,d,\mathfrak{m})$ be an $\mathrm{RCD}(-\epsilon(N,D),N)$ space with Alexandrov curvature bounded below by $-1$ and $\operatorname{diam}(X)\le D$. Suppose that every point of $X$ is $(k_0,\delta(N))$-regular, where $k_0$ is the rectifiable dimension of $X$. Then $X$ admits a finite cover of degree at most $C(N)$ that is a topological fiber bundle over a nilmanifold, with fibers having finite fundamental group. Moreover, if $\pi_1(X)$ is torsion-free, then $X$ itself is a topological fiber bundle over an infranilmanifold with simply connected fibers.
\end{thm1}

The proof of Theorem~A constructs globally defined equivariant maps
that are locally almost splitting, building on \cite{HHWZ,Wang2024}.
The construction in \cite{HHWZ} treats torsion-free abelian groups,
while that in \cite{Wang2024} assumes that the limit of the covering
spaces is Euclidean. Here the approximating groups are virtually
nilpotent, and the covering spaces may converge to a general
$\mathrm{RCD}(0,N)$ space admitting a cocompact nilpotent group action.
To control the group structure, we establish
Theorem~\ref{thm:large nil}, extending earlier results on approximate
groups in \cite{BGT,Zamora2020} to large scales.

We next recall the classical theory of collapsing manifolds with bounded sectional curvature; see \cite{Gromov1978,Ruh1982,CG1,CG2,Fukaya1987,Fukaya1988,Fukaya1989,CFG,Rong2019}.

Gromov proved that an almost flat $n$-manifold admits a finite cover of degree at most $C(n)$ that is diffeomorphic to a nilmanifold. Ruh subsequently showed that every almost flat manifold is itself diffeomorphic to an infranilmanifold. Rong later gave a new proof of the almost flat manifold theorem.

Fukaya proved that when a sequence of manifolds $M_i$ with uniformly bounded sectional curvature collapses to a compact manifold of lower dimension, the collapse is described by fiber bundle maps from $M_i$ to the limit. After smoothing the metrics if necessary, the fibers are almost flat and hence diffeomorphic to infranilmanifolds. When the limit is a compact singular space, Fukaya constructed $O(n)$-equivariant regular fibrations on the orthonormal frame bundles $F(M_i)$. Passing to the quotients yields singular fibrations from $M_i$ to the limit space. The singular behavior is governed by the isotropy groups of the induced $O(n)$-action on the limit of the frame bundles.

Cheeger--Fukaya--Gromov developed this theory further by constructing nilpotent collapsing structures on complete, possibly noncompact, manifolds with bounded sectional curvature. Roughly speaking, the $\epsilon$-collapsed region admits an atlas of compatible local singular fibrations. On suitable overlap neighborhoods, the fibration over a lower-dimensional base factors through the one over a higher-dimensional base. They also constructed a nearby metric admitting a nilpotent Killing structure, with respect to which the local nilpotent actions on suitable covers are isometric.

Rong introduced the $(\rho,\delta)$ local rewinding Reifenberg condition and conjectured that classical collapsing theory extends to manifolds with a lower Ricci curvature bound satisfying this condition. We formulate this condition below in the $\mathrm{RCD}$ setting.

\begin{defn}\label{def:Reifenberg}
Let $\rho,\delta>0$ and $n\in\mathbb{N}$. A pointed metric space $(X,p,d)$ satisfies the $(\rho,\delta)$ Reifenberg condition at $p$ if
\[
d_{\mathrm{GH}}\bigl(B_t(p),B_t(0^n)\bigr)\le t\delta
\qquad\text{for every }0<t\le\rho,
\]
where $0^n$ denotes the origin of $\mathbb{R}^n$.

For an $\mathrm{RCD}(-(n-1),n)$ space $(X,p,d,\mathcal H^n)$, let $\tilde{B}(x,r)$ denote the universal cover of $B_r(x)$ for $x\in X$ and $r>0$, equipped with the lifted length metric and measure. Let $\tilde p$ be a lift of $p$ to $\tilde{B}(p,\rho)$. We say that $X$ satisfies the $(\rho,\delta)$ local rewinding Reifenberg condition at $p$ if $\tilde{B}(p,\rho)$ satisfies the $(\rho/3,\delta)$ Reifenberg condition at $\tilde p$.
\end{defn}

Throughout the paper, covering spaces of open balls carry the lifted metric and measure. These spaces need not be complete, but they satisfy the RCD condition locally in the interior. We use relative covers and shrink the working radii so that the balls required for the local geometric and analytic estimates remain within the covering domain.

Under the local rewinding Reifenberg condition, analogues of Gromov's almost flat manifold theorem and Fukaya's regular fibration theorem have been established for manifolds with a lower Ricci curvature bound \cite{HKRX2020,Huang2020,Rong2022} and, more generally, for $\mathrm{RCD}$ spaces \cite{Wang2024,ZamoraZhu2024}. Some of these results assume weaker conditions that imply the local rewinding Reifenberg condition after suitable rescaling.

Our second main result establishes the topological singular
fibration conclusion predicted by Rong's conjecture in the
$\mathrm{RCD}$ setting.

In Fukaya's construction, a singular fibration is obtained as the quotient of an $O(n)$-equivariant regular fibration on the orthonormal frame bundle. Our construction works directly with local covering spaces.

Let $(X_i,p_i,d_i,\mathcal{H}^n)$ be a sequence of $\mathrm{RCD}(-(n-1),n)$ spaces such that
\[
\left(
X_i,p_i,d_i,
\frac{\mathcal{H}^n}{\mathcal{H}^n(B_1(p_i))}
\right)
\xrightarrow{\mathrm{mGH}}
(X,p,d,\mathfrak{m}).
\]
We call a Gromov--Hausdorff approximation
$\bar F_i:(X_i,p_i)\to(X,p)$ a singular fibration if it is locally the quotient of an orbit-preserving eGH homeomorphism on relative covering spaces, as follows.

For every $x\in X$ and every sequence $x_i\in X_i$ converging
to $x$, there exist $s>0$ and, for all sufficiently large $i$,
pointed local relative covering spaces
$(\tilde U_i,\tilde x_i,G_i)$ of neighborhoods of $B_{3s}(x_i)$,
with deck transformation groups $G_i$ and
$\pi(\tilde x_i)=x_i$, together with
equivariant Gromov--Hausdorff approximations
\[
(F_i,\phi_i,\psi_i):
(\tilde U_i,\tilde x_i,G_i)
\longrightarrow
(Y,\tilde x,G).
\]
For each sufficiently large $i$, the triple
$(Y,\tilde x,G)$ is chosen from the possible equivariant
limits of the relative local covers so that 
\[
B_{2s}(x)\subset U:=Y/G\subset X,
\qquad \pi(\tilde x)=x.
\] 
This choice may depend on $i$, and we suppress this dependence in the notation. 
The map $F_i$ restricts to a homeomorphism
from $B_{2s}(\tilde x_i)$ onto its image and maps the
intersection of each $G_i$-orbit with $B_{2s}(\tilde x_i)$
into a single $G$-orbit. Finally,
\[
\pi\circ F_i=\bar F_i\circ\pi
\qquad\text{on }B_s(\tilde x_i).
\]
 
\begin{thm2}[Singular fibration]
Fix $n\in\mathbb{N}$, $\rho,D>0$, and $0<\delta<\delta(n)$. Let $(X_i,p_i,d_i,\mathcal{H}^n)$ be a sequence of $\mathrm{RCD}(-(n-1),n)$ spaces satisfying the $(\rho,\delta)$ local rewinding Reifenberg condition and $\operatorname{diam}(X_i)\le D$. Suppose that
\[
\left(
X_i,p_i,d_i,
\frac{\mathcal{H}^n}{\mathcal{H}^n(B_1(p_i))}
\right)
\xrightarrow{\mathrm{mGH}}
(X,p,d,\mathfrak{m}).
\]
Then for all sufficiently large $i$, there exists a singular fibration $\bar{F}_i:X_i\to X$ that is an $\epsilon_i$-Gromov--Hausdorff approximation and whose fibers are homeomorphic to infranilmanifolds.
\end{thm2}


We first explain the local construction when the covering
spaces converge to Euclidean space. The general construction
applies the same strategy to a small neighborhood close to an equivariant tangent cone.

\begin{theorem}[Theorem~\ref{thm:R^n singular fiber bundle}]
\label{thm:singular}
Let $(Y_i,p_i,d_i,\mathcal{H}^n)$ be a sequence of $\mathrm{RCD}(-\epsilon_i,n)$ spaces, where $\epsilon_i\to0$, and let $G_i\le\mathrm{Isom}(Y_i)$ be closed discrete subgroups acting freely on $Y_i$. Suppose that
\[
\begin{CD}
(Y_i,p_i,G_i,\mathcal{H}^n)
@>{\mathrm{eGH}}>>
(\mathbb{R}^n,p,G,\mathcal{H}^n)\\
@V{\pi}VV @VV{\pi}V\\
\left(
Y_i/G_i,\bar p_i,
\mathcal{H}^n/\mathcal{H}^n(B_1(\bar p_i))
\right)
@>{\mathrm{mGH}}>>
(\mathbb{R}^n/G,\bar p,\mathfrak{m}),
\end{CD}
\]
where $p=0^n$, $\bar p_i=\pi(p_i)$, and $\bar p=\pi(p)$. Assume further that
\[
G=\operatorname{trans}(G)\rtimes G_p,
\qquad
\operatorname{trans}(G)\cong\mathbb{R}^k,
\]
where $\operatorname{trans}(G)$ is the subgroup of pure translations in $G$ and $G_p$ is the stabilizer of $p$.
Then, for every fixed $R>0$ and all sufficiently large $i$, there exists an orbit-preserving homeomorphism
\[
F_i:(B_R(p_i),G_i)\longrightarrow(B_R(p),G)
\]
onto its image. Moreover, $F_i$ induces a quotient map
\[
\bar{F}_i:B_{0.99R}(\bar p_i)\longrightarrow B_{0.99R}(\bar p),
\]
and $\bar{F}_i^{-1}(x)$ is homeomorphic to an infranilmanifold for every $x\in B_{0.98R}(\bar p)$.
\end{theorem}

We sketch the proof of Theorem~\ref{thm:singular}; Section~\ref{example B} gives a concrete example. For simplicity, assume that $G$ is connected. The structure theory of approximate groups implies that $G=G_0$ is nilpotent. The stabilizer $G_p=T$ is a torus and acts trivially on the translation factor. Thus $T\le O(n-k)$ and $G=\operatorname{trans}(G)\times T$.
For every fixed $\epsilon>0$, the set $G(p,\epsilon)$ defines $G$.
A gluing construction allows us to discard distant generators and
relations while preserving the geometry and group actions on large
balls about the basepoints. We may therefore assume that $G_i$ is
defined by $G_i(p_i,\epsilon)$.

Using nilprogressions and Malcev completion, followed by a quotient
by central periods, we construct a connected nilpotent Lie group
$N_i$ equipped with a left-invariant metric and a natural embedding
of $G_i$ as a discrete subgroup. There is a diffeomorphism $\varphi_i:N_i\to G$ realizing
$C^\infty$ convergence on bounded sets.
If $T_i$ denotes the maximal torus of $N_i$, then
$\varphi_i|_{T_i}:T_i\to T$ is an isomorphism and
$\varphi_i(gt)=\varphi_i(g)\varphi_i(t)$
for every $g\in N_i$ and $t\in T_i$.

Recall that $\mathbb{R}^n=G\times_T\mathbb{R}^{n-k}$,
where the associated bundle is defined by the equivalence relation $(gt,v)\sim(g,tv)$ for $g\in G$, $t\in T$, and $v\in\mathbb R^{n-k}$.
We construct the nilpotent model space
\[
M_i=N_i\times_{T_i}\mathbb{R}^{n-k}
\]
using the relation $(gt,v)\sim(g,\varphi_i(t)v)$ for $g\in N_i$, $t\in T_i$, and $v\in\mathbb R^{n-k}$.
Equipped with a suitable connection metric, $M_i$ is $C^\infty$-close to $\mathbb{R}^n$ on bounded sets and carries an isometric $N_i$-action. The torus compatibility of $\varphi_i$ ensures that
\[
J_i:M_i\longrightarrow\mathbb{R}^n,
\qquad
J_i[g,v]=[\varphi_i(g),v],
\]
is a well-defined diffeomorphism. It induces a bijection between the $N_i$-orbits in $M_i$ and the $G$-orbits in $\mathbb{R}^n$.

By construction, $G_i\le N_i$, and $(Y_i,p_i)$ and $(M_i,p_i')$ are close in the $G_i$-equivariant Gromov--Hausdorff sense, where $p_i'=[e_{N_i},0^{n-k}]$.
For every fixed $R>0$ and all sufficiently large $i$, the gluing and averaging arguments produce a $G_i$-equivariant and $\epsilon_i$-injective map
\[
f_i:B_R(G_ip_i)\longrightarrow B_R(G_ip_i')\subset M_i
\]
that is almost splitting on each $R$-ball centered on a point in $G_ip_i$. In particular, $f_i$ is a $G_i$-equivariant homeomorphism onto its image. 

Thus $F_i:=J_i\circ f_i$ is orbit-preserving and descends
to the desired singular fibration $\bar F_i$. For $y\in B_{0.98R}(p_i')$, the fiber over $\pi(J_i(y))$ is homeomorphic to $N_i(y)/G_i$, where $N_i(y)$ denotes the $N_i$-orbit of $y$. The stabilizer $(N_i)_y$ lies in the central torus $T_i$, so $N_i(y)$ can be identified with the connected nilpotent Lie group $N_i/(N_i)_y$. Its compact quotient by the image of $G_i$ is a nilmanifold. When $G$ is disconnected, the additional components introduce affine actions with finite holonomy, and the fibers are, in general, infranilmanifolds.

The nilpotent models smoothly approximate the limiting geometry and group action while carrying exact actions of the approximating discrete groups. This provides the setting required for equivariant gluing.

Our local construction also yields nilpotent collapsing structures
on the collapsed regions of possibly noncompact $\mathrm{RCD}$
spaces satisfying the local rewinding Reifenberg condition.
For an $\mathrm{RCD}(-(n-1),n)$ space
$(X,p,d,\mathcal H^n)$, we call $x\in X$ $\epsilon$-collapsed if
$\mathcal H^n(B_1(x))<\epsilon$.

The precise formulation appears in Theorem~\ref{Theorem C}.
As in the classical theory of \cite{CFG}, the collapsed region
admits an atlas of local singular fibrations with compatible fiber
partitions: locally on an overlap, one partition refines the other.
Our formulation uses equivariant maps from local covering spaces
to nilpotent models.

\begin{thm3}[Theorem~\ref{Theorem C}]
Let $(X,p,d,\mathcal{H}^n)$ be an $\mathrm{RCD}(-(n-1),n)$ space satisfying the $(\rho,\delta)$ local rewinding Reifenberg condition for some $\rho>0$ and $0<\delta<\delta(n)$. Then there exists $\epsilon=\epsilon(n,\rho)>0$ such that the $\epsilon$-collapsed region of $X$ admits a nilpotent collapsing structure.
\end{thm3}

In Theorem~\ref{thm:bounded sec}, we give an independent proof that two-sided sectional curvature bounds on complete manifolds yield a uniform positive lower bound on the injectivity radius at the basepoints of local universal covers. This estimate also follows from the main theorem of \cite{CFG}. It implies the local rewinding Reifenberg condition. Our methods therefore provide a new construction of the singular fibrations and nilpotent collapsing atlases of classical collapsing theory, together with explicit descriptions of their infranilmanifold fibers.

The paper is organized as follows. Section~2 collects the preliminaries. Section~3 presents examples illustrating Theorems~A and~B. Section~4 proves Theorem~A by combining the structure of approximate groups at large scales with equivariant gluing of almost splitting maps.

Section~5 constructs nilpotent models and local singular fibrations when the equivariant limit is a metric cone and the limit group splits into its translation and isotropy parts. Section~6 studies the group structure of equivariant tangent cones and constructs the local singular fibrations needed
for Theorem~B. Section~7 glues these local fibrations to prove Theorem~B.

Section~8 constructs compatible local nilpotent models on the collapsed region, proving Theorem~C. The appendix proves Theorem~\ref{thm:large nil} on nilprogressions at large scales.
   
{\bf Acknowledgments.}
We thank Xiaochun Rong for helpful discussions and for informing us that Cheng, Xu, and Zhang had proved the Riemannian case of Theorem~B using Ricci flow smoothing; we have not had access to their manuscript.
\section{Preliminaries}
\subsection{Equivariant Gromov--Hausdorff convergence}

We recall the notion of equivariant Gromov--Hausdorff convergence introduced by Fukaya and Yamaguchi \cite{Fukaya1986,FukayaYamaguchi1992}.

Let $(X,p)$ and $(Y,q)$ be pointed metric spaces, and let $H\le\mathrm{Isom}(X)$ and $K\le\mathrm{Isom}(Y)$ be closed subgroups. For $r>0$, define
\[
H(p,r)=\{h\in H:d_X(hp,p)\le r\},
\qquad
K(q,r)=\{k\in K:d_Y(kq,q)\le r\}.
\]

For $\epsilon>0$, a pointed $\epsilon$-equivariant Gromov--Hausdorff approximation, abbreviated as an $\epsilon$-eGHA, is a triple of maps $(f,\phi,\psi)$,
\[
f:B_{1/\epsilon}(p)\to B_{1/\epsilon}(q), \ \phi:H(p,1/\epsilon)\to K(q,1/\epsilon), \ \psi:K(q,1/\epsilon)\to H(p,1/\epsilon),
\]
satisfying the following conditions:
\begin{enumerate}
\item
$d(f(p),q)<\epsilon$, the image $f(B_{1/\epsilon}(p))$ is $2\epsilon$-dense in $B_{1/\epsilon}(q)$, and for all $x_1,x_2\in B_{1/\epsilon}(p)$, $\left| d_Y(f(x_1),f(x_2))-d_X(x_1,x_2) \right| \le\epsilon$.

\item
For every $h\in H(p,1/\epsilon)$ and every $x$ such that
$x,hx\in B_{1/\epsilon}(p)$,
\[
d_Y\bigl(\phi(h)f(x),f(hx)\bigr)<\epsilon.
\]

\item
For every $k\in K(q,1/\epsilon)$ and every $x$ such that
$x,\psi(k)x\in B_{1/\epsilon}(p)$,
\[
d_Y\bigl(kf(x),f(\psi(k)x)\bigr)<\epsilon.
\]
\end{enumerate}

The equivariant Gromov--Hausdorff distance $d_{\mathrm{eGH}}\bigl((X,p,H),(Y,q,K)\bigr)$ is the infimum of all $\epsilon>0$ for which such an $\epsilon$-eGHA exists. We say that $(X_i,p_i,G_i)$ converges to $(X,p,G)$ in the pointed equivariant Gromov--Hausdorff sense if
\[
d_{\mathrm{eGH}}\bigl((X_i,p_i,G_i),(X,p,G)\bigr) \le \epsilon_i \to 0.
\]

Throughout the paper, $\epsilon_i>0$ denotes a sequence tending to zero as $i\to\infty$. Its value may change from one occurrence to another, and it may depend on parameters held fixed in the argument. We enlarge such sequences when necessary without changing notation.

To simplify the notation, when we write $f:B_r(x)\to B_r(y)$ as a Gromov--Hausdorff approximation, an almost splitting map, a homeomorphism, or a singular fibration, we allow boundary collars controlled by the relevant approximation error to be omitted. This convention also applies to orbit neighborhoods and quotient maps.

We will use the following equivariant precompactness theorem.
\begin{theorem}[\cite{Fukaya1986,FukayaYamaguchi1992}]
\label{def:eGH}
Let $(X_i,p_i)$ be a sequence of pointed proper metric spaces converging to a pointed proper metric space $(X,p)$ in the pointed Gromov--Hausdorff sense. For each $i$, let $G_i\le\mathrm{Isom}(X_i)$ be a closed subgroup. Then, after passing to a subsequence,
\[
(X_i,p_i,G_i)\xrightarrow{\mathrm{eGH}}(X,p,G),
\]
where $G\le\mathrm{Isom}(X)$ is a closed subgroup. Moreover,
\[
(X_i/G_i,\bar p_i)\xrightarrow{\mathrm{GH}}(X/G,\bar p),
\]
where $\bar p_i$ and $\bar p$ denote the images of $p_i$ and $p$ in the respective quotient spaces.
\end{theorem}

When the actions have uniformly bounded displacement and the limit group is a Lie group, the forward group maps can be chosen to be homomorphisms.
\begin{lemma}[\cite{MazurRongWang,Harvey2016}]
\label{lemma:homomorphic eGHA}
Suppose that
\[
(X_i,p_i,G_i)\xrightarrow{\mathrm{eGH}}(X,p,G),
\]
where the ambient metric spaces are proper. Assume that $G$ is a Lie group and that there exists $R>0$ such that $G_i=G_i(p_i,R)$ for every $i$.
Then, for all sufficiently large $i$, the group maps
$\phi_i:G_i\to G$ in the equivariant approximations can be chosen to be continuous homomorphisms.
\end{lemma}

We also use the following convergence result for quotient actions.
\begin{lemma}[\cite{Wang2023}]
\label{lem:quotient_eGH}
Let $(X_i,p_i)$ be pointed proper geodesic metric spaces, and let $G_i\le\mathrm{Isom}(X_i)$ be closed subgroups. Suppose that $G_i'\trianglelefteq G_i$ is a closed normal subgroup for each $i$ and that
\[
(X_i,p_i,G_i,G_i')
\xrightarrow{\mathrm{eGH}}
(X,p,G,G').
\]
Then $G'$ is a normal subgroup of $G$, and the induced quotient actions converge:
\[
\begin{CD}
(X_i,p_i,G_i,G_i')
@>{\mathrm{eGH}}>>
(X,p,G,G')\\
@V{\pi}VV @VV{\pi}V\\
(X_i/G_i',\bar p_i,G_i/G_i')
@>{\mathrm{eGH}}>>
(X/G',\bar p,G/G').
\end{CD}
\]
\end{lemma}

In Lemma~\ref{lem:quotient_eGH}, the quotient groups are regarded as groups acting on the quotient spaces, and these actions need not be effective. The induced quotient group maps are almost homomorphic on large domains. Passing to the images of the quotient groups in the respective isometry groups gives the corresponding convergence of effective actions.

Equivariant Gromov--Hausdorff convergence provides natural discrete approximations; see the next subsection for the definitions of discrete approximations and small subgroups.

\begin{lemma}\label{lem:eGH discrete appro}
Let $(X_i,p_i)$ and $(X,p)$ be pointed proper metric spaces.
Let $\Gamma_i\le\mathrm{Isom}(X_i)$ be discrete subgroups,
and let $\Gamma\le\mathrm{Isom}(X)$ be a closed subgroup.
Suppose that
\[
(X_i,p_i,\Gamma_i)
\xrightarrow{\mathrm{eGH}}
(X,p,\Gamma)
\]
and that $\Gamma$ is a Lie group. Then the local group maps $\varphi_i:\Gamma_i\to\Gamma$ in the eGHAs may be chosen with $\varphi_i(e_i)=e_\Gamma$, and these normalized maps form a discrete approximation of $\Gamma$.
\end{lemma}

\begin{lemma}\label{lem:quotient appro}
In the setting of Lemma~\ref{lem:eGH discrete appro}, let $K$ be a compact normal subgroup of $\Gamma$.
Then the maps
$\bar\varphi_i:=\pi\circ\varphi_i:\Gamma_i\to\Gamma/K$
form a discrete approximation of $\Gamma/K$, where
$\pi:\Gamma\to\Gamma/K$ is the quotient homomorphism.

Moreover, a sequence of subgroups $H_i\le\Gamma_i$ is small
with respect to $\bar\varphi_i$ if and only if both of the
following conditions hold:
\begin{enumerate}
\item There exists $R>0$ such that
$H_i=\langle H_i\cap\Gamma_i(p_i,R)\rangle$ for every $i$.
\item Every subsequential equivariant limit of the actions
of $H_i$ is contained in $K$.
\end{enumerate}
\end{lemma}

Finally, we specify our convention for metrics on groups. Suppose that
\[
(X_i,p_i,G_i)\xrightarrow{\mathrm{eGH}}(X,p,G).
\]
For $\epsilon>0$, we say that $g,h\in G_i(p_i,1/(10\epsilon))$ are $\epsilon$-close if
\[
d_i(g^{-1}hx,x)<\epsilon
\qquad
\text{for every }x\in B_{1/\epsilon}(p_i).
\]

On the limit group $G$, fix a compatible left-invariant metric $d_G$. For each $i$, let $d_{G_i}$ be a left-invariant metric on $G_i$. Fix a large $R$ and assume that for all sequences $g_i,h_i\in G_i(p_i,R)$ converging to $g,h\in G(p,R)$, we have $d_{G_i}(g_i,h_i)\to d_G(g,h)$. Then there is a function $\Phi(\epsilon)\to0$ as $\epsilon\to0$ such that, for all sufficiently large $i$, the inequality $d_{G_i}(g_i,h_i)<\epsilon$ for $g_i,h_i\in G_i(p_i,R)$ implies that $g_i$ and $h_i$ are $\Phi(\epsilon)$-close in the preceding sense. Once these metrics have been fixed, we call \(g_i\) and \(h_i\) \(\epsilon\)-close if \(d_{G_i}(g_i,h_i)<\epsilon\). 

\subsection{Local groups and discrete approximations}\label{subsection:local group}
\begin{defn}
A symmetric local group is a topological space $G$ equipped with an identity element $e\in G$, an inversion map and a partially defined multiplication
\[
g\longmapsto g^{-1}:G\longrightarrow G, \ \cdot:\Omega\longrightarrow G, \ \Omega\subset G\times G,
\]
satisfying the following conditions:
\begin{enumerate}
\item
$\Omega$ is open and contains $(G\times\{e\})\cup(\{e\}\times G)$.

\item
The inverse map and the multiplication map are continuous.

\item
Whenever both $(gh)k$ and $g(hk)$ are well-defined, they are equal.

\item
For any $g\in G$, we have $eg=ge=g$.

\item
For any $g\in G$, $gg^{-1}$ and $g^{-1}g$ are well-defined and equal to $e$.
\end{enumerate}
\end{defn}

If $\Omega=G\times G$, then $G$ is a topological group.
\begin{defn}
Let $G$ be a local group and let $g_1,\dots,g_m\in G$.
We say that the product $g_1\cdots g_m$ is well-defined
if there exists a family $\{g_{[i,j]}:1\le i\le j\le m\}\subset G$
such that $g_{[j,j]}=g_j$ for any $1\le j\le m$,
and whenever $1\le i\le j<k\le m$, the product
$g_{[i,j]}g_{[j+1,k]}$ is well-defined and satisfies
\[
g_{[i,j]}g_{[j+1,k]}=g_{[i,k]}.
\]
In this case, we write $g_1\cdots g_m=g_{[1,m]}$.
\end{defn}

For subsets $A_1,\dots,A_m\subset G$, we say that
$A_1\cdots A_m$ is well-defined if $g_1\cdots g_m$ is
well-defined for every choice of $g_j\in A_j$.

\begin{defn}
A subset $A\subset G$ is called a multiplicative set
if it is symmetric and $A^{200}$ is well-defined.
\end{defn}

\begin{defn}
Let $A$ be a finite symmetric subset of a multiplicative set,
and let $C\in\mathbb{N}$. We call $A$ a
$C$-approximate group if $A^2$ can be covered by at
most $C$ left translates of $A$.
\end{defn}

\begin{defn}
A $C$-approximate group $A$ is called a
strong $C$-approximate group if there exists a symmetric
subset $S\subset A$ satisfying the following conditions:
\begin{enumerate}
\item
The product $\bigl(\{asa^{-1}:s\in S,\ a\in A^4\}\bigr)^{10^3C^3}$ is well-defined and contained in $A$.

\item
If $g,g^2,\dots,g^{1000}$ are well-defined and lie in $A^{100}$, then $g\in A$.

\item
If the powers $g,g^2,\dots,g^{10^6C^3}$ are well-defined and
belong to $A$, then $g\in S$.
\end{enumerate}
\end{defn}

For a multiplicative set $A\subset G$ and $g\in G$, define
the escape norm by
\[
\|g\|_A
:=
\inf\left\{
\frac{1}{m+1}:
m\in\mathbb{Z}_{\ge0},\
g^j\in A\text{ for every }0\le j\le m
\right\}.
\]

\begin{theorem}[Escape norm estimates, \cite{BGT}]
\label{thm:escape}
For every $C>0$, there exists $C'>0$ such that every strong
$C$-approximate group $A$ satisfies the following estimates:
\begin{enumerate}
\item
For $g_1,\dots,g_k\in A^{10}$, whenever $g_1\cdots g_k$
is well-defined,
\[
\|g_1\cdots g_k\|_A
\le C'\sum_{j=1}^k\|g_j\|_A.
\]

\item
For $g_1,g_2\in A^{10}$, $\|g_2g_1g_2^{-1}\|_A\le 10^3\|g_1\|_A$.

\item
For $g_1,g_2\in A^{10}$, $\|[g_1,g_2]\|_A\le C'\|g_1\|_A\|g_2\|_A$.
\end{enumerate}
\end{theorem}

To study discrete local groups, we use the following notion of discrete approximation; see \cite{BGT,Zamora2020}.

\begin{defn}\label{defn:discrete appoximation}
Let $\Gamma$ be a Lie group equipped with a proper compatible left-invariant
metric $d_\Gamma$, and let $\Gamma_i$ be a sequence of discrete
local groups. Maps $\varphi_i:\Gamma_i\to\Gamma$, normalized by
$\varphi_i(e_i)=e_\Gamma$, form a discrete approximation
of $\Gamma$ if there exist sequences $R_i\to\infty$ and
$\epsilon_i\to0$ such that, for all sufficiently large $i$,
the following conditions hold. Set $U_i:=\varphi_i^{-1}(B_{R_i}(e_\Gamma))$.
\begin{enumerate}
\item
$U_i$ is a finite multiplicative set.

\item
$\varphi_i(U_i)$ is $\epsilon_i$-dense in
$B_{R_i}(e_\Gamma)$.

\item
For every $g_1,g_2\in U_i$, the product $g_1g_2^{-1}$
is well-defined, and
\[
d_\Gamma\bigl(
\varphi_i(g_1g_2^{-1}),
\varphi_i(g_1)\varphi_i(g_2)^{-1}
\bigr)\le\epsilon_i.
\]
\end{enumerate}
\end{defn}

\begin{remark}
Only the local structure on $U_i$ is relevant to our arguments.
We may assume that $\Gamma_i=U_i^{200}$.
\end{remark}

The following consequence of \cite[Proposition~7.3]{BGT}
provides strong approximate groups from discrete approximations.

\begin{theorem}\label{thm:strong}
Suppose that the maps $\varphi_i:\Gamma_i\to\Gamma$ form a discrete
approximation. There exists $C>0$ with the following property.
For every $\delta'>0$, one can choose $0<\delta<\delta'$ such that,
for all sufficiently large $i$, there is a strong
$C$-approximate group $A_i$ satisfying $\varphi_i^{-1}(B_\delta(e_\Gamma)) \subset A_i \subset \varphi_i^{-1}(B_{\delta'}(e_\Gamma))$.
\end{theorem}

\begin{defn}\label{defn:small subgroups}
Let $\varphi_i:\Gamma_i\to\Gamma$ form a discrete approximation.
A sequence of genuine subgroups $H_i\le\Gamma_i$ with $H_i\subset U_i$ is called
small with respect to $\varphi_i$ if there exists
a sequence $\epsilon_i\to0$ such that $\varphi_i(H_i)\subset B_{\epsilon_i}(e_\Gamma)$.
We say that $(\Gamma_i)$ has no small subgroups with
respect to $\varphi_i$ if every such sequence $(H_i)$
is eventually trivial.
\end{defn}

Choose $\delta'>0$ sufficiently small that $\overline{B_{\delta'}(e_\Gamma)}$ contains no nontrivial subgroup and let $A_i$ be the strong approximate groups supplied by Theorem~\ref{thm:strong}. Then $(\Gamma_i)$ has no small subgroups if and only if, for all sufficiently large $i$, $\|g\|_{A_i}>0$ for every nonidentity $g\in\Gamma_i$.
By Theorem~\ref{thm:escape}, the set
\[
H_i:=\{g\in A_i:\|g\|_{A_i}=0\}
\]
is a finite subgroup normalized by $A_i^{10}$. Every Hausdorff limit
of $\varphi_i(H_i)$ is a subgroup contained in $\overline{B_{\delta'}(e_\Gamma)}$,
so $H_i$ is small. Conversely, every small subgroup is eventually
contained in $A_i$ and hence in $H_i$.

Since $\varphi_i$ is almost homomorphic on $U_i$, $H_i$ is normalized by $U_i$. By our convention that $\Gamma_i=U_i^{200}$, $H_i$ is normal in $\Gamma_i$ and
$\bar\varphi_i:\Gamma_i/H_i\longrightarrow\Gamma$
is a discrete approximation with no small subgroups.

We next recall the structure of discrete approximate groups.
Let $G$ be a discrete local group. For elements $u_1,\dots,u_k\in G$
and parameters $C_1,\dots,C_k>0$ for which all the products below are well-defined, set
\[
P(u_1,\dots,u_k;C_1,\dots,C_k)
:=
\left\{
u_1^{n_1}\cdots u_k^{n_k}:
n_j\in\mathbb{Z},\ |n_j|\le C_j
\text{ for }1\le j\le k
\right\}.
\]
The nilprogressions appearing below satisfy additional
triangular commutator relations.

\begin{theorem}[\cite{BGT,Zamora2020}]
\label{thm:nil}
Let $\Gamma$ be a $k$-dimensional Lie group equipped with
a left-invariant metric, and let
$\varphi_i:\Gamma_i\to\Gamma$ be a discrete approximation
by discrete local groups. Let $H_i$ be the small subgroups
constructed above, and let $\bar\varphi_i:\Gamma_i/H_i\longrightarrow\Gamma$
denote the corresponding quotient approximations.

Then there exist constants $C,\delta>0$ such that, for all
sufficiently large $i$, there are elements
$u_{i,1},\dots,u_{i,k}\in\Gamma_i/H_i$ and parameters
$C_{i,1},\dots,C_{i,k}>0$ for which the nilprogression
\[
P_i=P(u_{i,1},\dots,u_{i,k};C_{i,1},\dots,C_{i,k})
\]
is well-defined and has the following properties:
\begin{enumerate}
\item
For every $1\le j\le k$, $C_{i,j}\longrightarrow\infty$ as $i\to\infty$.

\item
For every $1\le l<j\le k$,
\[
[u_{i,l},u_{i,j}]
\in
P\left(
u_{i,j+1},\dots,u_{i,k};
\frac{C C_{i,j+1}}{C_{i,l}C_{i,j}},
\dots,
\frac{C C_{i,k}}{C_{i,l}C_{i,j}}
\right).
\]

\item
We have $\bar\varphi_i^{-1}(B_\delta(e_\Gamma))\subset P_i$ and $\bar\varphi_i(P_i)\subset B_C(e_\Gamma)$.

\item
The ordered representation is unique in these ranges. If integers
$n_j,n_j'$ satisfy $|n_j|,|n_j'|\le C_{i,j}$ for every $j$ and
$u_{i,1}^{n_1}\cdots u_{i,k}^{n_k}=u_{i,1}^{n_1'}\cdots u_{i,k}^{n_k'}$,
then $n_j=n_j'$ for every $1\le j\le k$.
\end{enumerate}
\end{theorem}

In \cite{Zamora2020}, the uniqueness assertion in Theorem~\ref{thm:nil} holds only on a subset of the nilprogression, called its grid part. To simplify the notation, we use the grid parts of these nilprogressions, on which the ordered representations are unique, and retain the notation \(P_i\). For all sufficiently large \(i\), the groupification of $P_i$ is a finitely generated torsion-free nilpotent group of rank $k$, embedded as a lattice in a $k$-dimensional connected, simply connected nilpotent Lie group $N_i$; see \cite[Theorem~2.53]{Zamora2020}.

To establish convergence of the Lie algebras of $N_i$ to
that of $\Gamma$, we first recall strong Malcev bases and
their associated coordinates.

\begin{defn}
Let $\mathfrak{n}$ be a nilpotent Lie algebra.
An ordered basis $(v_1,\dots,v_k)$ of $\mathfrak{n}$ is
called a strong Malcev basis if, for every
$1\le j\le k$, the subspace
\[
V_j:=\operatorname{span}\{v_{j+1},\dots,v_k\}
\]
is an ideal of $\mathfrak{n}$ and $v_j+V_j$ lies in the
center of $\mathfrak{n}/V_j$. Here $V_k=\{0\}$.
\end{defn}

\begin{theorem}\label{thm:nil str}
Let $N$ be a $k$-dimensional connected, simply connected
nilpotent Lie group, and let $\mathfrak{n}$ be its Lie
algebra. Suppose that $(v_1,\dots,v_k)$ is a strong Malcev
basis of $\mathfrak{n}$. Then:
\begin{enumerate}
\item
The exponential map $\exp_N:\mathfrak{n}\longrightarrow N$ is a diffeomorphism.

\item
The map $\psi:\mathbb{R}^k\longrightarrow N$, $\psi(x_1,\dots,x_k)=\exp_N(x_1v_1)\cdots\exp_N(x_kv_k)$, is a diffeomorphism.

\item
After identifying $\mathfrak{n}$ with $\mathbb{R}^k$
through the given basis, the mutually inverse maps
$\exp_N^{-1}\circ\psi:\mathbb{R}^k\longrightarrow\mathbb{R}^k$ and $\psi^{-1}\circ\exp_N:\mathbb{R}^k\longrightarrow\mathbb{R}^k$
are polynomial maps whose degrees are bounded solely in
terms of $k$.
\end{enumerate}
\end{theorem}

Consider the nilprogressions in Theorem~\ref{thm:nil},
\[
P_i=P(u_{i,1},\dots,u_{i,k};C_{i,1},\dots,C_{i,k}).
\]
For each $1\le j\le k$, define
$g_{i,j}:=u_{i,j}^{\lfloor C_{i,j}/C\rfloor}$ and
$v_{i,j}:=\log_{N_i}(g_{i,j})\in\mathfrak n_i$.
Then $(v_{i,1},\dots,v_{i,k})$ is a strong Malcev basis of $\mathfrak n_i$. By the choice of $C_{i,j}$, there is $\epsilon'>0$ such that $\bar\varphi_i(g_{i,j})\notin B_{\epsilon'}(e_\Gamma)$ for every $j$ and all sufficiently large $i$. After passing to a subsequence, we may assume that
$\bar\varphi_i(g_{i,j})\longrightarrow g_j\in\Gamma_0\setminus\{e_\Gamma\}$ for any
$1\le j\le k$.

Let $\mathfrak g$ be the Lie algebra of $\Gamma$. The corresponding limiting one-parameter segments determine vectors $v_j\in\mathfrak g$ with $\exp_\Gamma(v_j)=g_j$. By \cite{BGT,Zamora2020}, $(v_1,\dots,v_k)$ is a strong Malcev basis of $\mathfrak g$.

Let $\widetilde\Gamma_0$ be the universal covering group of $\Gamma_0$. Use these bases to define the ordered exponential charts $\psi_i:\mathbb R^k\to N_i$ and $\psi:\mathbb R^k\to\widetilde\Gamma_0$ as in Theorem~\ref{thm:nil str}. In these coordinates, the multiplication laws are
\[
Q_i(x,y):=\psi_i^{-1}(\psi_i(x)\psi_i(y)),
\qquad
Q(x,y):=\psi^{-1}(\psi(x)\psi(y)).
\]
Their derivatives at the origin determine the corresponding Lie algebra structures.

\begin{theorem}[Lie algebra structure convergence,
{\cite[Proposition~8.2]{Zamora2020}}]
\label{thm:Q conv}
For all sufficiently large $i$, the maps $Q_i$ and $Q$ are polynomial maps of degree at most $d(k)$. Moreover, the coefficients of $Q_i$ converge to the corresponding coefficients of $Q$.
\end{theorem}

Theorem~\ref{thm:nil} describes approximations to a sufficiently
small neighborhood of $e_\Gamma$ containing no nontrivial
subgroup. We also need the following result concerning
approximate groups on large scales.

\begin{theorem}\label{thm:large nil}
Let $\Gamma$ be a compactly generated nilpotent Lie group
equipped with a left-invariant metric. Assume that $\Gamma$
has no nontrivial compact subgroup, and write $\dim\Gamma=k$ and $\operatorname{rank}(\Gamma/\Gamma_0)=l$.
Let $\varphi_i:\Gamma_i\to\Gamma$ be a discrete approximation,
and let $\bar\varphi_i:\Gamma_i/H_i\longrightarrow\Gamma$
be the quotient approximations from Theorem~\ref{thm:nil}.
Let $\pi:\Gamma\to\Gamma/\Gamma_0$ be the quotient map.
Choose elements $w_{i,1},\dots,w_{i,l}\in\Gamma_i/H_i$
such that
\[
\bar\varphi_i(w_{i,s})\longrightarrow g_s\in\Gamma
\qquad (1\le s\le l),
\]
where $(\pi(g_1),\dots,\pi(g_l))$ is a Malcev basis of
$\Gamma/\Gamma_0$.

Then there exist elements $u_{i,1}',\dots,u_{i,k}'$ and
parameters $C_{i,1}',\dots,C_{i,k}'>0$ such that
\[
P(u_{i,1}',\dots,u_{i,k}';C_{i,1}',\dots,C_{i,k}')
\]
satisfies all the conclusions of Theorem~\ref{thm:nil}.
Moreover, there exist $C>0$ and a function
$\Lambda:(0,\infty)\to(0,\infty)$ with $\lim_{r\to\infty}\Lambda(r)=\infty$ such that, for every fixed $r>0$ and all sufficiently large
$i$ (depending on $r$), the progression
\[
P_i'(r):=
P\left(
w_{i,1},\dots,w_{i,l},u_{i,1}',\dots,u_{i,k}';
\underbrace{r,\dots,r}_{l\text{ terms}},
rC_{i,1}',\dots,rC_{i,k}'
\right)
\]
has the following properties:
\begin{enumerate}
\item
We have $\bar\varphi_i^{-1}(B_{\Lambda(r)}(e_\Gamma))
\subset P_i'(r)$ and $\bar\varphi_i(P_i'(r))
\subset B_{(r+1)C}(e_\Gamma)$.

\item
The ordered representation is unique within the prescribed ranges.
More precisely, if integers $m_s,m_s',n_j,n_j'$ satisfy
\[
|m_s|,|m_s'|\le r,
\qquad |n_j|,|n_j'|\le rC_{i,j}',
\]
\[
w_{i,1}^{m_1}\cdots w_{i,l}^{m_l}
(u_{i,1}')^{n_1}\cdots(u_{i,k}')^{n_k}
=
w_{i,1}^{m_1'}\cdots w_{i,l}^{m_l'}
(u_{i,1}')^{n_1'}\cdots(u_{i,k}')^{n_k'},
\]
then $m_s=m_s'$ and $n_j=n_j'$ for every $s,j$.

\item
For every $1\le s<j\le l$,
\begin{equation}
[w_{i,s},w_{i,j}]
\in
P\left(
w_{i,j+1},\dots,w_{i,l},u_{i,1}',\dots,u_{i,k}';
C,\dots,C,CC_{i,1}',\dots,CC_{i,k}'
\right).
\end{equation}

\item
For every $1\le s\le l$ and $1\le j\le k$,
\begin{equation}\label{e0}
[w_{i,s},u_{i,j}']
\in
P\left(
u_{i,j+1}',\dots,u_{i,k}';
C_{i,j+1}',\dots,C_{i,k}'
\right).
\end{equation}
\end{enumerate}
\end{theorem}
A complete proof of Theorem~\ref{thm:large nil} is given in Appendix~A.

For every sufficiently large fixed $r$ and all sufficiently large $i$, the progression $P_i'(r)$ in
Theorem~\ref{thm:large nil} can also be identified with a defining subset of a lattice in a connected, simply connected nilpotent Lie group of dimension $k+l$. The coefficient convergence in Theorem~\ref{thm:Q conv} extends to these enlarged progressions, with the limiting group given by the connected, simply connected nilpotent completion of $\Gamma$.

\subsection{Geometry of RCD spaces}
We recall the geometric facts about $\mathrm{RCD}(K,N)$ spaces used below. These spaces have Lie isometry groups \cite{Sosa2018,GuSan2019} and satisfy the Margulis lemma \cite{DSZZ2023}.
\begin{theorem}[Margulis lemma, \cite{DSZZ2023}]
\label{thm:Margulis}
For every $K\in\mathbb R$ and $N\ge1$, there exist constants $C=C(K,N)>0$ and $\epsilon_0=\epsilon_0(K,N)>0$ with the following property. Let $(X,d,\mathfrak m,p)$ be an $\mathrm{RCD}(K,N)$ space of rectifiable dimension $k$. Then the image of the natural homomorphism
\[
\pi_1(B_{\epsilon_0}(p),p)
\longrightarrow
\pi_1(X,p)
\]
contains a normal nilpotent subgroup of index at most $C$. Moreover, this subgroup admits a nilpotent basis of length at most $k$. 
\end{theorem}

Under the local rewinding Reifenberg condition, the following local version holds; see \cite[Lemma 2.6]{Wang2024} and compare to \cite{KapovitchWilking2011}.

\begin{lemma}\label{thm:Margulis2}
For every $n\ge1$, there exist constants $\delta(n),\epsilon(n),C(n)>0$ with the following property. Let $0<\rho\le1$ and $0<\delta<\delta(n)$, and let $(X,d,\mathcal H^n,p)$ be an $\mathrm{RCD}(-(n-1),n)$ space satisfying the $(\rho,\delta)$ local rewinding Reifenberg condition at $p$. Then, for every $0<\epsilon<\epsilon(n)$, the image of the natural homomorphism
\[
\pi_1(B_{\epsilon\rho}(p),p)
\longrightarrow
\pi_1(B_\rho(p),p)
\]
contains a normal nilpotent subgroup of index at most $C(n)$. Moreover, this subgroup admits a nilpotent basis of length at most $n$.
\end{lemma}

An $\mathrm{RCD}(K,N)$ space $(X,d,\mathfrak m)$ is noncollapsed if $N \in \mathbb N$ and $\mathfrak m=\mathcal H^N$. Throughout, $0^k$ denotes the origin of $\mathbb R^k$.

\begin{defn}
Let $(X,d,\mathfrak m)$ be an $\mathrm{RCD}(K,N)$ space. For $\epsilon,r>0$ and $k\in\mathbb N$, define
\[
\mathcal R_{\epsilon,r}^k(X)
:=
\left\{
x\in X:
d_{\mathrm{GH}}\bigl(B_s(x),B_s(0^k)\bigr)<\epsilon s
\text{ for every }0<s<2r
\right\},
\]
and set
\[
\mathcal R_\epsilon^k(X)
:=
\bigcup_{r>0}\mathcal R_{\epsilon,r}^k(X),
\qquad
\mathcal R^k(X)
:=
\bigcap_{\epsilon>0}\mathcal R_\epsilon^k(X).
\]
Points in $\mathcal R_\epsilon^k(X)$ are called $(k,\epsilon)$-regular, and points in $\mathcal R^k(X)$ are called $k$-regular.
\end{defn}

By \cite{BS2023}, every $\mathrm{RCD}(K,N)$ space admits an integer $k\le N$ such that $\mathfrak m\bigl(X\setminus\mathcal R^k(X)\bigr)=0$.
We call $k$ the rectifiable dimension of $(X,d,\mathfrak m)$.

The notion of $\delta$-splitting maps originates in Cheeger--Colding--Naber theory
\cite{CheegerColding1997,CheegerColding2000a,CheegerNaber2015},
and their $\mathrm{RCD}$ extensions appear in \cite{BNS}.
\begin{defn}\label{def:analytic splitting}
Let $(X,d,\mathfrak m)$ be an $\mathrm{RCD}(-1,N)$ space, $p\in X$, $r,\delta>0$, and $1\le k\le N$. A map $u=(u^1,\dots,u^k):B_r(p)\to\mathbb R^k$ is a $(\delta,k)$-splitting map if its coordinates lie in the domain of the local Laplacian on $B_r(p)$ and satisfy
\[
\operatorname{Lip}(u)\le C(N),
\qquad
\|\Delta u^a\|_{L^\infty(B_r(p))}<\delta
\quad(1\le a\le k),
\]
\[
\fint_{B_r(p)}
\bigl|\langle\nabla u^a,\nabla u^b\rangle-\delta_{ab}\bigr|
\,d\mathfrak m<\delta
\quad(1\le a,b\le k).
\]
\end{defn}

If every coordinate is harmonic, we call $u$ a harmonic $(\delta,k)$-splitting map. The pointed notation
$u:(B_r(p),p)\to(\mathbb R^k,0^k)$ additionally requires $u(p)=0^k$.

We next recall the volume convergence theorem for noncollapsed $\mathrm{RCD}(K,n)$ spaces.
\begin{theorem}[Volume convergence, \cite{DeGigli2019}]
Let $(X_i,d_i,\mathcal H^n,p_i)$ be a sequence of noncollapsed $\mathrm{RCD}(K,n)$ spaces converging to $(X,d,\mathfrak m,p)$ in the pointed measured Gromov--Hausdorff sense. If $\mathcal H^n(B_1(p_i))\ge v>0$ for every $i$, then $\mathfrak m=\mathcal H^n$ and
\[
\mathcal H^n(B_r(p_i))
\longrightarrow
\mathcal H^n(B_r(p))
\qquad\text{for every }r>0.
\]
\end{theorem}

We next state the transformation theorem \cite{BNS,CheegerNaberJiang2021,CheegerNaber2015,HondaPeng2024,HuangHuang2024}. We shall use the bounded-Laplacian formulation in \cite{Huang2024}; see also \cite[Theorem~2.6]{HHWZ}. Here a ball $B_s(x)$ is $(n,\eta)$-symmetric if
$d_{\mathrm{GH}}(B_s(x),B_s(0^n))\le\eta s$.

\begin{theorem}[Transformation theorem]\label{transformation}
Let $n\in\mathbb N$. For every $\delta>0$, there exist constants
$C=C(n,\delta)>0$ and $\epsilon=\epsilon(n,\delta)>0$ with the
following property. Let $(X,d,\mathcal H^n)$ be an
$\mathrm{RCD}(-\epsilon^2(n-1),n)$ space, let $x\in X$, and fix $0<r_0\le1$.
Suppose that $B_s(x)$ is $(n,\epsilon^2)$-symmetric for every
$r_0\le s\le1$ and that
$u:B_2(x)\to\mathbb R^n$ is an $(\epsilon,n)$-splitting map.
Then, for each $r_0\le s\le1$, there is an invertible lower
triangular matrix $T_s\in\mathbb R^{n \times n}$ with positive diagonal entries such that:
\begin{enumerate}
\item
The map $s^{-1}T_su: \bigl(B_s(x),s^{-1}d\bigr)\longrightarrow\mathbb R^n$ is a $(\delta,n)$-splitting map.

\item
For every $1\le a,b\le n$, $\fint_{B_s(x)} \bigl\langle\nabla(T_su)^a,\nabla(T_su)^b\bigr\rangle \,d\mathcal H^n =\delta_{ab}$.

\item
For every $r_0\le s\le\tfrac12$, $\|T_sT_{2s}^{-1}-\mathrm{I}\|\le\delta$.

\item
For every $r_0\le s\le1$, $\|T_s\|\le Cs^{-\delta}$ and $\|T_s^{-1}\|\le C$.
\end{enumerate}
\end{theorem}

One application of the transformation theorem is the construction
of bi-H\"older coordinate maps.
$\Phi(t\mid\cdot)$ always denotes a nonnegative function tending to zero as $t\to0$, with the other indicated parameters held fixed. 
The transformation theorem yields the following canonical Reifenberg theorem; see \cite{BNS,CheegerNaberJiang2021,HondaPeng2024}.
\begin{theorem}\label{thm:Rei}
Let $(X,p,d,\mathcal H^n)$ be an $\mathrm{RCD}(-\epsilon^2(n-1),n)$ space and
$u:B_2(p)\to\mathbb R^n$ be an $(\epsilon,n)$-splitting map.
For sufficiently small $\epsilon>0$, we have
\[
(1-\Phi(\epsilon\mid n))
d(x,y)^{1+\Phi(\epsilon\mid n)}
\le |u(x)-u(y)|
\le (1+\Phi(\epsilon\mid n))d(x,y)
\]
for all $x,y\in B_1(p)$.
Moreover, if $X$ is a smooth Riemannian $n$-manifold with
$\operatorname{Ric}\ge-\epsilon^2(n-1)$, then $du_x:T_xX\longrightarrow\mathbb R^n$ is invertible for every $x\in B_1(p)$.
\end{theorem}

The transformation theorem also holds under the generalized
Reifenberg condition \cite{HuangHuang2024}. Rather than
introduce this condition, we state the consequences needed below.

\begin{theorem}
For every $n\in\mathbb N$, there exists $\epsilon=\epsilon(n)>0$ with the following property.
Let $(M,p)$ be a pointed Riemannian $n$-manifold such that $\overline B_2(p)$ is compact and
\[
\operatorname{sec}_M\ge-1,
\qquad
\operatorname{Ric}_M\ge-\epsilon,
\]
and let $1\le k\le n$. If $u:B_2(p)\to B_2(0^k)\subset\mathbb R^k$ is an
$(\epsilon,k)$-splitting map, then $u$ is a submersion
on $B_1(p)$.
\end{theorem}

The following topological analogue holds in the nonsmooth setting; see \cite[Proposition~1.17]{HuangHuang2024} and \cite[proof of Theorem~4.1]{HHWZ}.
\begin{lemma}
\label{lemma:topological submersion}
For every $N\in[1,\infty)$, there exist constants
$\delta(N)>0$ and $\epsilon(N)>0$
with the following property.
Let $(X,d,\mathfrak m)$ be an $\mathrm{RCD}(-\epsilon(N),N)$ space
with Alexandrov curvature bounded below by $-1$,
and let $n\le N$ be its rectifiable dimension. Fix $p\in X$.
Suppose that every point of $B_2(p)$ is $(n,\delta(N))$-regular.
Let $1\le k\le n$.
If $u:B_2(p)\to B_2(0^k)\subset\mathbb R^k$ is an
$(\epsilon(N),k)$-splitting map, then $u$ is a topological
submersion on $B_1(p)$.
\end{lemma}

\subsection{Covering spaces and defining subsets}\label{subsec:universal covers}
Let $G$ be a group and let $S \subset G$ be a subset
containing the identity. Let $F$ be the free group on the symbols
$e_g$, $g\in S$. Quotient $F$ by the normal subgroup generated by
\[
e_{g_1}e_{g_2}e_{g_1g_2}^{-1},
\qquad g_1,g_2,g_1g_2\in S.
\]
We denote the resulting group by $\hat S$ and call it the
groupification of $S$.

The map
\[
\iota:S\longrightarrow\hat S,\qquad \iota(g)=[e_g],
\]
preserves every product defined within $S$. There is also a
canonical homomorphism
\[
\pi:\hat S\longrightarrow G,\qquad \pi([e_g])=g.
\]
Since $\pi\circ\iota$ is the inclusion of $S$ into $G$, the map
$\iota$ is injective. Moreover, every local homomorphism
$\phi:S\to H$ to a group $H$ extends uniquely to a homomorphism
$\hat\phi:\hat S\to H$ satisfying $\hat\phi\circ\iota=\phi$.

If $S$ generates $G$, then $\pi$ is surjective. We say that $S$
defines $G$ if $\pi$ is an isomorphism. Thus the elements
of $S$ generate $G$, and the product relations visible within $S$
generate all relations among them.

Every $\mathrm{RCD}(K,N)$ space $(X,d,\mathfrak m)$ admits a
universal cover $\widetilde X$, which is again an
$\mathrm{RCD}(K,N)$ space with the lifted metric and measure
\cite{MondinoWei2019}. Moreover, $X$ is semilocally simply
connected \cite{Wang2022,Wang2021}, so $\pi_1(X,p)$ is naturally
isomorphic to the deck transformation group of $\widetilde X$.

Suppose now that $(X,p,d,\mathfrak m)$ is compact and
$\operatorname{diam}(X)\le D$. Let $(\widetilde X,\tilde p)$ be
its universal cover, and identify $G=\pi_1(X,p)$ with its deck
group. Then the displacement set
$G(\tilde p,r)$
defines $G$ for every $r \ge 20D$; see \cite{FukayaYamaguchi1992,SantosZamora,Wang2023}. 

We will use the following consequence of the proof of \cite[Theorem~3.4]{Wang2023}.
\begin{theorem}\label{theorem:conv to G_0}
Let $(X_i,p_i)$ be pointed proper geodesic metric spaces, and
let $G_i\le\mathrm{Isom}(X_i)$ be closed subgroups. Suppose that
\[
(X_i,p_i,G_i)\xrightarrow{\mathrm{eGH}}(X,p,G),
\]
where $G$ is a Lie group with identity component $G_0$.
Assume that there exists $R>0$ such that $G$ is generated by $G(p,R/2)$ and the set $G_i(p_i,R)$ defines $G_i$ for all sufficiently large $i$. Then there are closed normal subgroups
$\Gamma_i\trianglelefteq G_i$ such that
\[
(X_i,p_i,G_i,\Gamma_i)
\xrightarrow{\mathrm{eGH}}(X,p,G,G_0),
\qquad
G_i/\Gamma_i\cong G/G_0
\]
for all sufficiently large $i$. Moreover, $\Gamma_i$ is generated by elements that are $\epsilon_i$-close to the identity for all sufficiently large $i$.
\end{theorem}

We next consider local relative covers of an $\mathrm{RCD}(K,N)$
space $(X,d,\mathfrak m)$. As the examples in
\cite{SormaniWei2004} show, the universal covers of fixed-radius
balls need not admit a pointed Gromov--Hausdorff convergent
subsequence, even when the ambient compact manifolds have a
uniform lower sectional curvature bound. Passing to relative
covers of smaller balls restores precompactness, as shown by
Xu \cite{Xu2023}.

Fix $p\in X$ and $0<r_1<r_2$. Write $\tilde B(p,r_2)$ for the universal cover of $B_{r_2}(p)$, and let $\tilde p$ be a lift of $p$ to this cover. We denote by $(\tilde{B}(p,r_1,r_2),\tilde p)$ the connected component of the preimage of $B_{r_1}(p)$
containing $\tilde p$, equipped with its intrinsic length
metric and lifted measure $\tilde{\mathfrak m}$. 

\begin{theorem}[Precompactness of relative covers, \cite{Xu2023}]
\label{thm:RCD precompact}
Fix $K\in\mathbb R$, $1\le N<\infty$, and $0<r_1<r_2$.
There exists a positive nondecreasing function
$A_{K,N,r_1,r_2}:(0,\infty)\to(0,\infty)$ such that every
relative cover $\tilde{B}(p,r_1,r_2)$ as above satisfies
\[
0<\tilde{\mathfrak m}(B_r(x))
\le
A_{K,N,r_1,r_2}(r)\,
\tilde{\mathfrak m}(B_{r/2}(x))
\]
for every $x\in\tilde{B}(p,r_1,r_2)$ and $r>0$.

Let $G=\operatorname{im}\bigl[
\pi_1(B_{r_1}(p),p)\longrightarrow
\pi_1(B_{r_2}(p),p)
\bigr]$,
viewed as the deck transformation group of the relative cover.
Then the family of triples $(\tilde{B}(p,r_1,r_2),\tilde p,G)$ is precompact in the pointed equivariant Gromov--Hausdorff
topology.
\end{theorem}

\begin{remark}
In Theorem~\ref{thm:RCD precompact}, the quotient
$\tilde{B}(p,r_1,r_2)/G$ is isometric to $B_{r_1}(p)$
equipped with its intrinsic length metric. This metric may
differ from the restriction of $d$, but the two metrics agree
on $B_{r_1/3}(p)$. Whenever we identify the quotient metric
with $d$, we therefore restrict to such a smaller ball and
usually suppress this distinction in the notation.
\end{remark}

We will also use the following lemma.

\begin{lemma}[{\cite[Lemma~2.10]{Wang2023}}]
\label{lemma:defining_quotient}
Let $S$ be a defining subset of a group $G$, and let
$H\trianglelefteq G$ be a normal subgroup. If $S\cap H$
generates $H$, then $\pi(S)$ defines $G/H$, where
$\pi:G\to G/H$ is the quotient homomorphism.
\end{lemma}

\subsection{Gluing almost splitting maps}

Pointed equivariant Gromov--Hausdorff convergence controls geometry
on bounded regions. When the quotient diameters are uniformly
bounded and the acting groups are identified, local approximations
can be extended to globally defined maps while retaining local approximation estimates.

We call a map $f:X\to Y$
$\epsilon$-injective if $f(x)=f(y)$ implies
$d(x,y)<\epsilon$ for any $x,y\in X$.

\begin{lemma}[Extension property, \cite{Wang2023}]
\label{thm:gluing GHA}
Let $(X,p)$ and $(Y,q)$ be proper geodesic metric spaces equipped
with isometric $G$-actions. Suppose that
\[
10\epsilon<D<\frac{1}{10\epsilon},
\qquad
\operatorname{diam}(X/G)\le D,
\qquad
\operatorname{diam}(Y/G)\le D.
\]
Let $f:B_{4D}(p)\longrightarrow B_{4D}(q)$
be an $\epsilon$-almost $G$-equivariant $\epsilon$-eGHA. 
Then there exists a map $f':X\to Y$ such that:
\begin{enumerate}
\item
$d_Y(f'(x),f(x))\le2\epsilon$ for every $x\in B_{4D}(p)$;

\item
$f'$ is a $2\epsilon$-GHA on any $D$-ball and uniformly $2\epsilon$-almost $G$-equivariant on $X$,
\[
d_Y\bigl(f'(gx),gf'(x)\bigr)\le 2\epsilon
\qquad(g\in G,\ x\in X).
\]

\item
$f'$ is $2\epsilon$-injective and $f'(X)$ is $3\epsilon$-dense in $Y$.
\end{enumerate}

\end{lemma}

We recall the extension argument in Lemma \ref{thm:gluing GHA}.
For each $x\in X$, choose $g\in G$ such that
$gx\in B_{2D}(p)$, and define $f'(x):=g^{-1}f(gx)$.
Choosing a different $g$ changes the value of $f'(x)$ by at most $\epsilon$. The map $f'$ satisfies the properties in the lemma. 

We now turn to gluing almost splitting maps. Fix a sufficiently
small $\epsilon>0$. For a pointed metric measure space
$(X,p,d,\mathfrak m)$ and a Riemannian $k$-manifold $(K,g_K)$,
we impose the following geometric conditions:
\begin{equation}\label{eq:condition}
\begin{gathered}
(X,d,\mathfrak m)\in
\mathrm{RCD}(-(N-1)\epsilon^2,N),
\qquad
\dim K=k\le N,
\qquad
\operatorname{inj}_K\ge10\epsilon^{-1},\\
\text{for every }y\in K,\quad
(\exp_y)^*g_K
\text{ is }\epsilon^2\text{-}C^\infty\text{-close to the flat metric}\\
\text{on }B_{10/\epsilon}(0_y)\subset T_yK.
\end{gathered}
\end{equation}

For a pair $(X,K)$ satisfying \eqref{eq:condition}, a continuous
map $f:X\to K$ is called locally $(\delta,k)$-splitting
if, for every $x\in X$ with $y=f(x)$, we have $f(B_{1}(x))\subset B_{1+\delta}(y)$
and the map
\[
\exp_y^{-1}\circ f:B_{1}(x)\longrightarrow T_yK\cong\mathbb R^k
\]
is a $(\delta,k)$-splitting map.

An averaging argument gives the following theorem; see \cite{Huang2020,Wang2024}.
\begin{theorem}\label{thm:ave splitting}
There exists $\epsilon(N)>0$ such that the following holds
for every $0<\epsilon\le\epsilon(N)$.
Let $(X,K)$ satisfy~\eqref{eq:condition} with parameter
$\epsilon$. Suppose that $f':X\to K$ has $\epsilon$-dense
image and that, for every $x\in X$, there is an
$\epsilon$-Gromov--Hausdorff approximation onto its image
\[
(u,v):B_{100}(x)\longrightarrow\mathbb R^k\times Y,
\]
where $Y$ is either a point or an $\mathrm{RCD}(0,N-k)$
space, with $N-k\ge1$ in the latter case. Assume that $u$
is an $(\epsilon,k)$-splitting map and
\[
\sup_{B_{100}(x)}
\bigl|u-\exp_y^{-1}\circ f'\bigr|\le\epsilon,
\qquad y=f'(x),
\]
after an orthonormal identification $T_yK\cong\mathbb R^k$.

Then there exists a continuous map $f:X\to K$ that is
uniformly $\Phi(\epsilon\mid N)$-close to $f'$ and locally
$(\Phi(\epsilon\mid N),k)$-splitting.
\end{theorem}

The next theorem makes this construction equivariant.
\begin{theorem}[$G$-equivariant locally almost splitting map]
\label{thm:inv}
In the setting of Theorem~\ref{thm:ave splitting}, suppose
that a virtually nilpotent Lie group $G$ acts on $X$ by measure-preserving
isometries and on $K$ by isometries, and that both actions identify $G$
with a closed subgroup of the corresponding isometry group.
Assume that $G$ is compactly generated and that $f'$ is
uniformly $\epsilon$-almost $G$-equivariant:
\[
d_K\bigl(f'(gx),gf'(x)\bigr)\le\epsilon
\qquad(g\in G,\ x\in X).
\]
Then there exists a continuous $G$-equivariant map
$h:X\to K$ that is uniformly $\Phi(\epsilon\mid N)$-close
to $f'$ and locally $(\Phi(\epsilon\mid N),k)$-splitting.
\end{theorem}

\begin{proof}
Choose a compact symmetric generating neighborhood
$S\subset G$ containing the identity, and let $\mathrm m$
be a Haar measure on $G$. Since $G$ is virtually nilpotent,
it is unimodular, so $\mathrm m$ is both left- and
right-invariant. For discrete $G$, we use counting measure.

Let $f:X\to K$ be the map supplied by
Theorem~\ref{thm:ave splitting}. Its uniform closeness to
$f'$ implies that it is $\Phi(\epsilon\mid N)$-almost $G$-equivariant.

For $j\ge1$, define $f_j(x)$ to be the center of mass of the family $\{g^{-1}f(gx):g\in S^j\}$ with respect to normalized Haar measure on $S^j$.
The center-of-mass estimates used in the equivariant
averaging constructions of \cite{Wang2024,HHWZ} show that $f_j$ is locally
$(\Phi(\epsilon\mid N),k)$-splitting and $\Phi(\epsilon\mid N)$-close to $f$, where $\Phi(\epsilon\mid N)$ is independent of $j$.

In particular, the maps $f_j$ have uniform local Lipschitz
bounds and uniform bounds for their coordinate Laplacians.
After passing to a diagonal subsequence, they converge
locally uniformly to a continuous map $h:X\to K$.
The Laplacian bounds also give strong
$H^{1,2}_{\mathrm{loc}}$ convergence in fixed target
charts; see \cite[Theorem~4.4]{AH18}. Consequently, the
splitting estimates pass to the limit. Thus $h$ is
uniformly $\Phi(\epsilon\mid N)$-close to $f'$ and locally
$(\Phi(\epsilon\mid N),k)$-splitting.

It remains to prove equivariance. The volume asymptotic
for groups of polynomial growth gives constants
$d(G)\in\mathbb Z_{\ge0}$ and $c(S)>0$ such that
\[
\lim_{j\to\infty}
\frac{\mathrm m(S^j)}{j^{d(G)}}=c(S);
\]
see \cite[Theorem~1.1]{Breuillard2014}. In particular, $\frac{\mathrm m(S^{j+1})}{\mathrm m(S^j)}\longrightarrow1$.
For every $g_0\in S$, both $S^j$ and $S^jg_0$ are
contained in $S^{j+1}$ and have the same measure. Hence
\[
\frac{\mathrm m(S^j\triangle S^jg_0)}
     {\mathrm m(S^j)}
\le
2\,\frac{\mathrm m(S^{j+1})-\mathrm m(S^j)}
        {\mathrm m(S^j)}
\longrightarrow0.
\]

The invariance of the center of mass under isometries and the right-invariance of $\mathrm m$ show that $g_0^{-1}f_j(g_0x)$ is the center of mass of the family $\{g^{-1}f(gx):g\in S^jg_0\}$ with respect to normalized Haar measure on $S^jg_0$. The preceding
symmetric-difference estimate and stability of centers
of mass therefore imply
\[
\sup_{x\in X}
d_K\bigl(g_0^{-1}f_j(g_0x),f_j(x)\bigr)
\longrightarrow0.
\]
Passing to the locally uniform limit gives $h(g_0x)=g_0h(x)$ for $g_0\in S$ and $x\in X$.
Since $S$ generates $G$, the map $h$ is $G$-equivariant.
\end{proof}

\begin{remark}
The proof of Theorem~\ref{theorem:smooth bundle over nil} requires a smooth equivariant almost splitting map. In the setting of Theorem~\ref{thm:inv}, suppose that $X$ is a smooth $n$-manifold with almost nonnegative Ricci curvature and $\operatorname{diam}(X/G)\le D$. Following the extension construction, we construct $f$ using smooth harmonic local splitting maps and smooth cutoffs.
For each $k\ge0$, the resulting map has uniform $C^k$ bounds on every $10D$-ball.
These bounds may depend on $k$, the spaces, and the chosen map $f$, but not on the center of the ball. For each fixed $X$, the averaged maps $f_j$ satisfy corresponding bounds independent of $j$. A diagonal subsequence therefore converges in $C^\infty_{\mathrm{loc}}$ to a smooth $G$-equivariant almost splitting map $h$.
\end{remark}

We give a corollary of the gluing argument. 
\begin{lemma}\label{lemma:no compact limit}
Let $(X_i,p_i,d_i,\mathcal H^n)$ be noncollapsed
$\mathrm{RCD}(-\epsilon_i,n)$ spaces with $\epsilon_i\to0$,
and suppose that
\[
(X_i,p_i,d_i,\mathcal H^n)
\xrightarrow{\mathrm{pmGH}}
(\mathbb R^n,0^n,d_{\mathrm{Eucl}},\mathcal H^n).
\]
If finite groups $G_i\le\mathrm{Isom}(X_i)$ act freely
and satisfy $G_i=G_i(p_i,R)$
for some fixed $R>0$, then $G_i$ is trivial for all
sufficiently large $i$.
\end{lemma}

\begin{proof}
After passing to a subsequence, take an equivariant limit
\[
(X_i,p_i,G_i)\xrightarrow{\mathrm{eGH}}
(\mathbb R^n,0^n,G).
\]
The displacement bound makes $G$ compact. 

By Lemma~\ref{lemma:homomorphic eGHA}, we may choose the group maps
$\phi_i:G_i\to G$ to be homomorphisms. Volume convergence implies
that these maps are injective for all sufficiently large $i$.
We therefore identify each $G_i$ with its image in $G$.

Choose a fixed point $q\in\mathbb R^n$ of $G$ and points $q_i\to q$. Construct almost splitting maps on $B_{10}(q_i)$ and average over the finite groups $G_i$. The local averaging estimates and the canonical Reifenberg theorem give a $G_i$-equivariant homeomorphism from $B_9(q_i)$ to $B_9(q)$. The preimage of $q$ is fixed by every element of $G_i$, contradicting freeness whenever $G_i$ is nontrivial.
\end{proof}

\subsection{Properties of nilpotent groups}
We collect some basic facts about nilpotent groups. We write $C_g(h)=ghg^{-1}$ for conjugation and $\operatorname{Ad}(g)=d_eC_g$ for the induced map on the Lie algebra. 
\begin{defn}\label{defn:infra}
Let $G$ be a connected and simply connected nilpotent Lie group, and let $\Gamma\le G\rtimes\mathrm{Aut}(G)$ be a torsion-free discrete subgroup. If $[\Gamma:\Gamma\cap G]< \infty$ and $G/\Gamma$ is compact, then $G/\Gamma$ is an infranilmanifold.
\end{defn}

\begin{lemma}\label{thm:G_0}
Let $G$ be a nilpotent Lie group with identity component $G_0$. Then $[K,G_0]=\{e\}$ for every compact subgroup $K\le G$.
\end{lemma}

\begin{lemma}
Let $G$ be a finitely generated virtually nilpotent discrete group and let $H\trianglelefteq G$. Then $\operatorname{rank}(G)=\operatorname{rank}(H)+\operatorname{rank}(G/H)$.
\end{lemma}

\begin{lemma}\label{lem:pure translation}
Let $H\le\mathrm{Isom}(\mathbb R^l)$ be a closed nilpotent subgroup with compact quotient $\mathbb R^l/H$ and $\dim H=k$. Then $H\cong\mathbb R^k\times\mathbb Z^{l-k}$ and $H$ consists of pure translations.
\end{lemma}
\begin{proof}
By the proof of \cite[Corollary~4.1.13]{Thurston}, the identity
component $H_0$ is abelian. Let $E$ be the unique maximal affine
subspace on which $H_0$ acts by translations;
see \cite[Proposition~4.2.4]{Thurston}.
Since $H_0$ is normal in $H$, uniqueness implies that $E$ is
$H$-invariant. If $E\ne\mathbb R^l$, distance to $E$ would descend
to an unbounded continuous function on the compact quotient
$\mathbb R^l/H$, a contradiction.
Thus $E=\mathbb R^l$, and $H_0$ consists of pure translations.

Let $V\subset\mathbb R^l$ be the translation subspace of $H_0$,
and write $k=\dim V$. Then $H_0\cong\mathbb R^k$ acts transitively
on each affine subspace parallel to $V$.
For every $h\in H$, its linear part preserves $V$.
Nilpotence implies that the restriction to $V$ is unipotent.
Since it is also orthogonal, it is the identity.

After removing the $V$-component of each element's translation
part using $H_0$, we identify $H/H_0$ with a closed discrete
subgroup of $\mathrm{Isom}(V^\perp)$.
Its action is cocompact, so it is a nilpotent crystallographic
group of dimension $l-k$.
By Bieberbach's theorem, it contains a normal translation lattice of finite index.
Nilpotence implies that conjugation by every element of $H/H_0$
is unipotent on $V^\perp$.
This action is also orthogonal and hence trivial.
Therefore every element of $H/H_0$ acts by pure translations
on $V^\perp$.
Consequently, every element of $H$ acts by pure translation
on $\mathbb R^l$, and $H\cong\mathbb R^k\times\mathbb Z^{l-k}$.
\end{proof}

\begin{lemma}\label{lemma:group splitting}
Let $G$ be a nilpotent Lie group with finitely many components, and suppose that $G_0$ is simply connected. Then there is a compact subgroup $K$ such that $G\cong G_0\times K$.
\end{lemma}
\begin{proof}
By \cite[Lemma 3.2(b)]{Wilking2000}, there is a maximal compact subgroup $K$ such that $G\cong G_0\rtimes K$. For $k\in K$, let $A:=\operatorname{Ad}(k)$ be the induced action on the Lie algebra of $G_0$. Nilpotence gives $(A-\mathrm{Id})^c=0$ for some $c\in\mathbb N$.
Finiteness of $K$ gives $A^m=\mathrm{Id}$ for some $m\in\mathbb N$.
Thus $A$ is both unipotent and of finite order, so $A=\mathrm{Id}$.
Therefore $G\cong G_0\times K$.
\end{proof}

\begin{corollary}
A compactly generated nilpotent Lie group $G$ has no nontrivial compact subgroup if and only if $G_0$ is simply connected and $G/G_0$ is torsion-free.
\end{corollary}
\begin{proof}
The reverse implication is immediate, so we prove the forward implication. Suppose that $G$ has no nontrivial compact subgroup. Then $G_0$ is simply connected. If $G/G_0$ contained a nontrivial torsion element, there would be $g\notin G_0$ with $g^m\in G_0$ for some $m>1$. The subgroup $H:=\langle g,G_0\rangle$ has finitely many components, so Lemma~\ref{lemma:group splitting} gives $H\cong G_0\times K$. Here $K$ is a nontrivial compact subgroup of $G$, a contradiction.
\end{proof}

\section{Examples}
We use the following example to illustrate the constructions in Theorems~A and~B. Let
\[
 N=\left\{n(a,b,d)=
 \begin{pmatrix}1&a&d\\0&1&b\\0&0&1\end{pmatrix}:a,b,d\in\mathbb R\right\},
 \qquad
 \Gamma=\{n(a,b,d):a,b,d\in\mathbb Z\}.
\]
The lattice generators $x=n(1,0,0)$, $y=n(0,1,0)$ and
$z=n(0,0,1)$ satisfy $[x,y]=z$, with $z$ central. Their logarithms
$X,Y,Z$ satisfy $[X,Y]=Z$. Give $N$ the left-invariant metric $g_i$
with these vectors orthogonal and
\[
 |X|_{g_i}=|Y|_{g_i}=i^{-2},\qquad
 |Z|_{g_i}=i^{-5}.
\]
In the $g_i$-orthonormal basis $(i^2X,i^2Y,i^5Z)$, we have
\[
[i^2X,i^2Y]=i^{-1}(i^5Z).
\]
Thus, $|\sec_{N,g_i}|\to0$,
$(N,g_i,e_N)$ converges to $(\mathbb R^3,0^3)$, and
$\operatorname{diam}(N/\Gamma,g_i)\to0$.

Next consider the unit sphere
\[
S^3=\{(w_1,w_2)\in\mathbb C^2:|w_1|^2+|w_2|^2=1\}.
\]
The torus $T^2=\mathbb R^2/\mathbb Z^2$ acts on $S^3$ by rotating the two complex coordinates separately. The quotient is $S^3/T^2=[0,\pi/2]$. We choose the quotient projection $\pi:S^3\to[0,\pi/2]$ so that $\pi(0,1)=0$.

\subsection{A fiber bundle over a nilmanifold}
Let $\widetilde X_i=(N,g_i)\times S^3$.
Let $\hat x_i$ and $\hat y_i$ act on $N$ by left translations
by $x$ and $y$, and on $S^3$ by the torus elements
$(1/i,0)$ and $(0,1/i)$, respectively.
Let $\hat z_i$ act by left translation by $z$ on $N$ and
trivially on $S^3$, while $h$ acts trivially on $N$ and by
$(1/2,1/2)$ on $S^3$.

Set $G_i:=\langle\hat x_i,\hat y_i,\hat z_i,h\rangle$
and $X_i:=\widetilde X_i/G_i$.
Choose $\tilde p_i:=(e_N,(0,1))$ and
$p_i:=\pi(\tilde p_i)$. The corresponding limit basepoint in the covering spaces is
$\tilde p:=(0^3,(0,1))\in\mathbb R^3\times S^3$.

On $(N,g_i)$, both $\hat x_i$ and $\hat y_i$ move $e_N$ by $i^{-2}$, while on $S^3$ each has maximal displacement of order $i^{-1}$. We obtain
\[
\begin{CD}
(\widetilde X_i,\tilde p_i,G_i) @>{\mathrm{eGH}}>> (\mathbb R^3\times S^3,\tilde p,G)\\
@VV\pi V @VV\pi V\\
(X_i,p_i) @>{\mathrm{GH}}>> S^3/T^2=[0,\pi/2],
\end{CD}
\]
where $G=\mathbb R^3\times T^2$. In particular, the diameters of $X_i$ are uniformly bounded.

We now illustrate the proof of Theorem~A for this convergence.
Since $G_i$ is nilpotent, we may take $\Gamma_i=G_i$ and
$\hat X_i=\widetilde X_i/\Gamma_i=X_i$ in the notation of
Section~\ref{section:model}. 

Let $\psi_i:G_i\to G$ denote the forward group maps in the
equivariant approximations, and let
$\rho_1:G\to\mathbb R^3$ be the projection.
The composites $\varphi_i:=\rho_1\circ\psi_i$ form a discrete approximation of $\mathbb R^3$.
The small-subgroup construction gives
the finite normal subgroups
$H_i=\{e,h\}\trianglelefteq\Gamma_i$.
After passing to a subsequence, these converge
equivariantly to $H=\{e,h\}\le T^2$.

We next consider the equivariant convergence of the quotient actions
\[
(\widetilde X_i/H_i,\bar p_i,\Gamma_i/H_i)
\longrightarrow
(\mathbb R^3\times(S^3/H),\bar p,G/H).
\]
Here $\bar p_i$ is the image of $\tilde p_i$, and
$G/H=\mathbb R^3\times(T^2/H)$.
Let $\bar\varphi_i:\Gamma_i/H_i\to\mathbb R^3$ denote the induced quotient discrete approximations, which have no small subgroups.
For every fixed $R>0$ and all sufficiently large $i$,
Theorem~\ref{thm:large nil} gives a nilprogression of rank $3$
containing the inverse image of $B_R(0^3)$ under this
quotient approximation.

Since $\widetilde X_i$ is simply connected and
$\hat X_i=\widetilde X_i/\Gamma_i$ has uniformly bounded
diameter, a sufficiently large displacement set defines
$\Gamma_i$. The same holds for $\Gamma_i/H_i$, since the
subgroups $H_i$ are uniformly bounded.
Combining this defining-set property with the nilprogression
structure shows that $\Gamma_i':=\Gamma_i/H_i$ is a
torsion-free nilpotent group.

Let $N_i$ be the Malcev completion of $\Gamma_i'$. We can choose a left-invariant metric on $N_i$ that is $C^\infty$-close to the Euclidean metric on bounded sets. For $g\in\Gamma_i'(\bar p_i,R)$, its action on $N_i$ is close to translation by $\bar{\varphi}_i(g)$ on $\mathbb R^3$.

Both $(\widetilde X_i/H_i)/\Gamma_i'$ and $N_i/\Gamma_i'$
have uniformly bounded diameter.
The extension and gluing arguments therefore give a globally defined
smooth $\Gamma_i'$-equivariant map
$f_i:\widetilde X_i/H_i\to N_i$ that is almost splitting
on every fixed $R$-ball.
Huang--Huang's nondegeneracy theorem implies that $df_i$ is
surjective everywhere. Thus $f_i$ descends to a smooth
fiber bundle map $\bar f_i:X_i\longrightarrow N_i/\Gamma_i'$.

The groups $\Gamma_i'$ and $N_i$ are constructed solely from the equivariant convergence and the algebraic structure theorems. The construction does not use the product structure of $\widetilde X_i$, which is unavailable in the general proof of Theorem~A.

\subsection{A singular fibration over an interval}\label{example B}
We keep the same universal covers $\widetilde X_i$ and now take $G_i:=\langle\hat x_i,\hat y_i,\hat z_i\rangle$, omitting the antipodal action. In this subsection, write $X_i:=\widetilde X_i/G_i$ and let $p_i$ be the image of $\tilde p_i$. The same equivariant convergence holds, with limit group
$G=\mathbb R^3\times T^2$. We construct singular fibrations $\bar{F}_i:X_i\to[0,\pi/2]$ from this equivariant convergence. Although the general proof of Theorem~B uses local universal covers, the global universal covers $(\widetilde X_i,G_i)$ suffice in this example.

We begin near the endpoint $0\in[0,\pi/2]$.
Take $\tilde p=(0^3,0,1)\in\mathbb R^3\times S^3$, so that
$\pi(\tilde p)=0$.
We first describe the construction on the equivariant tangent
cone at $\tilde p$, then transfer the resulting singular
fibration to a neighborhood of $0$ in the quotient.

Choose $\tilde p_i\in\widetilde X_i$ converging to $\tilde p$, and let $r_i\to\infty$ sufficiently slowly. Set $Y_i:=r_i\widetilde X_i$ and write $G_i^{r_i}$ for the action of $G_i$ on $Y_i$. Then
\[
(Y_i,\tilde p_i,G_i^{r_i})
\longrightarrow
(\mathbb R^6,q,H=\mathbb R^4\times S^1).
\]
Apply Theorem~\ref{thm:nil} to the discrete approximation $G_i^{r_i}\to H$ on a small $2r$-neighborhood of $e_H$. We obtain a nilprogression $P_i$ of rank $5$ containing an $r$-neighborhood of the identity in $G_i^{r_i}$. Write $a_i,b_i,c_i,d_i,e_i$ for the formal Malcev generators corresponding respectively to $\hat x_i,\hat y_i,\hat x_i^i,\hat y_i^i,\hat z_i$. 

Consider the truncated set of powers
$\{\hat x_i^j:|j|<K_i\}$, where $K_i$ is the first positive
integer for which $\hat x_i^{K_i}$ leaves the $r$-ball
about the identity in $G_i^{r_i}$.
By our choice of metric, this set converges to a subset
of the $S^1$-factor of $H$.
The corresponding truncated sets of powers of
$\hat y_i,\hat x_i^i,\hat y_i^i,\hat z_i$ converge to intervals
in four mutually perpendicular directions in $\mathbb R^4$.

The groupification of $P_i$ is a torsion-free nilpotent group
$\widetilde\Gamma_i$ of rank $5$. Its Malcev completion
$\widetilde N_i$ admits a left-invariant metric that is
$C^\infty$-close to the Euclidean metric on bounded sets.
Moreover, the action of each $g\in P_i$ on an $r$-ball in
$\widetilde N_i$ is close to the action of its image on an
$r$-ball about $e_H$, locally identified with a ball in
$\mathbb R^5$.

The natural projection $\pi:\widetilde\Gamma_i\to G_i^{r_i}$ has a rank-$2$ kernel generated by $a_i^ic_i^{-1}$ and $b_i^id_i^{-1}$.
For the rescaling sequence under consideration, the relation $a_i^i=c_i$ is visible on a bounded ball and converges to the period relation defining the kernel of $\mathbb R^5\to H$. The relation $b_i^i=d_i$ occurs only far from the identity and is omitted from the local construction. Let $\hat G_i$ be the groupification of $S_i:=G_i^{r_i}(\tilde p_i,6R_i)$, where $R_i\to\infty$ sufficiently slowly. Form the gluing space
\[
\hat X_i:=\hat G_i\times_{S_i}\overline B_{3R_i}^{Y_i}(\tilde p_i),
\]
using the construction of Section~\ref{subsection:same orbits}; let $\hat p_i=[e,\tilde{p}_i] \in \hat{X}_i$ be the basepoint. The purpose of $\hat X_i$ is to remove the relation realized far from the identity. In this example, $\hat X_i$ is the universal cover of a small neighborhood of the local $G_i$-orbit through $\tilde p_i$. The projection of this orbit to the $S^3$-factor is close to
the second-coordinate circle $\{(0,w_2)\in S^3:|w_2|=1\}$.

Then $(\hat X_i,\hat p_i,\hat G_i)$ agrees with $(Y_i,\tilde p_i,G_i^{r_i})$ on the $R_i$-ball, and hence
\[
(\hat X_i,\hat p_i,\hat G_i)
\longrightarrow
(\mathbb R^6,q,H=\mathbb R^4\times S^1).
\]
The local construction near \(\hat p_i\) transfers to a neighborhood of \(\tilde p_i\) and induces a singular fibration on the quotient.

Set $D_i:=\langle a_i^ic_i^{-1}\rangle$, $N_i:=\widetilde N_i/D_i$ and $\Gamma_i:=\widetilde\Gamma_i/D_i$.
With the induced metric, $N_i$ is $C^\infty$-close to $H$ on bounded sets. In the notation of Section~\ref{sec 6}, the group obtained by gluing is $\hat G_i=\Gamma_i$ in this example, since $H$ is connected.

Let $T_i$ be the maximal torus of $N_i$, and identify it
with the $S^1$-factor of $H$.
The orthogonal action of this circle on the normal slice
$\mathbb R^2$ defines the associated bundle
\[
M_i:=N_i\times_{T_i}\mathbb R^2.
\]
Equipped with the induced connection metric, $M_i$ is
$C^\infty$-close to $\mathbb R^6$ on bounded sets and
carries an isometric $\hat G_i$-action.
With $E_i=N_i$, $(M_i,N_i,\hat G_i,E_i)$
is the nilpotent model used near $\hat p_i$.
Write $p_i'=[e_{N_i},0^2]\in M_i$.

The group $\hat G_i$ acts on $M_i$ through the $N_i$-factor.
For a sufficiently large fixed $R$, the gluing construction
gives an $\epsilon_i$-injective, $\hat G_i$-equivariant map
\[
f_i:B_R(\hat G_i\hat p_i)
\longrightarrow B_R(\hat G_ip_i')\subset M_i
\]
that is almost splitting on each $R$-ball.
The canonical Reifenberg theorem and $\epsilon_i$-injectivity
imply that $f_i$ is a homeomorphism onto its image.

The construction of $N_i$ provides a diffeomorphism
$\varphi_i:N_i\to H$ realizing the $C^\infty$ approximation
on bounded sets.
Its restriction to $T_i$ is an isomorphism onto the
$S^1$-factor of $H$, and it satisfies
\[
\varphi_i(g\tau)=\varphi_i(g)\varphi_i(\tau),
\qquad g\in N_i,\ \tau\in T_i.
\]
Consequently, $\varphi_i$ induces a diffeomorphism
\[
J_i:M_i\to\mathbb R^6,\qquad
J_i([g,v])=[\varphi_i(g),v],
\]
taking $N_i$-orbits bijectively to $H$-orbits.
Thus $F_i:=J_i\circ f_i$ is the desired orbit-preserving
homeomorphism. We return to the original scale and pass to the quotient to obtain a local singular fibration on $X_i$.

For $t$ sufficiently close to $0$, let $y_i\in M_i$
represent the corresponding quotient point.
The fiber over $t$ is homeomorphic to $N_i(y_i)/\hat G_i$.
It is a four-dimensional nilmanifold when $t=0$ and a
five-dimensional nilmanifold when $t>0$.
The circle isotropy at $t=0$ accounts for the drop
in fiber dimension.

The same construction at $\pi/2$ and regular
fibration charts over $(0,\pi/2)$ give local singular
fibrations on a finite open cover of $[0,\pi/2]$.
Each endpoint chart restricts to a regular fibration over
its intersection with the open interval.
On overlaps, we work in a common local covering space
and compare the corresponding maps into the nilpotent models.
The equivariant gluing argument of Section~7 then yields
a global singular fibration
$\bar F_i:X_i\to[0,\pi/2]$.

\section{Construction of fiber bundles}\label{section:model}
\subsection{Fiber bundles over nilmanifolds}
We prove the first part of Theorem~A. Consider a sequence of metric measure spaces $(X_i,p_i,d_i,\mathfrak m_i)\in\mathrm{RCD}(-\epsilon_i,N)\cap\mathrm{Alex}(-1)$ with $\epsilon_i\to0$ and $\operatorname{diam}(X_i)\le D$, satisfying the pointwise regularity assumption of Theorem~A. After passing to a subsequence, we may assume that the spaces $X_i$ have the same rectifiable dimension $k_0\le N$. Thus every point of $X_i$ is $(k_0,\delta(N))$-regular. 
It suffices to prove the conclusion for all sufficiently large $i$.

By \cite{MondinoWei2019}, $X_i$ has a universal cover $(\widetilde X_i,\tilde p_i)$ that is an $\mathrm{RCD}(-\epsilon_i,N)$ space with the lifted metric and measure. The lifted metric also has Alexandrov curvature bounded below by $-1$. The deck group $G_i$ of this cover is isomorphic to $\pi_1(X_i,p_i)$. For all sufficiently large $i$, the Margulis lemma gives a normal nilpotent subgroup $\Gamma_i\trianglelefteq G_i$ of index at most $C(N)$, with a nilpotent basis of length at most $k_0$. 

After renormalizing the measures and passing to a subsequence, we obtain
\begin{center}
 $\begin{CD}
		(\widetilde{X}_i,\tilde{p}_i,G_i,\Gamma_i,\tilde{\mathfrak{m}}_i) @>{\mathrm{emGH}}>> (Y,\tilde{p},G,G',\tilde{\mathfrak{m}})\\
		@VV\pi V @VV\pi V\\
		(\hat{X}_i=\widetilde{X}_i/\Gamma_i,\hat{p}_i,G_i/\Gamma_i,\hat{\mathfrak{m}}_i) @>{\mathrm{emGH}}>> (\hat{X},\hat{p},G/G',\hat{\mathfrak{m}})\\
		@VV\pi V @VV\pi V\\
		(X_i,p_i,\mathfrak{m}_i) @>{\mathrm{mGH}}>> (X,p,\mathfrak{m}).
	\end{CD}$
\end{center} 
Then $\operatorname{diam}(\hat{X}) \le C(N)D$. The following theorem implies the first part of Theorem A.

\begin{theorem}[First part of Theorem~A]\label{thm:first}
For all sufficiently large $i$, the space $\hat X_i$ is the
total space of a topological fiber bundle over a nilmanifold,
with fibers having finite fundamental group.
\end{theorem}

The spaces $\widetilde X_i$ satisfy $\mathrm{RCD}(-\epsilon_i,N)$, so their limit $Y$ satisfies $\mathrm{RCD}(0,N)$. Since $Y/G'$ is compact, the splitting theorem gives $Y=\mathbb R^k\times Z$ with $Z$ compact; see \cite{Gigli2013}. 

Let $\rho_1:\mathrm{Isom}(Y)\to\mathrm{Isom}(\mathbb R^k)$
be the projection, and let
$\psi_i:(G_i,\Gamma_i)\to(G,G')$
denote the group maps in the equivariant approximations. Define $\varphi_i:=\rho_1\circ\psi_i$ and $\Gamma:=\rho_1(G')$. Compactness of $\hat X$ implies compactness of $\mathbb R^k/\Gamma$. Since the groups $\Gamma_i$ admit nilpotent bases of uniformly
bounded length, their limit $G'$ and its image $\Gamma$ are
nilpotent. Lemma~\ref{lem:pure translation} therefore shows that
$\Gamma$ consists of pure translations and
$\Gamma\cong\mathbb R^{k_1}\times\mathbb Z^{k_2}$,
where $k_1+k_2=k$.

Since $Z$ is compact, so is $\mathrm{Isom}(Z)$. Thus the maps $\varphi_i:\Gamma_i\to\Gamma$ form a discrete
approximation. The small-subgroup construction in
Subsection~\ref{subsection:local group} gives finite maximal
small subgroups $H_i\le\Gamma_i$ with respect to $\varphi_i$. The group $G_i$ is generated by $G_i(\tilde p_i,3D)$; thus $H_i$ is normal in $G_i$. After passing to a subsequence, assume that $H_i$ converges
equivariantly to $H\le\{e\}\times\mathrm{Isom}(Z)$. Set $R:=2\operatorname{diam}(Z)+1$. For all sufficiently large $i$, we have $H_i\subset\Gamma_i(\tilde p_i,R)$.   

The quotient actions converge as follows:
\[
(\widetilde{X}_i/H_i,\bar{p}_i,\Gamma_i/H_i) \overset{\mathrm{emGH}}\longrightarrow ( \mathbb{R}^k \times (Z/H), \bar{p},G'/H).
\]
Denote the induced quotient discrete approximations by
$\bar\varphi_i:\Gamma_i/H_i\to\Gamma$. 
We have $\operatorname{diam}(Z/H)\le R/2$. The action of $G'/H$ on $\mathbb R^k\times(Z/H)$ preserves the product decomposition. Let $\bar\rho_1:G'/H\to\Gamma$ be the induced projection. Then
\[
\bar\rho_1^{-1}(B_r(e_\Gamma))
\subset
(G'/H)(\bar p,r+R/2)
\qquad(r>0).
\] 

Equip $\Gamma$ with the metric induced by its translation action on $\mathbb R^k$. Apply Theorem~\ref{thm:large nil} to the quotient approximations $\bar\varphi_i:\Gamma_i/H_i\to\Gamma$. For each sufficiently large fixed $r$ and all sufficiently large $i$ depending on $r$, we obtain a nilprogression $P_i'$ of rank $k$ satisfying
\[
\bar\varphi_i^{-1}(B_{\Lambda(r)}(e_\Gamma))\subset P_i',
\qquad
\bar\varphi_i(P_i')\subset B_{C(r+1)}(e_\Gamma).
\] 

For sufficiently large fixed $r$ and all sufficiently large $i$,
the groupification of $P_i'$ is a torsion-free nilpotent group
$\Gamma_i'$ of rank $k$. Malcev's theorem embeds $\Gamma_i'$ as a lattice in a $k$-dimensional connected, simply connected nilpotent Lie group $N_i$. By Theorem~\ref{thm:Q conv}, the Lie algebra structures converge to the abelian structure.

Write the ordered generators of $P_i'$ as $(u_{i,1},\ldots,u_{i,k})$. For $1\le j\le k_2$, the images $\bar\varphi_i(u_{i,j})$ converge to elements $u_j\in\Gamma$ whose component classes form a basis of $\Gamma/\Gamma_0$. For $k_2+1\le j\le k$, the images $\bar\varphi_i(u_{i,j})$ converge to $e_\Gamma$. The following lemma is proved in \cite{Wang2024}; we include a brief proof for the reader's convenience. 

\begin{lemma}\label{lemma:metric}
The group $N_i$ admits a left-invariant metric that is $C^\infty$-close to the Euclidean metric on bounded sets. For every $g_i\in P_i'$, its action on $N_i$ is close to the translation by $\bar\varphi_i(g_i)$ on $\mathbb R^k$.
\end{lemma}
\begin{proof}
First fix a sufficiently large auxiliary scale $L$, then choose $r$ sufficiently large, and finally take $i$ sufficiently large. Let $s_{i,j}\in\mathbb N$ be the smallest positive integer such that $\bar\varphi_i(u_{i,j}^{s_{i,j}})\notin B_L(e_\Gamma)$. After passing to a subsequence, suppose that $\bar\varphi_i(u_{i,j}^{s_{i,j}})\to v_j\in\Gamma$. Choose $v_{i,j}$ in the Lie algebra of $N_i$ with $\exp_{N_i}(v_{i,j})=u_{i,j}^{s_{i,j}}$, and define a left-invariant metric by
\[
g(v_{i,j},v_{i,j'})=v_j\cdot v_{j'}
\qquad(1\le j,j'\le k).
\]
The conclusion follows from Theorem~\ref{thm:Q conv}.  
\end{proof}

We next show that $\Gamma_i'$ is naturally isomorphic to $\Gamma_i/H_i$. Put $R_0:=\max\bigl(\{0\}\cup\{|u_j|:1\le j\le k_2\}\bigr)$, where $|u_j|$ is the Euclidean translation length. 

\begin{lemma}\label{lemma:radius control}
For all sufficiently large $i$, the displacement of $u_{i,j}$ at $\bar p_i$ is at most $R_0+R$ for $1\le j\le k_2$, and at most $R$ for $k_2<j\le k$. 
\end{lemma}
\begin{proof}
Recall that $R\ge2\operatorname{diam}(Z/H)+1$. For $1\le j\le k_2$, the element $u_{i,j}\in\Gamma_i/H_i$ is equivariantly close to an element with action $(u_j,h_j)$ on $\mathbb R^k\times(Z/H)$. The isometry $(u_j,h_j)$ displaces $\bar p$ by at most $R_0+R/2$, so $u_{i,j}\in(\Gamma_i/H_i)(\bar p_i,R_0+R)$ for all sufficiently large $i$.

The bound for \(k_2<j\le k\) follows in the same way from $\bar\varphi_i(u_{i,j})\to e_\Gamma$.
\end{proof}

Choose \(r\) sufficiently large that $\Lambda(r)\ge 20C(N)D+R_0+R+1$.
\begin{lemma}
$P_i'$ defines $\Gamma_i/H_i$. Thus $\Gamma_i/H_i$ is naturally isomorphic to $\Gamma_i'$.
\end{lemma}
\begin{proof}
By the choice of $\Lambda(r)$, the set $S_i:=(\Gamma_i/H_i)(\bar p_i,\Lambda(r)-1)$ defines $\Gamma_i/H_i$. The inclusion $\bar\varphi_i^{-1}(B_{\Lambda(r)}(e_\Gamma))\subset P_i'$ gives $S_i\subset P_i'$. By Lemma~\ref{lemma:radius control}, the ordered generators of $P_i'$ belong to $S_i$. The inclusion $S_i\subset P_i'$ therefore induces a surjective homomorphism $\phi_i:\Gamma_i/H_i\longrightarrow\Gamma_i'$.

Conversely, $\bar\varphi_i(P_i')\subset B_{C(r+1)}(e_\Gamma)$. The argument of Lemma~\ref{lemma:radius control} gives $P_i'\subset(\Gamma_i/H_i)(\bar p_i,C(r+1)+R)$. This larger displacement set also defines $\Gamma_i/H_i$. Thus the inclusion of $P_i'$ into this defining set induces a homomorphism $\phi_i':\Gamma_i'\longrightarrow\Gamma_i/H_i$.

The composition $\phi_i'\circ\phi_i:\Gamma_i/H_i\to\Gamma_i/H_i$ is induced by the inclusion of $S_i$ into the larger defining set, so it is the identity. Hence $\phi_i$ is injective as well as surjective, and therefore is an isomorphism.
\end{proof}

We henceforth identify $\Gamma_i/H_i$ with $\Gamma_i'$. In particular, $k=\operatorname{rank}(\Gamma_i/H_i)\le k_0$. If $k=0$, then $\Gamma_i/H_i$ is trivial and $\Gamma_i$
is finite, so the map $\hat X_i\to\{\mathrm{pt}\}$ gives
the first conclusion. In the torsion-free case, $G_i$
is finite and hence trivial, giving the second conclusion.
We therefore assume $k\ge1$.  We can now construct the fiber bundle over a nilmanifold.
\begin{proof}[Proof of Theorem~\ref{thm:first}]
Let $\bar{\psi}_i:\Gamma_i/H_i\to G'/H$ be the group maps in the quotient $\epsilon_i$-eGHAs. By construction, $N_i$ is $C^\infty$-close to $\mathbb R^k$ on bounded sets. For every $g\in(\Gamma_i/H_i)(\bar p_i,1/\epsilon_i)$, its action on $N_i$ is close to the action of $\bar\rho_1(\bar{\psi}_i(g))$ on $\mathbb R^k$.

Equivariant convergence first gives an almost
$\Gamma_i/H_i$-equivariant almost splitting map $f_i'$
from a large ball in $\widetilde X_i/H_i$ into $N_i$.
The extension construction in the proof of Lemma~\ref{thm:gluing GHA} gives a globally defined
almost $\Gamma_i/H_i$-equivariant map, still denoted by $f_i'$,
that retains the local almost splitting estimates.

By Theorem~\ref{thm:inv}, we may choose $f_i'$ to be $\Gamma_i/H_i$-equivariant and locally almost splitting. The Alexandrov lower bound and the assumed $(k_0,\delta(N))$-regularity pass to the covering domains. With $\delta(N)$ chosen sufficiently small, Lemma~\ref{lemma:topological submersion} therefore shows that $f_i'$ is a topological submersion. It descends to a topological submersion
\[
f_i'':
\hat X_i=\widetilde X_i/\Gamma_i
\longrightarrow
N_i/(\Gamma_i/H_i)=N_i/\Gamma_i'.
\]
Since $\hat X_i$ is compact, $f_i''$ is proper and hence is
a topological fiber bundle projection;
see \cite[Essay~II, Chapter~1]{KirbySiebenmann}.

Under the natural identifications of fundamental groups,
the equivariant lift identifies $(f_i'')_*$ with the quotient
homomorphism $\Gamma_i\to\Gamma_i/H_i$.
Its surjectivity implies that the fibers are connected.
Since the nilmanifold base is aspherical, the long exact
sequence of homotopy groups gives $\pi_1(F)\cong H_i$
for each fiber $F$. Thus every fiber has finite fundamental
group.

In the smooth case of Theorem~\ref{theorem:smooth bundle over nil}, the construction in the proof of Theorem~\ref{thm:inv} allows us to choose $f_i'$ smooth. By \cite{HuangHuang2024}, its differential has rank $k$ everywhere. Thus $f_i''$ is a proper smooth submersion
and hence a smooth fiber bundle projection by
Ehresmann's theorem.
\end{proof}

\subsection{Fiber bundles over infranilmanifolds}\label{subsection:4.2}
Keep the notation above and assume that $G_i$ is torsion-free. 
We have shown that $Y=\mathbb R^k\times Z$ and that $\rho_1(G')$ consists of pure translations. 

\begin{lemma}\label{lemma:isom quotient}
After replacing $\Gamma_i$ by a larger nilpotent subgroup of $G_i$ if necessary, the equivariant approximations induce an isomorphism $G_i/\Gamma_i\cong G/G'$.
\end{lemma}
\begin{proof}
Since $|G_i/\Gamma_i|\le C(N)$, we have $|G/G'|\le C(N)$. Fix a sufficiently large $R$ and then take $i$ sufficiently large. Let $G_i'\le G_i$ be generated by the elements that are $\epsilon_i$-close to elements of $G'(\tilde p,R)$. By \cite[Theorem~3.10]{FukayaYamaguchi1992}, the group $G_i'$ is normal in $G_i$, converges equivariantly to $G'$, and satisfies $G_i/G_i'\cong G/G'$, with the isomorphism induced by the equivariant approximations.

By construction, $\Gamma_i\le G_i'$, so
$[G_i:G_i']\le C(N)$.
The map $\rho_1\circ\psi_i:G_i'\to\rho_1(G')$
is a discrete approximation.
The small subgroups $H_i$ are finite and hence trivial
because $G_i$ is torsion-free.
Theorem~\ref{thm:large nil}, together with the defining-set
argument above, therefore identifies $G_i'$ with a lattice
in a $k$-dimensional nilpotent Lie group.
Since $\Gamma_i$ has finite index in $G_i'$, these groups
have the same Malcev completion.
Replacing $\Gamma_i$ by $G_i'$ proves the lemma.
\end{proof}

\begin{proof}[Proof of the second part of Theorem A]
Apply the argument of the preceding subsection. Since $G_i$ is torsion-free, the finite subgroup $H_i$ is trivial. Thus $\Gamma_i$ is identified with a lattice in $N_i$, and $f_i':\widetilde X_i\to N_i$ is a $\Gamma_i$-equivariant locally almost splitting map.

We now extend the translation action of $\Gamma_i$ on $N_i$
to an affine action of $G_i$, following
\cite[proof of Theorem~A]{Wang2024}.
The construction will make $f_i'$ almost $G_i$-equivariant.

Let $J=[G:G']\le C(N)$.
Choose representatives $g_1,\ldots,g_J\in G$, with $g_1=e$,
such that
\[
G=\bigsqcup_{j=1}^J G'g_j.
\]
Write $\rho_1(g_j)=(b_j,A_j)$, where $b_j\in\mathbb R^k$ and $A_j\in O(k)$. For each $j$, choose $g_{ij}\in G_i$ and $b_{ij}\in N_i$ converging to $g_j\in G$ and $b_j\in\mathbb R^k$, respectively. Since $G_i/\Gamma_i\cong G/G'$, we have $G_i=\bigsqcup_{j=1}^J\Gamma_i g_{ij}$.

By Malcev rigidity, conjugation by $g_{ij}$ on the lattice $\Gamma_i$ extends uniquely to an automorphism $C_{g_{ij}}\in\mathrm{Aut}(N_i)$. Under the chosen identifications of $N_i$ with $\mathbb R^k$,
the actions of these automorphisms converge on bounded sets
to the linear actions of $A_j$. In particular, each $C_{g_{ij}}$ is almost isometric for all sufficiently large $i$.
Since $g_1=e$, we choose $g_{i1}=e$, $b_{i1}=e_{N_i}$, and $C_{g_{i1}}=\mathrm{Id}$. 

Define $\phi_i:G_i\to N_i\rtimes\mathrm{Aut}(N_i)$
by the following formula for $h\in\Gamma_i$ and
$1\le j\le J$:
\[
\phi_i(h g_{ij}) := (h b_{ij}, C_{b_{ij}^{-1}} \circ C_{g_{ij}}).
\]
The restriction of $\phi_i$ to $\Gamma_i$ is the original
translation action. On $G_i$, the map is almost multiplicative:
for $1\le j_1,j_2\le J$, the product
$\phi_i(g_{ij_1})\phi_i(g_{ij_2})$ is close to
$\phi_i(g_{ij_1}g_{ij_2})$, with defect of the form
$(g,C_{g^{-1}})$ for some $g\in B_{\epsilon_i}(e_{N_i})$.
These defects arise from the finite quotient $G_i/\Gamma_i$.
Apply the finite-extension averaging construction in
\cite[Lemma~3.1.2]{Dekimpe} to adjust the affine images of
the coset representatives while fixing the translation
action on $\Gamma_i$. We continue to denote the resulting
homomorphism by $\phi_i$.
The correction is bounded by $C\epsilon_i$, since the
nilpotency classes are uniformly bounded and $N_i$ is
$C^\infty$-close to $\mathbb R^k$.

To prove injectivity, let $g\in\ker(\phi_i)$.
Since $\Gamma_i$ has finite index in $G_i$, some positive
power $g^m$ belongs to $\Gamma_i$.
The injectivity of $\phi_i|_{\Gamma_i}$ gives $g^m=e$,
and torsion-freeness then gives $g=e$.
Thus $\phi_i$ is injective.
We henceforth identify $G_i$ with this discrete subgroup
of $N_i\rtimes\mathrm{Aut}(N_i)$.

The action of $\Gamma_i$ on $N_i$ is by translations and hence isometric. By construction, the $G_i$-action on $N_i$ is close, on large balls, to the $G$-action on $\mathbb R^k$. It is therefore almost isometric. Since $|G_i/\Gamma_i|\le C(N)$, we can average the left-invariant metric over $G_i/\Gamma_i$. The resulting metric is still \(C^\infty\)-close to the original one and $G_i$-invariant. We henceforth use this metric, so that $G_i$ acts isometrically on $N_i$.

The map $f_i':\widetilde X_i\to N_i$ is
$\Gamma_i$-equivariant, almost $G_i$-equivariant, and locally
almost splitting.
Averaging over the finite quotient $G_i/\Gamma_i$ makes it
$G_i$-equivariant while preserving the local splitting
estimates.
Lemma~\ref{lemma:topological submersion} then shows that
$f_i'$ is a topological submersion.
As in the first part, it descends to a topological fiber
bundle projection from $X_i=\widetilde X_i/G_i$ onto the
infranilmanifold $N_i/G_i$.

Under the natural identifications of both fundamental groups
with $G_i$, the induced homomorphism is the identity.
Since the base is aspherical, the long exact sequence of
homotopy groups shows that the fibers are connected and
simply connected.

In the smooth case, we use smooth almost splitting maps,
and Ehresmann's theorem gives a smooth fiber bundle projection.
\end{proof}

\section{Nilpotent models}\label{sec 6}
In this section, we construct local orbit-preserving homeomorphisms when the limit space is a metric cone close to $\mathbb{R}^n$ and the limit group is generated by pure translations and the isotropy group at the basepoint. We first determine the local group presentation, then construct its nilpotent model and action, and finally transfer the orbit structure to the original spaces.

Let $Y=\mathbb R^l\times C(Z)$ be an $\mathrm{RCD}(0,n)$
space with measure $\mathcal H^n$, where $C(Z)$ does not
split off a Euclidean line.
Let $p=(0^l,o)$, where $o$ is the cone vertex, and let
$G\le\mathrm{Isom}(Y)$ be closed.
We say that $G$ satisfies property~Q if:
\begin{itemize}\label{prop:Q}
\item $G_0$ is nilpotent.
\item $G=\operatorname{trans}(G)\rtimes G_p$, where $\operatorname{trans}(G)\cong\mathbb R^k$ is a pure translation subgroup.
\end{itemize}
If $G$ satisfies property~Q, then
$G_0=\mathbb R^k\times T$, where $T=(G_p)_0$ is a torus
acting on $\mathbb R^{l-k}\times C(Z)$.
Set $s:=\dim T$.

\begin{theorem}\label{thm:R^n singular fiber bundle}
Let $(Y_i,p_i,d_i,\mathcal H^n)$ be a sequence of $\mathrm{RCD}(-\epsilon_i,n)$ spaces with $\epsilon_i\to0$, and let $G_i\le\mathrm{Isom}(Y_i)$ be closed discrete subgroups acting freely. Suppose that 
\begin{center}
 $\begin{CD}
(Y_i,p_i,G_i,\mathcal{H}^n) @>{\mathrm{eGH}}>> (Y=\mathbb{R}^l \times C(Z),p=(0^l,o),G,\mathcal{H}^n)\\
		@VV\pi V @VV\pi V\\
(Y_i/G_i,\bar{p}_i,\mathcal{H}^n/\mathcal{H}^n(B_1(\bar{p}_i))) @>{\mathrm{mGH}}>> (Y/G,\bar{p},\mathfrak{m}).
\end{CD}$
\end{center}
Assume that $G$ satisfies property Q and that either $Y=\mathbb R^n$, or $l<n$, $C(Z)$ does not split, and $Z$ is $\delta$-close to $S^{n-l-1}$ for some $0<\delta<\delta(n)$.
Then, for all sufficiently large $i$, there exists an orbit-preserving homeomorphism $F_i:B_1(p_i)\to B_1(p)$. It induces a quotient map $\bar F_i:B_{0.99}(\bar p_i)\to B_{0.99}(\bar p)$, whose fibers over $B_{0.98}(\bar p)$ are homeomorphic to infranilmanifolds.  
\end{theorem}

The main difficulty in Theorem~\ref{thm:R^n singular fiber bundle} is that the group maps in the eGHAs need not be homomorphisms, so equivariant averaging does not apply directly.

Nilpotent models provide smooth spaces preserving the discrete actions on local covering spaces. They allow us to glue equivariant almost splitting maps and recover the singular fibers from the resulting orbit quotients. We now give their precise definition. 

\begin{defn}\label{defn:nilpotent models}
We call a quadruple $(M,N,\hat{G},E)$ an $\epsilon$-nilpotent model if the following conditions hold.

\begin{enumerate}
\item
$N$ is a connected nilpotent Lie group equipped with a left-invariant metric. For some integers $0\le k\le n$ and $s\ge0$, we have
\[
\dim N=k+s,
\qquad
\dim T_N=s,
\]
where $T_N$ is the maximal torus of $N$.
There is an $s$-dimensional torus $T\le O(n-k)$ such that exponential coordinates induce an $\epsilon$-$C^\infty$-close diffeomorphism
\[
\varphi:N\longrightarrow G_0=\mathbb{R}^k\times T.
\]
The restriction $\varphi|_{T_N}:T_N\to T$ is an isomorphism, and
\[
\varphi(gt)=\varphi(g)\varphi(t)
\qquad
\text{for all }g\in N,\ t\in T_N.
\]

\item
The model space is the associated bundle
$M=N\times_{T_N}\mathbb{R}^{n-k}$,
defined by the equivalence relation $(gt,v)\sim(g,\varphi(t)v)$.
The group $N$ acts on $M$ by left translation on the first factor.
Define
\[
A=
\left\{
(a,b)\in\operatorname{Aut}(N)\times O(n-k):
\varphi(a(t))=b\varphi(t)b^{-1}
\forall t\in T_N,
\right\}
\]
\[
K=
\left\{
(t^{-1},\mathrm{Id}_N,\varphi(t)):
t\in T_N
\right\}.
\]
The natural action
\[
(u,a,b)\cdot[g,v]=[u\,a(g),bv]
\]
descends to an action of $(N\rtimes A)/K$ on $M$.
The group $\hat{G}$ is a discrete subgroup of $(N\rtimes A)/K$, and $E=N\hat{G}$,
where $N$ is identified with its natural image in $(N\rtimes A)/K$. We require $E$ to be closed in $(N\rtimes A)/K$.

\item
The space $M$ carries a Riemannian metric that is
$\epsilon$-$C^\infty$-close to the Euclidean metric on
$\mathbb{R}^{n}$. Both $N$ and $\hat{G}$ act isometrically on
$M$, and $\hat{G}$ is $\epsilon$-dense in $E$ with respect to
the chosen metric on $E$.
\end{enumerate}
\end{defn}

Fix $R>10$. Let $(M,N,\hat G,E)$ be an $\epsilon$-nilpotent model, and let $(X,p,d,\mathcal{H}^n)$ be an $\mathrm{RCD}(-\epsilon,n)$ space equipped with an isometric $\hat G$-action. Write $(\bar X,\bar p)=(X/\hat G,[p])$ and $p'=[e_N,0^{n-k}]\in M$.
We call a map
\[
f:B_R(\hat Gp)\longrightarrow B_R(\hat Gp')
\]
an $\epsilon$-singular fiber bundle map with model $(M,N,\hat G,E)$ if $f$ is $\hat G$-equivariant, $\epsilon$-injective, and $(\epsilon,n)$-splitting on every ball $B_R(gp)$, $g\in\hat G$.
We always take $\epsilon$ sufficiently small and usually omit it from the terminology.

For sufficiently small $\epsilon$, the local homeomorphism property and \(\epsilon\)-injectivity imply that $f$ is a homeomorphism onto its image. By equivariance, it induces a homeomorphism between corresponding neighborhoods in $X/\hat G$ and $M/\hat G$. Composing with the natural projection to $M/E$ gives the quotient singular fibration $\bar f:B_1(\bar p)\to M/E$:
\[
B_1(\bar p) \to M/\hat G \to M/E.
\]

We should point out that if $G$ is a pure translation group and $Y=\mathbb{R}^n$, our construction for Theorem \ref{thm:R^n singular fiber bundle} agrees with the construction in \cite{Wang2024}.

\subsection{Construction of the acting groups}\label{subsection:same orbits}
In Theorem~\ref{thm:R^n singular fiber bundle}, for any $\epsilon > 0$, $G(p,\epsilon)$ defines $G$. Since $G_i$ converge to $G$, we obtain the following gap phenomenon.

\begin{lemma}\label{lemma:same abstract group}
In the setting of Theorem~\ref{thm:R^n singular fiber bundle}, fix $0<r_1<r_2<\infty$ and set $S_{i,1}=G_i(p_i,r_1)$ and $S_{i,2}=G_i(p_i,r_2)$. Let $\hat{S}_{i,1}$ and $\hat{S}_{i,2}$ be the groups defined by these sets. For all sufficiently large $i$, the inclusion $S_{i,1}\subset S_{i,2}$ induces an isomorphism $\hat{S}_{i,1}\to\hat{S}_{i,2}$.  
\end{lemma}

\begin{proof}
Since $G=\operatorname{trans}(G)\rtimes G_p$, we identify
$Gp$ with $\mathbb R^k$, taking $p$ to $0^k$. Let
$\iota_{i,1}:S_{i,1}\to\hat S_{i,1}$ and $\iota_{i,2}:S_{i,2}\to\hat S_{i,2}$ be the embeddings. Then let
$T_i:\hat S_{i,1}\to\hat S_{i,2}$ be induced by the inclusion $S_{i,1} \subset S_{i,2}$.
We construct an inverse to $T_i$.

Fix $L=10(r_2+1)$. Take the eGHA $\phi_i:G_i \to G$ with $\phi_i(e_{G_i})=e_G$. Then
$\phi_i(G_i(p_i,L))p$ is $\epsilon_i$-dense in $B_{L}(0^k)$.

Choose a fixed integer $m$ such that $r_2/m<r_1/20$.
For each $g\in S_{i,2}$, choose
$g_0=e_{G_i},g_1,\dots,g_m=g$ satisfying
\[
\left|\phi_i(g_l)p-\frac{l}{m}\phi_i(g)p\right|
\le\epsilon_i
\qquad(1 \le l<m).
\]
For all sufficiently large $i$, every
$g_{l-1}^{-1}g_l$ belongs to $G_i(p_i,r_1/10)$. By the choice of $g_l$, we have all $g_l \in S_{i,2}$.
Moreover, if $g\in S_{i,1}$, then all $g_l$ belong to $S_{i,1}$. We call $g_0=e,g_1,\dots,g_m=g$ a chain from $e_{G_i}$ to $g$.

Define $W_i(g):= \iota_{i,1}(g_0^{-1}g_1) \cdots \iota_{i,1}(g_{m-1}^{-1}g_m)$.
The multiplication relations along these chains give
\[
W_i(g)=\iota_{i,1}(g)\quad(g\in S_{i,1}),
\qquad
T_i(W_i(g))=\iota_{i,2}(g)\quad(g\in S_{i,2}).
\]

It remains to show that $W_i$ preserves products defined within $S_{i,2}$.
We verify this by filling the corresponding polygon in the
Euclidean orbit and pulling a fine triangulation back to
the approximating action, as in \cite{Wang2021,Wang2022}.

Take $g,h,gh\in S_{i,2}$. Concatenate the chosen chain from
$e$ to $g$, the left translate by $g$ of the chain from
$e$ to $h$, and the reversed chain from $e$ to $gh$.
This gives a closed polygon in $G_ip_i$ with $3m$ edges,
each shorter than $r_1/10$, whose vertices are represented
by elements of $G_i(p_i,2r_2)$.

Apply $\phi_i$ to the group elements representing the vertices, then
evaluate the resulting elements of $G$ at $p$. The resulting polygon lies in
$B_{2r_2+1}(0^k)$, and its edges have length less than
$r_1/8$ for large $i$. Fill the polygon by coning it to $0^k$, and triangulate
the filling finely enough that every edge has length
less than $r_1/3$.

Using the density of $\phi_i(G_i(p_i,L))p$, choose elements in $G_i(p_i,L)$ for the interior vertices approximating their
Euclidean positions. For large $i$, all lifted edges are shorter
than~$r_1$. Thus, if $a,b,c \in G_i(p_i,L)$
are the vertices of an oriented triangle, its edge
labels belong to $S_{i,1}$ and satisfy
\[
\iota_{i,1}(a^{-1}b)\,
\iota_{i,1}(b^{-1}c)
=
\iota_{i,1}(a^{-1}c).
\]
The triangular identities show that the relation
$g\cdot h=gh$ follows from multiplication relations
within $S_{i,1}$. Hence $W_i(g)W_i(h)=W_i(gh)$.

Therefore $W_i$ extends to a homomorphism
$\hat S_{i,2}\to\hat S_{i,1}$. The two identities above
show that this homomorphism is inverse to $T_i$.
\end{proof}

By Lemma~\ref{lemma:same abstract group}, a diagonal choice gives $R_i\to\infty$ such that, for all sufficiently large $i$, the groups defined by $G_i(p_i,t)$ are naturally isomorphic for every $1/R_i \le t\le6R_i$. We denote this group by $\hat G_i$. 

We next define a gluing space following the construction in \cite{PanWang2021}, 
\[
\hat{Y}_i: =  \hat{G}_i   \times_{G_i(p_i,6R_i)} \bar{B}_{3R_i}(p_i)= ( \hat{G}_i \times  \bar{B}_{3R_i}(p_i)) / \sim,
\]
where $(g_1,x_1) \sim (g_2,x_2)$ in $\hat{G}_i \times  \bar{B}_{3R_i}(p_i)$ if there exists $g \in G_i(p_i,6R_i)$ such that $gx_1=x_2$ and $g_1\iota(g^{-1})=g_2$; $\iota:G_i(p_i,6R_i) \to \hat{G}_i$ is the natural injective local homomorphism.

Equip $\hat Y_i$ with the length metric induced by $\bar B_{3R_i}(p_i)$. Then $\hat G_i$ acts isometrically on $\hat Y_i$. Let $\hat{p}_i$ be the image of $(e_{\hat G_i},p_i)$ in $\hat{Y}_i$. By \cite[Lemma~5.2]{PanWang2021}, the pairs $(B_{R_i}(p_i),G_i(p_i,R_i))$ and $(B_{R_i}(\hat p_i),\hat G_i(\hat p_i,R_i))$ are naturally equivariantly isometric. The quotient spaces $Y_i/G_i$ and $\hat Y_i/\hat G_i$ also agree on the $R_i/3$-balls about their basepoints. Since $R_i\to\infty$, the glued triples have the same eGH limit as the original triples:
\begin{equation}\label{eq:convergence}
\begin{CD}
(\hat{Y}_i,\hat{p}_i,\hat{G}_i,\mathcal{H}^n)
@>{\mathrm{eGH}}>>
(Y,p,G,\mathcal{H}^n)\\
@VV\pi V @VV\pi V\\
(\hat{Y}_i/\hat{G}_i,\bar{p}_i,\mathcal{H}^n/\mathcal{H}^n(B_1(\bar p_i)))
@>{\mathrm{mGH}}>>
(Y/G,\bar{p},\mathfrak{m}).
\end{CD}
\end{equation}

It therefore suffices to prove Theorem~\ref{thm:R^n singular fiber bundle} for $(\hat{Y}_i,\hat{p}_i,\hat{G}_i)$. We replace $G_i$ by $\hat{G}_i$ so that a local subset defines the whole group, allowing us to use a gluing argument similar to that in the proof of Theorem~A. The additional difficulty is the presence of compact isotropy.

Since $\hat G_i$ is defined by $\hat G_i(\hat p_i,1)$, Theorem~\ref{theorem:conv to G_0} gives a normal subgroup $\Gamma_i\triangleleft\hat G_i$ such that $\Gamma_i$ converges to $G_0$ and $\hat G_i/\Gamma_i\cong G/G_0=G_p/T$, where $T$ is the identity component of $G_p$. Moreover, $\Gamma_i$ is generated by $\Gamma_i(\hat p_i,\epsilon_i)$. 
\begin{lemma}\label{lemma:rank of Gamma_i}
For every fixed $\epsilon>0$ and all sufficiently large $i$, $\Gamma_i$ is a torsion-free nilpotent group of rank $k$ defined by $\Gamma_i(\hat p_i,\epsilon)$.
\end{lemma}
\begin{proof}
Since $G_0=\mathbb R^k\times T$, apply Theorem~\ref{thm:large nil} to $\Gamma_i\to G_0\to\mathbb R^k$. The finite small subgroups of $\Gamma_i$ have bounded orbits
and are trivial by Lemma~\ref{lemma:no compact limit},
provided $\delta$ is sufficiently small. Thus $\Gamma_i(\hat p_i,\epsilon)$ is contained in a grid nilprogression $P_i'$ of rank $k$, with $P_i'\subset\Gamma_i(\hat p_i,R)$ for some fixed $R>0$ and all sufficiently large $i$.
Let $\hat{P}_i$ be the groupification of $P_i'$. Then $\hat{P}_i$ is a torsion-free nilpotent group of rank $k$. Since $\Gamma_i$ is generated by $\Gamma_i(\hat p_i,\epsilon_i)$, there is a natural surjective homomorphism $\pi:\hat{P}_i \to \Gamma_i$. Since $\Gamma_i(\hat p_i,\epsilon) \subset P_i'$, there is an injective local homomorphism $\iota: \Gamma_i(\hat p_i,\epsilon) \to \hat{P}_i$.

Let $\{g_j:1\le j\le J\}\subset G_p$ be a set of coset
representatives of $G/G_0=G_p/T$, with $g_1=e_G$, and
choose $g_{i,j}\in\hat G_i$ close to $g_j$, with
$g_{i,1}=e_{\hat G_i}$. Then every $g\in\hat G_i$ has a
unique expression $g=g_{i,j}h$, with $1\le j\le J$ and
$h\in\Gamma_i$.

Let $Q=G_p/T=\hat G_i/\Gamma_i$, and denote the image of
$g\in\hat G_i$ in $Q$ by $\bar g$. Define a section
$t:Q\to\hat G_i$ of the quotient map by
$t_{\bar g_{i,j}}=g_{i,j}$. Define
\[
c:Q\times\hat G_i(\hat p_i,\epsilon/2)
\longrightarrow\Gamma_i(\hat p_i,\epsilon),
\qquad
c(q,s):=t_{\bar s q}^{-1}s\,t_q.
\]
Since each $g_{i,j}$ displaces $\hat p_i$ by at most
$\epsilon_i$, the correction has displacement at most
$\epsilon/2+2\epsilon_i<\epsilon$ for large $i$.
For each $s\in\hat G_i(\hat p_i,\epsilon/2)$, define
\[
F_s:Q\times\hat P_i\longrightarrow Q\times\hat P_i,
\qquad
F_s(q,\hat g)=
\bigl(\bar s q,\iota(c(q,s))\hat g\bigr).
\]
For $s,t,st\in\hat G_i(\hat p_i,\epsilon/2)$, we have
\[
c(\bar t q,s)c(q,t)
=t_{\overline{st}q}^{-1}st\,t_q
=c(q,st).
\]
All three correction terms lie in the domain of $\iota$,
so $F_s\circ F_t=F_{st}$. Also,
$F_e=\mathrm{Id}$ and $F_s^{-1}=F_{s^{-1}}$.
Thus $s\mapsto F_s$ is a local homomorphism into
$\mathrm{Perm}(Q\times\hat P_i)$ with ordinary composition.
Since $\hat G_i(\hat p_i,\epsilon/2)$ defines $\hat G_i$
for all sufficiently large $i$, this map extends uniquely
to a homomorphism on $\hat G_i$.

For any $g\in\Gamma_i$, the extended map $F_g$ preserves
the fiber $\{e_Q\}\times\hat P_i$ and acts on it by left
translation. These translations determine a homomorphism
\[
\psi:\Gamma_i\longrightarrow\hat P_i,
\qquad
F_g(e_Q,\hat g)=(e_Q,\psi(g)\hat g).
\]
For $g\in\Gamma_i(\hat p_i,\epsilon/2)$, we have
$c(e_Q,g)=g$ and hence $\psi(g)=\iota(g)$.
The ordered generators of $P_i'$ lie in
$\Gamma_i(\hat p_i,\epsilon_i)$, so $\psi\circ\pi$
fixes them and is the identity on $\hat P_i$.
Therefore $\pi$ is an isomorphism, and $\Gamma_i$
is torsion-free nilpotent of rank $k$.

Finally, the argument in the proof of
Lemma~\ref{lemma:same abstract group}, applied to the
restricted $\Gamma_i$-actions, shows that
$\Gamma_i(\hat p_i,\epsilon)$ defines $\Gamma_i$.
\end{proof}

\subsection{Nilpotent model spaces}\label{sec:nil model space}
In this subsection, we use $\Gamma_i \le \hat{G}_i$ to construct a connected nilpotent group that converges to $G_0=\mathbb{R}^k \times T$ in the $C^\infty$ sense. 

Applying the structure of approximate groups to the eGHA $\psi_i:\Gamma_i \to G_0$, we conclude that there are $\delta>0$ and a nilprogression $P_i \subset \Gamma_i$ of rank $k+s$ so that $\psi_i^{-1}(B_{2\delta}(e_{G})) \subset P_i$ and $\psi_i(P_i) \subset B_{1}(e_{G})$. Here $P_i$ is obtained from an approximation to a sufficiently small neighborhood of $e_{G_0}$, whereas $P_i'$ in Lemma~\ref{lemma:rank of Gamma_i} is obtained from an approximation to $\mathbb R^k$ on large scales.

By Malcev's theorem, for all sufficiently large $i$, $P_i$
embeds in a lattice $\tilde\Gamma_i$ in a connected,
simply connected nilpotent Lie group $\tilde N_i$
of dimension $k+s$.
Lemma~\ref{lemma:metric} and Theorem~\ref{thm:Q conv}
provide left-invariant metrics for which $\tilde N_i$
converges in $C^\infty$ on bounded sets to $\mathbb R^{k+s}$.
After choosing local coordinates on $T$, the comparison
maps in exponential coordinates agree asymptotically
with $\psi_i$ on $P_i$.

By Theorem~\ref{theorem:conv to G_0}, $\Gamma_i$ is generated
by an $\epsilon_i$-ball centered at the identity for some
$\epsilon_i\to0$. Therefore the inclusion $P_i\subset\Gamma_i$
induces a surjective homomorphism
$\pi:\tilde\Gamma_i\to\Gamma_i$. Set $D_i=\ker(\pi)$. Since $\psi_i^{-1}(B_{2\delta}(e_{G})) \subset P_i$ and $\pi$ is injective on $P_i \subset \tilde{\Gamma}_i$, the distance between any two distinct elements of $D_i$ is at least $\delta$. After passing to a subsequence, we obtain \[ \tilde\psi_i:(\tilde\Gamma_i,D_i) \xrightarrow{\mathrm{GH}}(\mathbb R^{k+s},D), \] where $D$ is a discrete subgroup of $\mathbb R^{k+s}$. The maps $\tilde\psi_i$ restrict to local isomorphisms between large balls in $D_i$ and $D$. The quotients therefore satisfy \[ \Gamma_i=\tilde\Gamma_i/D_i \longrightarrow\mathbb R^{k+s}/D. \] Thus $\mathbb R^{k+s}/D=\mathbb R^k\times T=G_0$.

\begin{lemma}\label{lemma:same kernel}
For all sufficiently large $i$, $D_i$ is central in
$\tilde N_i$, and the local group approximations
$D_i\to D$ extend to a group isomorphism.
\end{lemma}
\begin{proof}
By Lemma~\ref{lemma:rank of Gamma_i}, for small $\epsilon>0$ and large $i$, $\Gamma_i$ is a torsion-free nilpotent group of rank \(k\) and is defined by $\Gamma_i(\hat p_i,\epsilon)$. Since $\tilde{\Gamma}_i$ has rank $k+s$, we conclude that $D_i$ has rank $s$.

Choose a basis $\{u_1,\dots,u_s\}$ of $D$ and elements $u_{i,j}\in D_i$ such that $u_{i,j}\to u_j$ for $1\le j\le s$. We first show that, for each $1 \le j \le s$, $u_{i,j}$ is in the center of $\tilde N_i$. $D_i=\ker(\pi)$ is a normal subgroup of $\tilde{\Gamma}_i$. Since $\tilde{\Gamma}_i$ is $\epsilon_i$-dense in $\tilde{N}_i$, it is generated by elements close to $e_{\tilde N_i}$. If $g \in \tilde{\Gamma}_i$ is close to the identity, then $gu_{i,j}g^{-1} \in D_i$ is close to $u_{i,j}$. Any two distinct elements of $D_i$ are at distance at least some $\delta>0$. Thus $gu_{i,j}g^{-1}=u_{i,j}$. The elements of $\{u_{i,1},\dots,u_{i,s}\}$ lie in the center of $\tilde\Gamma_i$ and hence also in the center of $\tilde{N}_i$.

We next show that the $u_{i,j}$ are linearly independent for all sufficiently large $i$. Since they are central in $\tilde N_i$, we work in additive exponential coordinates. Suppose that there are real coefficients $a_{i,j}$ satisfying \[ \sum_{j=1}^s a_{i,j}^2=1, \qquad \sum_{j=1}^s a_{i,j}u_{i,j}=0. \] After passing to a subsequence, assume that $a_{i,j}\to a_j$ for every $j$. Theorem~\ref{thm:Q conv} then gives \[ \sum_{j=1}^s a_ju_j=0^{k+s}, \qquad \sum_{j=1}^s a_j^2=1. \] This contradicts the fact that $\{u_1,\dots,u_s\}$ is a basis of $D$.

Finally, we show that $D_i$ is generated by $\{u_{i,1},\dots,u_{i,s}\}$ for all sufficiently large $i$. Assume $g_i \in D_i$ does not belong to the subgroup generated by $\{u_{i,1},\dots,u_{i,s}\}$. Since $\operatorname{rank}(D_i)=s$, there exists a nonzero $t \in \mathbb{Z}$ such that $g_i^t=\prod_{j=1}^s u_{i,j}^{a_j}$ for some $a_j \in \mathbb{Z}$. In particular, $g_i$ is in the center of $\tilde N_i$ and $g_i=\sum_{j=1}^s \frac{a_j}{t}u_{i,j}$. After adding suitable integer multiples of the $u_{i,j}$ to $g_i$ in these coordinates, we may assume $|\frac{a_j}{t}| \le 1$. In particular, passing to a subsequence, $g_i$ converges to $g \in D$. Then $g=\sum_{j=1}^s b_ju_{j}$ for some $b_j \in \mathbb{Z}$. Then $g_i=\sum_{j=1}^s b_ju_{i,j}$ by uniform discreteness, a contradiction.

Sending $u_{i,j}$ to $u_j$ therefore defines an isomorphism
$D_i\to D$ compatible with the local approximations.
\end{proof}

Now we define $N_i=\tilde{N}_i/D_i$ with the induced left-invariant metric. Since $\Gamma_i \cong \tilde{\Gamma}_i/D_i$, we can identify $\Gamma_i$ with a discrete subgroup of $N_i$ by a natural injective homomorphism $\phi_i:\Gamma_i \to N_i$. Let $T_i$ be the maximal torus of $N_i$.
\begin{corollary}\label{cor:homo along torus}
For all sufficiently large $i$, there is a $C^\infty$-close
diffeomorphism $\varphi_i:N_i\to G_0$ satisfying
$\varphi_i(T_i)=T$. Moreover, for every $g\in N_i$ and $t\in T_i$,
\[
\varphi_i(gt)=\varphi_i(g)\varphi_i(t),
\qquad
\varphi_i(t^{-1})=\varphi_i(t)^{-1}.
\] 
\end{corollary}
\begin{proof}
Since $D_i$ is central in $\tilde N_i$, its Malcev completion
$N_i'\le\tilde N_i$ is also central. Let $N' \le \mathbb{R}^{k+s}$ be the subspace generated by $D$. Recall that there is a $C^{\infty}$-close diffeomorphism $\tilde{\varphi}_i:\tilde{N}_i\to\mathbb{R}^{k+s}$ induced by a linear identification of their Lie algebras as vector spaces. By Lemma~\ref{lemma:same kernel}, we may slightly modify this diffeomorphism so that it maps $N_i'$ isomorphically onto $N'$ and maps $D_i$ to $D$. Moreover, since $N_i'$ is in the center of $\tilde{N}_i$, for any $g \in N_i'$ and $h \in \tilde{N}_i$, $\tilde{\varphi}_i(gh)=\tilde{\varphi}_i(g)+\tilde{\varphi}_i(h)$.

Then $\tilde{\varphi}_i$ descends to a \(C^\infty\)-close diffeomorphism $\varphi_i:N_i \to \mathbb{R}^{k+s}/D =G_0$. Since $T_i=N_i'/D_i$, we have $\varphi_i(T_i)=T$. Moreover, for every $g\in N_i$ and $t\in T_i$, we have $\varphi_i(gt)=\varphi_i(g)\varphi_i(t)$.
\end{proof}

We next construct a $G$-equivariant almost splitting map on $Y$.
\begin{lemma}\label{lemma:limit R^n}
The group $G_p$ can be identified with a subgroup of $O(n)$ through an injective homomorphism $\iota:G_p\to O(n)$. Moreover, there is a $G$-equivariant map $f:B_2(Gp)\to\mathbb R^n$ that is $\Phi(\delta\mid n)$-almost splitting on every
ball $B_2(gp)$, $g\in G$.
\end{lemma}

\begin{proof}
Since $C(Z)$ does not split, $G_p \le O(l) \times \mathrm{Isom}(Z)$. Since $Z$ is sufficiently close to $S^{n-l-1}$, there is an injective homomorphism $\iota_Z:\mathrm{Isom}(Z) \to O(n-l)$. Together with the identity on \(O(l)\), this induces an injective homomorphism $\iota:G_p \to O(l) \times O(n-l)\le O(n)$.

Let $H$ be the image of $G_p$ under the projection onto
$\mathrm{Isom}(Z)$. By averaging, we obtain an $H$-equivariant
$\Phi(\delta\mid n)$-almost splitting map
$f_z:B_2(o)\to\mathbb R^{n-l}$. Since $G$ satisfies property~Q,
taking the product with the identity on the $\mathbb R^l$-factor
gives a $G$-equivariant map $f:B_2(Gp)\to\mathbb R^n$
that is almost splitting on each $B_2(gp)$.
\end{proof}

Now we construct the nilpotent model space.
Since $G_0$ is nilpotent, we have $G_0=\operatorname{trans}(G) \times T =\mathbb{R}^k \times T$ with $\iota(T) \le O(n-k)$. 
We can identify 
\[
\mathbb{R}^n =G_0 \times_T \mathbb{R}^{n-k} =  G_0 \times \mathbb{R}^{n-k} / \sim,
\] 
by $(gt,x) \sim (g,\iota(t)x)$ for any $g \in G_0, t\in T$ and $x \in \mathbb{R}^{n-k}$. To simplify notation, we identify $T$ with $\iota(T)$ and write $T\le O(n-k)$.  

We now define the nilpotent model space as the associated bundle
\[
M_i:=N_i \times_{T_i} \mathbb{R}^{n-k} = N_i \times \mathbb{R}^{n-k} / \sim,
\]
where $(g,x) \sim (h,y)$ if there exists $t \in T_i$ so that $gt=h$ and $\varphi_i(t)y=x$. Equip $M_i$ with the connection metric induced by the
left-invariant metric on $N_i$ and the Euclidean metric
on $\mathbb R^{n-k}$.
Its horizontal metric is inherited from $N_i/T_i$,
and its vertical metric is Euclidean.
The diffeomorphism $\varphi_i:N_i\to G_0$ induces a
$C^\infty$-close diffeomorphism
\[
M_i\longrightarrow\mathbb R^n,
\qquad [g,x]\longmapsto[\varphi_i(g),x].
\]
The group $N_i$ acts isometrically by left translation
on the first factor.

Fix $R>0$. For all sufficiently large $i$ and every $g_i\in\Gamma_i(\hat p_i,R)$, the action of $g_i$ on $\hat Y_i$ is close to that of $\varphi_i(g_i)$ on $\mathbb R^n$. The latter action is also close to that of $g_i$ on $M_i$.
Thus the $\Gamma_i$-actions on $\hat Y_i$ and $M_i$ are equivariantly close on bounded balls.
In the next subsection, we extend the action on $M_i$
to $\hat G_i$ while retaining this local approximation.

\subsection{Group actions on nilpotent models}\label{sec:nil group}
The preceding subsection identifies $\Gamma_i$ with a subgroup of $N_i$ through the natural injective homomorphism $\phi_i:\Gamma_i\to N_i$. We now study the relation between $\hat{G}_i$ and $N_i$.  For $g_i\in\hat G_i$, let $C_{g_i}$ denote conjugation by $g_i$ on $\Gamma_i$.
\begin{lemma}\label{lem:conjugation extend}
For all sufficiently large $i$, $C_{g_i}$ extends naturally to an automorphism of $\tilde N_i$ preserving $D_i$. It therefore induces an automorphism of $N_i$, also denoted by $C_{g_i}$.
\end{lemma}
\begin{proof}
Since $\hat{G}_i$ is generated by $\hat{G}_i(\hat{p}_i,\epsilon_i)$ for some $\epsilon_i \to 0$, we may assume $g_i \in \hat{G}_i(\hat{p}_i,\epsilon_i)$.
 
For a small $\epsilon$, conjugation by $g \in G(p,\epsilon)$ maps a small neighborhood of $e_{G}$ to a small neighborhood of $e_G$. Thus conjugation by $g_i \in \hat G_i(\hat p_i,\epsilon_i)$ maps a small neighborhood of the identity to a small neighborhood as well. In particular, conjugation by $g_i$ induces a homomorphism $C_{g_i}:\tilde\Gamma_i\to\tilde\Gamma_i$. Its inverse is induced by conjugation by $g_i^{-1}$, so $C_{g_i}$ is an automorphism. Malcev rigidity extends it to an automorphism of $\tilde N_i$. Since $\Gamma_i$ is normal in $\hat G_i$, conjugation by $g_i$ also induces an automorphism of $\Gamma_i$. By construction, the projection $\pi:\tilde\Gamma_i\to\Gamma_i$ intertwines these two automorphisms. Hence $C_{g_i}(D_i)=D_i$.
\end{proof}

We now extend the translation action of $\Gamma_i$
on $M_i$ to an action of $\hat G_i$.
We first determine which pairs consisting of an
automorphism of $N_i$ and an orthogonal transformation
of $\mathbb R^{n-k}$ induce diffeomorphisms of $M_i$.

Let $A_i:=\{ (a,b) \in \mathrm{Aut}(N_i) \times O(n-k) | \varphi_i(a(t))=b\varphi_i(t)b^{-1}, \forall t \in T_i \}$.
\begin{lemma}
The set $A_i$ is a group acting on $M_i$ by diffeomorphisms. 
\end{lemma}

\begin{proof}
For $[g,x]\in M_i=N_i\times_{T_i}\mathbb R^{n-k}$ and $(a,b)\in A_i$, define $(a,b):[g,x]\mapsto[a(g),bx]$. To check that this is well-defined, note that $[g,x]=[gt,\varphi_i(t^{-1})x]$ for every $t\in T_i$. Then 
\[
[a(gt),b\varphi_i(t^{-1})x]=[a(g)a(t),b\varphi_i(t^{-1})b^{-1}b(x)]=[a(g),bx], 
\]
since $\varphi_i(a(t))= b\varphi_i(t)b^{-1}$ and $\varphi_i$ is a homomorphism on $T_i$.
Thus the $A_i$-action is well-defined on $M_i$.

For $(a,b)\in A_i$, the inverse action is induced by $(a^{-1},b^{-1})\in A_i$, so $(a,b)$ acts by a diffeomorphism. If $(a_1,b_1),(a_2,b_2)\in A_i$, then their natural product $(a_1a_2,b_1b_2)$ also lies in $A_i$. Therefore $A_i$ is a group. 
\end{proof}

Define the semidirect product $N_i \rtimes A_i$ as follows. For any $g_1,g_2 \in N_i$ and $(a_1,b_1),(a_2,b_2) \in A_i$, let $(g_1,a_1,b_1)(g_2,a_2,b_2)=(g_1a_1(g_2),a_1a_2,b_1b_2)$.
We next define the action of $N_i \rtimes A_i$ on $M_i$. For any $(g,a,b) \in N_i \rtimes A_i$ and $[h,x] \in M_i=N_i \times_{T_i} \mathbb{R}^{n-k}$,
\[
(g,a,b):[h,x]\mapsto[ga(h),bx].
\]
Then $N_i \rtimes A_i$ acts by diffeomorphisms on $M_i$.

Define $K_i \le N_i \rtimes A_i$ by $K_i:=\{(t^{-1},\mathrm{Id}_{N_i},\varphi_i(t)) | t \in T_i\}$. 
\begin{lemma}
The subgroup $K_i$ is the kernel of the action of $N_i\rtimes A_i$ on $M_i$. 
\end{lemma}
\begin{proof}
Every element of $K_i$ acts trivially on $M_i$.
Conversely, suppose that $(g,a,b)\in N_i\rtimes A_i$
acts trivially. Fix $h\in N_i$. For every
$x\in\mathbb R^{n-k}$, the equality
$[ga(h),bx]=[h,x]$ gives $t_h\in T_i$ such that
\[
ga(h)t_h=h,\qquad bx=\varphi_i(t_h)x.
\]
The first equation determines $t_h$ uniquely,
independently of $x$. Hence $b=\varphi_i(t_h)$.
Since $\varphi_i|_{T_i}$ is injective, $t_h=t$
is also independent of $h$.
Taking $h=e_{N_i}$ gives $g=t^{-1}$.
Centrality of $T_i$ then gives $a(h)=h$ for every
$h\in N_i$. Thus
$(g,a,b)=(t^{-1},\mathrm{Id}_{N_i},\varphi_i(t))\in K_i$.
\end{proof}
For any $(g,a,b) \in N_i \rtimes A_i$, we use the notation $[g,a,b] \in (N_i \rtimes A_i)/K_i$ for its quotient image.

The natural inclusions
$\Gamma_i\le N_i\le(N_i\rtimes A_i)/K_i$ are injective.
We next extend their composition to an injective homomorphism
from $\hat G_i$ to $(N_i\rtimes A_i)/K_i$.

For all large $i$, we have $\hat{G}_i/\Gamma_i \cong G/G_0 \cong G_p/T$ which is a finite group of order $J$. Choose $\{g_j\}_{j=1}^J \subset G_p$ to be a set of coset representatives of $ G_p/T$ with $g_1=e_{G}$. For each $j$ and large $i$, choose $g_{i,j} \in \hat{G}_i$ close to $g_j \in G_p$. Then $\hat{G}_i$ is generated by $\Gamma_i$ and $\{g_{i,j} | 1 \le j \le J\}$.

Notice that $\operatorname{trans}(G)=\mathbb{R}^k$ is $G$-invariant. Thus $G_p \le O(k) \times O(n-k)$. For each $j$, assume $g_j=(b_{1j},b_{2j}) \in O(k) \times O(n-k)$. Therefore for each $t \in T \le O(n-k)$, $g_jtg_j^{-1}=b_{2j}tb_{2j}^{-1}$. 

Fix a large $i$. Using Lemma~\ref{lem:conjugation extend}, we denote the conjugation action of $g_{i,j} \in\hat G_i$ on $N_i$ by $C_{i,j}$.
Identify $T_i$ with $T$ using $\varphi_i$.
Since $g_{i,j}$ is eGH-close to $g_j$, their conjugation actions on this torus are close.
The injectivity radii of $T$ and $T_i$ have a common positive lower bound for all sufficiently large $i$.
These conjugation actions therefore coincide.
Equivalently, $\varphi_i(C_{i,j}(t))=b_{2j}\varphi_i(t)b_{2j}^{-1}$ for every $t\in T_i$.
In particular, from the definition of $A_i$, $(C_{i,j},b_{2j}) \in A_i$.
\begin{lemma}\label{lemma:auto conv}
The action of $(C_{i,j},b_{2j})$ on $M_i$ converges to the action of $g_j$ on $\mathbb{R}^n$.
\end{lemma}
\begin{proof}
Recall that $\mathbb R^n=G_0\times_T\mathbb R^{n-k}$ and that $N_i$ converges in $C^\infty$ to the first factor $G_0$.

Let $C_j$ denote conjugation by $g_j$ on $G_0$.
Under the chosen identifications, the automorphisms
$C_{i,j}$ converge to $C_j$.
It therefore suffices to show that the $(C_j,b_{2j})$-action on $(\mathbb R^k\times T^s)\times_{T^s}\mathbb R^{n-k}$ agrees with the $g_j$-action on $\mathbb R^n$.

$g_j=(b_{1j},b_{2j}) \in O(k) \times O(n-k)$. Then the conjugation $C_j$ on $G_0$ is $(b_{1j}, C_{b_{2j}})$, where $b_{1j}$ acts orthogonally on the \(\mathbb R^k\)-factor of $G_0$ and $C_{b_{2j}}$ acts on $T$. A point $(x,y) \in \mathbb{R}^k \times \mathbb{R}^{n-k}$ corresponds to $(x,e_{T},y) \in (\mathbb{R}^{k} \times T) \times_{T} \mathbb{R}^{n-k}$. Then
\[
(C_j,b_{2j})(x,e_{T},y) = (b_{1j}(x),e_{T},b_{2j}(y)) \in (\mathbb{R}^{k} \times T) \times_{T} \mathbb{R}^{n-k}
\]
is exactly $g_j(x,y)$, as required.
\end{proof}

We next construct $\phi_i:\hat{G}_i \to (N_i \rtimes A_i)/K_i$ which extends the natural inclusion $\Gamma_i \le N_i$.
\begin{lemma}\label{lem:embed G_i}
We can extend $\Gamma_i\le N_i\le (N_i\rtimes A_i)/K_i$ to an injective homomorphism $\phi_i:\hat G_i\to (N_i\rtimes A_i)/K_i$. Moreover, for every fixed $R>0$ and all sufficiently large $i$, the action of $g\in\hat G_i(\hat p_i,R)$ on $B_R(\hat p_i)$ is close to the action of $\phi_i(g)$ on $M_i$.
\end{lemma}
\begin{proof}
Every $g\in\hat G_i$ has a unique expression
$g=hg_{i,j}$ with $h\in\Gamma_i$.
Define
\[
\phi_i(g)=[\phi_i(h),C_{i,j},b_{2j}].
\]
For each fixed $R>0$, the corresponding actions are close
on $R$-balls whenever $g\in\hat G_i(\hat p_i,R)$. Therefore $\phi_i$ is an almost homomorphism. 

Since $\phi_i$ is already a homomorphism on $\Gamma_i$,  the error comes from the finite group $\hat G_i/\Gamma_i\cong G/G_0$. For any $g_1,g_2,g_1g_2\in\hat G_i(\hat p_i,R)$, 
$\phi_i(g_1)\phi_i(g_2)\phi_i(g_1g_2)^{-1}=[g,C_{g^{-1}},\mathrm{Id}_{\mathbb R^{n-k}}] \in  (N_i \rtimes A_i)/K_i$ for some $ g \in B_{\epsilon_i}(e_{N_i})$. Since $N_i$ has a positive lower bound on its injectivity radius, we may identify $B_{\epsilon_i}(e_{N_i})$ with an $\epsilon_i$-ball in the universal cover $\tilde{N}_i$. Then we may apply the averaging technique to modify $\phi_i$ to be an exact homomorphism; see also subsection \ref{subsection:4.2}. $\phi_i$ is injective because $M_i$ and $\hat Y_i$ are non-collapsing.
\end{proof}
\begin{remark}
The homomorphism $\phi_i$ maps $g_{i,j}$ to $[h_{i,j},C_{h_{i,j}}^{-1}C_{i,j},b_{2j}]$, where $h_{i,j} \in N_i$ is the correction term which is $\epsilon_i$-close to the identity. The factor $C_{h_{i,j}}^{-1}$ ensures that conjugation by $g_{i,j}$ on $\Gamma_i$ corresponds, under $\phi_i$, to conjugation by $\phi_i(g_{i,j})$ on $\phi_i(\Gamma_i)$. Since $h_{i,j} \in N_i$,  $C_{h_{i,j}}^{-1}$ is trivial on $T_i$. Thus the correction leaves the slice component $b_{2j}$ unchanged.
\end{remark}
Using $\phi_i$, we may simply identify $\hat{G}_i$ with a subgroup of $(N_i \rtimes A_i)/K_i$. This subgroup acts on $M_i$ by almost isometries. Averaging the metric over $\hat{G}_i/\Gamma_i$ if necessary, we may assume $\hat{G}_i$ acts isometrically on $M_i$. 

Let $E_i\le (N_i\rtimes A_i)/K_i$ be generated by $N_i$ and
$\hat G_i$. Then $(M_i,N_i,\hat G_i,E_i)$ is an
$\epsilon_i$-nilpotent model.

\subsection{Singular fibrations from nilpotent models}\label{section:gluing}
Let $p_i'=[e_{N_i},0^{n-k}] \in M_i$ and consider the isotropy $(E_i)_{p_i'}$. Then we have $T_i \le (E_i)_{p_i'}$ and $E_i=N_i(E_i)_{p_i'} = N_i \times_{T_i} (E_i)_{p_i'}$. Moreover, by the proof of Lemma~\ref{lemma:auto conv}, $(E_i)_{p_i'}$ converges to $G_p$. By the argument of Lemma~\ref{lemma:homomorphic eGHA},
the isomorphism $\varphi_i|_{T_i}:T_i\to T$ extends to
a group isomorphism
$\varphi_i':(E_i)_{p_i'}\to G_p$
that is also a Gromov--Hausdorff approximation.

We extend $\varphi_i:N_i\to G_0$ to a map
$\varphi_i'':E_i\to G$ as follows.
For $g=g_1g_2$, with $g_1\in N_i$ and
$g_2\in(E_i)_{p_i'}$, set
$\varphi_i''(g)=\varphi_i(g_1)\varphi_i'(g_2)$.
\begin{lemma}\label{lemma:extend diff}
The map $\varphi_i''$ is a well-defined $C^\infty$ diffeomorphism. 
\end{lemma}
\begin{proof}
By Corollary~\ref{cor:homo along torus}, $\varphi_i$
intertwines the right torus actions, and $\varphi_i'$
intertwines the corresponding left actions. Hence the
product diffeomorphism $\varphi_i\times\varphi_i'$
descends through
\[
E_i=N_i\times_{T_i}(E_i)_{p_i'},
\qquad
G=G_0\times_T G_p
\]
to the map $\varphi_i''$. Its inverse descends in the
same way, so $\varphi_i''$ is a diffeomorphism.
Its $C^\infty$ closeness follows from that of the two factors.
\end{proof}
By Lemma~\ref{lemma:extend diff}, we henceforth denote this extension by $\varphi_i:E_i \to G$.
By Lemma~\ref{lemma:limit R^n}, we may identify $G_p$ with a subgroup of $O(n)$. We next show that $\varphi_i:E_i\to G$ induces a diffeomorphism from $M_i$ to $\mathbb R^n$ preserving the full orbit relations.

\begin{lemma}\label{lem:orbit-preserving}
For all sufficiently large $i$, there is a natural
diffeomorphism $J_i:M_i\to\mathbb R^n$ such that $J_i(E_i x)=GJ_i(x)$ for any $x\in M_i$.
\end{lemma}
\begin{proof}
We define $J_i:M_i = N_i \times_{T_i} \mathbb{R}^{n-k} \to G_0 \times_{T} \mathbb{R}^{n-k} = \mathbb{R}^n$ by
\[
J_i([g,x])=[\varphi_i(g),x],
\]
for $g \in N_i$ and $x \in \mathbb{R}^{n-k}$. This is well defined since, for $t\in T_i$,
\[
J_i([gt,\varphi_i(t^{-1})x])
 =[\varphi_i(g)\varphi_i(t),\varphi_i(t^{-1})x]
 =[\varphi_i(g),x].
\]
Since $\varphi_i$ is a $C^{\infty}$-close diffeomorphism, $J_i$ is also a $C^{\infty}$-close diffeomorphism.   

By definition, $J_i$ maps $N_i$-orbits to $G_0$-orbits, and $J_i^{-1}$ maps $G_0$-orbits to $N_i$-orbits. Let $\{g_j\}_{j=1}^J\subset G_p$ be a set of coset representatives of $G_p/T$, with $g_1=e_G$, and choose $g_{i,j}\in\hat G_i$ converging to $g_j$. Write $g_j=(b_{1j},b_{2j})\in O(k)\times O(n-k)$. Then $\phi_i(g_{i,j})=[h_{i,j},C_{h_{i,j}}^{-1}C_{i,j},b_{2j}]$, where $h_{i,j}\in N_i$ is $\epsilon_i$-close to the identity; see the remark after Lemma~\ref{lem:embed G_i}.

Then for any $1\le j\le J$ and $[g,x]\in M_i$, 
\[
J_i(\phi_i(g_{i,j})[g,x])
 =J_i([C_{i,j}(g)h_{i,j},b_{2j}x])
 =[\varphi_i(C_{i,j}(g)h_{i,j}),b_{2j}x].
\]
The action of $g_j$ on $\mathbb{R}^n=G_0 \times_T \mathbb{R}^{n-k}$ can be expressed as $(e_{G_0},C_j,b_{2j})$. Then 
\[
g_jJ_i([g,x])=[C_j(\varphi_i(g)),b_{2j}x].
\]  
The two displayed points have the same slice component, so they lie on the same $G_0$-orbit. Thus $J_i(\phi_i(g_{i,j})[g,x])$ and $J_i([g,x])$ lie on the same $G$-orbit for every $j$. Together with the $N_i$-orbit correspondence, this shows that $J_i$ maps each $E_i$-orbit into a $G$-orbit.

Similarly, we can show that $J_i^{-1}$ maps every $G$-orbit to an $E_i$-orbit.
\end{proof}

We now construct a map $f_i$ and its quotient $\bar{f}_i$ fitting into the diagram
\begin{equation}\label{eq:singular fiber}
\begin{CD}
(\hat{Y}_i,\hat{p}_i,\hat{G}_i) @>{f_i}>> (M_i,p_i',E_i)\\
@V{\pi}VV @VV{\pi}V\\
(\hat{Y}_i/\hat{G}_i,\bar{p}_i) @>{\bar{f}_i}>> (M_i/E_i,\bar{p}_i').
\end{CD}
\end{equation}
Under the identifications $M_i\to\mathbb R^n$, both the
$\hat G_i$- and $E_i$-actions converge to the induced
$G$-action on $\mathbb R^n$. Consequently, $\hat G_i$
is $\epsilon_i$-dense in $E_i$.

Since $\Gamma_i$ is nilpotent of rank $k$ and $|\hat{G}_i/\Gamma_i|$ is bounded, we can apply Theorem \ref{thm:inv} to obtain the following equivariant almost splitting maps.
\begin{lemma}\label{lem:G_i equi}
Fix $R>10$. For sufficiently small $\delta$ and all sufficiently large $i$, there is an $\epsilon_i$-injective and $\hat G_i$-equivariant map
\[
f_i:B_R(\hat G_i\hat p_i)\longrightarrow B_R(\hat G_i p_i')\subset M_i,
\]
which is $(\Phi(\delta\mid n)+\epsilon_i)$-almost splitting on any $B_R(g\hat{p}_i)$, $g \in \hat{G}_i$. In particular, $f_i$ is a homeomorphism.
\end{lemma}

The quotient of $f_i$ is the singular fibration
\[
\bar{f}_i:B_R(\bar p_i)\longrightarrow B_R(\bar p_i')\subset M_i/E_i,
\]
where $\bar p_i$ and $\bar p_i'$ are the quotient basepoints. It factors through the homeomorphism induced by $f_i$ on the quotients by $\hat G_i$, followed by the projection $M_i/\hat G_i\to M_i/E_i$. 

For fixed $R>10$ and all sufficiently large $i$, the group-metric convention
in Section~2.1 and the $\epsilon_i$-density of $\hat G_i$ in $E_i$ imply that,
for every $x\in B_{0.99R}(p_i')$, the orbit quotient $E_i(x)/\hat G_i$ lies
in an $\epsilon_i$-ball in the model quotient $M_i/\hat G_i$.

\begin{lemma}\label{lem:infranil fiber}
Let $\bar x\in B_{0.98R}(\bar p_i')\subset M_i/E_i$
and choose a representative $x\in B_{0.99R}(p_i')$.
Then the fiber $\bar f_i^{-1}(\bar x)$ is homeomorphic
to $E_i(x)/\hat G_i$, an infranilmanifold of dimension
$k+s-\dim((E_i)_x)$, where $(E_i)_x$ is the isotropy
subgroup at $x$.
\end{lemma}
\begin{proof}
For large $i$, the orbit $E_i(x)$ lies in $B_R(\hat{G}_ip_i')$. Under the homeomorphism induced by $f_i$, the fiber is identified with $E_i(x)/\hat G_i$. Since $E_i=\hat{G}_iN_i$ and $N_i/\Gamma_i$ is compact, the fiber is compact and connected. The action of $\hat G_i$ on the model neighborhood is free; hence the fiber dimension is $\dim(E_i)-\dim((E_i)_x)=k+s-\dim((E_i)_x)$.

To show that $E_i(x)/\hat G_i$ is an infranilmanifold, it suffices to consider the connected component of $E_i(x)$ containing $x$, since the fiber is connected. This component is $N_i(x)$, because $E_i/N_i=G_p/T$ is finite and $N_i$ is connected. Let $H_i=\{g \in \hat{G}_i | g(N_i(x))=N_i(x)\}$. Then $E_i(x)/\hat{G}_i = N_i(x)/H_i$.

Recall that $N_i=\tilde N_i/D_i$, with $\tilde N_i$ simply connected and $D_i$ central. Let $\pi:\tilde N_i\to N_i$ be the quotient map, set $S_i=(N_i)_x\le T_i$, and let $\tilde S_i=\pi^{-1}(S_i)$ be its full preimage. The closed subgroup $S_i$ need not be connected. The subgroup $\tilde S_i$ is central, contains $D_i$, and lies in the Malcev completion of $D_i$. Thus $N_i(x)=\tilde N_i/\tilde S_i$.

To study the $H_i$-action on $N_i(x)$, we use the notation of Lemma~\ref{lem:embed G_i} and the remark following it. Since $\Gamma_i \le H_i$, it suffices to determine which representatives $g_{i,j}$ preserve $N_i(x)$.
Suppose that $x=[g,y] \in M_i=N_i \times_{T_i} \mathbb{R}^{n-k}$. Then $g_{i,j}(x)=[C_{i,j}(g)h_{i,j},b_{2j}y]$.
$g_{i,j}(x) \in N_i(x)$ if and only if there exists $t_{i,j} \in T_i$ such that 
\begin{equation}\label{equation:exists t}
\varphi_i(t_{i,j})^{-1} b_{2j}y=y.
\end{equation}
In this case, $g_{i,j}(x)=[C_{i,j}(g)h_{i,j}t_{i,j},y]$. Fix $j$ for which such a $t_{i,j}\in T_i$ exists.
The action of $g_{i,j}$ on $N_i(x)$ is induced by the affine transformation $(h_{i,j}t_{i,j},C_{h_{i,j}^{-1}}C_{i,j})\in N_i\rtimes\mathrm{Aut}(N_i)$.
This transformation lifts to $\tilde g_{i,j}\in\tilde N_i\rtimes\mathrm{Aut}(\tilde N_i)$.
The lift is unique up to $D_i$-translations.
Its automorphism part maps $D_i$ to itself.

For $j$ satisfying \eqref{equation:exists t}, the points $x$ and $g_{i,j}x$ lie on the same $N_i$-orbit and have the same closed stabilizer $S_i$. Hence the automorphism part of $\tilde g_{i,j}$ preserves its full preimage $\tilde S_i$.

Let $\tilde{S}_{i,0}$ be the identity component of $\tilde{S}_i$ and define $\bar{N}_i=\tilde{N}_i/\tilde{S}_{i,0}$. Let $\bar D_i$ be the image of $\tilde S_i$ in $\bar N_i$. $\bar{N}_i$ is simply connected and $\bar{D}_i$ is a discrete subgroup in the center of $\bar{N}_i$. Thus $\tilde g_{i,j}$ induces an affine transformation $\bar g_{i,j}\in\bar N_i\rtimes\mathrm{Aut}(\bar N_i)$.
Its automorphism part maps $\bar D_i$ to itself. Let $\bar{\Gamma}_i$ be the image of $\tilde{\Gamma}_i$ in $\bar N_i$. Define 
\[
\bar{G}_i=\langle \bar{\Gamma}_i,\bar{D}_i, \bar{g}_{i,j} \text{ where }\eqref{equation:exists t}\text{ holds} \rangle.
\]
Thus $N_i(x)/H_i$ is homeomorphic to $\bar N_i/\bar G_i$.
The group $\bar G_i$ acts freely and discretely on
$\bar N_i$.
Its image in $\mathrm{Aut}(\bar N_i)$, the holonomy group
of this presentation, is finite; indeed,
\[
[\bar G_i:\bar G_i\cap\bar N_i]
\le
\bigl|\bar G_i/\langle\bar\Gamma_i,\bar D_i\rangle\bigr|
\le J=|G_p/T|.
\]
Therefore the fiber is homeomorphic to an infranilmanifold.
\end{proof}

Now we can finish the proof of Theorem~\ref{thm:R^n singular fiber bundle}.
\begin{proof}[Proof of Theorem~\ref{thm:R^n singular fiber bundle}]
Fix $R>10$ and carry out the preceding construction on the corresponding larger neighborhoods. By Subsection~\ref{subsection:same orbits}, the triples $(Y_i,p_i,G_i)$ and $(\hat Y_i,\hat p_i,\hat G_i)$ are equivariantly isometric on neighborhoods of their basepoints. Lemma~\ref{lem:G_i equi} gives the map $f_i$ to the nilpotent model $(M_i,E_i)$, and Lemma~\ref{lem:infranil fiber} identifies its quotient fibers over $B_{0.98R}(\bar p_i')$ as infranilmanifolds. Lemmas~\ref{lem:orbit-preserving} and~\ref{lemma:limit R^n} give homeomorphisms between neighborhoods of the basepoints
\[
(M_i,E_i,p_i')\longrightarrow(\mathbb R^n,G,0^n)\longrightarrow(Y,G,p)
\]
that map $E_i$-orbits bijectively to $G$-orbits. Composing these maps with $f_i$ gives the orbit-preserving homeomorphism $F_i$, with quotient map $\bar F_i$.

We restrict $F_i$ to $B_1(p_i)$ and $\bar F_i$ to $B_{0.99}(\bar p_i)$. Choosing $\delta$ sufficiently small and then $i$ sufficiently large ensures that every fiber over a point
of $B_{0.98}(\bar p)$ lies in $B_{0.99}(\bar p_i)$.
\end{proof}

\begin{corollary}\label{cor:C(n)-nilpotent}
In Theorem~\ref{thm:R^n singular fiber bundle}, suppose
in addition that each $G_i$ contains a nilpotent subgroup
of index at most $C$ and nilpotency step at most $c$.
Then the constructed infranilmanifold fibers admit
presentations whose holonomy groups have order at most $C$.
\end{corollary}
\begin{proof}
By assumption, each $\hat G_i$ contains a nilpotent subgroup
$G_i'$ of index at most $C$ and nilpotency step at most $c$.
After passing to a subsequence, these subgroups converge
to a nilpotent subgroup $G'\le G$ of index at most $C$.
Set $G_p':=G'\cap G_p$. Then
$G'=\operatorname{trans}(G)\rtimes G_p'$.
Since $G'$ is nilpotent, $G_p'$ acts trivially on the
translation factor and hence lies in $O(n-k)$.
The subgroup induced by \(G_i'\) in each fiber presentation has trivial holonomy for all sufficiently large \(i\).
Consequently, the holonomy group of the fiber presentation
in Lemma~\ref{lem:infranil fiber} has order at most
$[\hat G_i:G_i']\le C$.
\end{proof}

\begin{remark}
Corollary~\ref{cor:C(n)-nilpotent} and the generalized
Margulis lemma~\ref{thm:Margulis2} imply that the
infranilmanifold fibers in Theorems~B and~C admit
presentations whose holonomy groups have order at most
$C(n)$.
\end{remark}
\subsection{Injectivity radius in local universal covers}
The collapsing theory of \cite{CFG} gives a uniform
positive injectivity radius bound at the basepoints of
local universal covers of complete manifolds with
two-sided sectional curvature bounds.
This implies the local rewinding Reifenberg condition.
We give an independent proof using the local group
structure.

\begin{theorem}\label{thm:bounded sec}
There exists $\rho=\rho(n)>0$ such that the following holds. Let $N$ be a complete $n$-manifold with $|\operatorname{sec}|\le1$, and let $p\in N$. If $\tilde p\in\tilde{B}(p,\rho)$ is a lift of $p$, then $\operatorname{inj}(\tilde p)\ge\rho/3$.
\end{theorem}
\begin{proof}
We argue by contradiction. Suppose that there are complete pointed $n$-manifolds $(N_i,p_i)$ with $|\operatorname{sec}_{N_i}|\le1$, radii $r_i\to0$, and lifts $\tilde p_i\in\tilde B(p_i,r_i)$ such that $\operatorname{inj}(\tilde p_i)<r_i/3$. After a slow rescaling and relabeling, we may also assume that $|\operatorname{sec}_{N_i}|\le\epsilon_i\to0$; for example, replace $(N_i,r_i)$ by $(r_i^{-1/2}N_i,r_i^{1/2})$.

We construct a pseudocover and its local group of pseudodeck transformations as in \cite{BuserKarcher}. Choose $R_i\to\infty$ with $R_i<\epsilon_i^{-1/4}$. The exponential map
\[
\exp_i:B_{R_i}(0^n)\subset T_{p_i}N_i\longrightarrow N_i
\]
has nondegenerate differential. Equip $B_{R_i}(0^n)$ with the pullback metric and denote the resulting space by $N_i'$, with origin $q_i$. Radial geodesics give $\operatorname{inj}_{N_i'}(q_i)\ge R_i$.

For each $v\in B_{R_i/10}(q_i)$ with $\exp_i(v)=p_i$, there is a unique local isometry
\[
g:B_{R_i/10}(q_i)\longrightarrow B_{R_i/5}(q_i)
\]
satisfying $g(q_i)=v$ and $\exp_i\circ g=\exp_i$. Let $G_i$ be the resulting local group. Products are understood on common domains and extended by unique continuation. In particular, if $g,h\in G_i(q_i,R_i/30)$, then $gh\in G_i(q_i,R_i/10)$ is well-defined.

After passing to a subsequence, we have
\[
\begin{CD}
(N_i',q_i,G_i) @>{\mathrm{eGH}}>> (\mathbb R^n,q=0^n,G)\\
@V{\pi}VV @VV{\pi}V\\
(N_i,p_i) @>{\mathrm{GH}}>> (X,p).
\end{CD}
\]
The quotient identification is local: $(N_i'/G_i,\bar q_i)$ and $(N_i,p_i)$ are isometric on the $R_i/100$-balls about their basepoints. After a further slow rescaling and relabelling, we may assume that $G$ satisfies property Q. Let $T$ be the maximal torus of $G_0$, and write $G_0/T=\mathbb R^k$.

Fix $R>0$ and set $S_i=G_i(q_i,6R)$. Let $\hat G_i$ be the abstract group defined by $S_i$, and form the gluing space
\[
\hat N_i=\hat G_i\times_{S_i}\overline B_{3R}(q_i),
\qquad
\hat q_i=[e_{\hat G_i},q_i].
\]
We show that the natural map $\iota_i:S_i\to\hat G_i$ is injective by embedding $S_i$ into a genuine group.

Set $Q=G/G_0=G_q/T$. Choose a fixed comparison radius $L>30R$, large enough for the local constructions below. For all sufficiently large $i$, separation of the finitely many components of $G$ gives a well-defined component map
\[
\beta_i:G_i(q_i,L)\longrightarrow Q.
\]
Equivariant convergence implies that $\beta_i$ preserves every local product used below. Let $\Gamma_i$ denote its local kernel, consisting of elements close to $G_0$. The following construction is similar to that in the proof of Lemma~\ref{lemma:rank of Gamma_i}. Applying the structure of approximate groups at large scales to $\Gamma_i \to G_0=\mathbb{R}^k \times T \overset{\pi}\to \mathbb{R}^k$, we obtain a nilprogression $P_i\subset\Gamma_i$ of rank $k$ containing $\Gamma_i(q_i,10R)$. Let $\hat P_i$ be the abstract group defined by $P_i$. For large $i$, the natural local homomorphism $\iota:\Gamma_i(q_i,10R) \longrightarrow\hat P_i$ is injective.

Choose representatives $\{g_j:1\le j\le J\}\subset G_q$
for $Q$, with $g_1=e_G$, and choose $g_{i,j}\in G_i$
converging to $g_j$, with $g_{i,1}=e_{G_i}$. Write
$\bar g_j$ for the image of $g_j$ in $Q$, and define a section
\[
t:Q\longrightarrow S_i,
\qquad t_{\bar g_j}=g_{i,j}.
\]
Then $t_{e_Q}=e_{G_i}$, $\beta_i(t_a)=a$, and each $t_a$
displaces $q_i$ by at most $\epsilon_i$.

Define
\[
c:Q\times S_i\longrightarrow\Gamma_i(q_i,10R),
\qquad
c(a,s):=t_{\beta_i(s)a}^{-1}s\,t_a.
\]
The correction has displacement at most
$6R+2\epsilon_i<10R$ for large $i$.
For $s\in S_i$, define
\[
F_s:Q\times\hat P_i\longrightarrow Q\times\hat P_i,
\qquad
F_s(a,\hat g)=
\bigl(\beta_i(s)a,\iota(c(a,s))\hat g\bigr).
\]
If $s,u,su\in S_i$, then
\[
c(\beta_i(u)a,s)c(a,u)
=t_{\beta_i(su)a}^{-1}su\,t_a
=c(a,su).
\]
All three correction terms lie in
$\Gamma_i(q_i,10R)\subset P_i$, so this relation also
holds after applying $\iota$. Consequently,
\[
F_s\circ F_u=F_{su}.
\]
In particular, $F_e=\mathrm{Id}$ and
$F_s^{-1}=F_{s^{-1}}$, so each $F_s$ is a bijection.
With ordinary composition of permutations, the map
\[
F:S_i\longrightarrow\mathrm{Perm}(Q\times\hat P_i),
\qquad F(s)=F_s,
\]
is therefore a local group homomorphism.

To check injectivity, observe that
\[
F_s(e_Q,e_{\hat P_i})
=\bigl(\beta_i(s),\iota(t_{\beta_i(s)}^{-1}s)\bigr).
\]
Equality of these values for $s,u\in S_i$ gives
$\beta_i(s)=\beta_i(u)$ and then, by injectivity of
$\iota$ and left cancellation, $s=u$. Thus $F$ is
injective. By the defining property of $\hat G_i$,
the local homomorphism $F$ extends to a homomorphism
on $\hat G_i$. It follows that
$\iota_i:S_i\to\hat G_i$ is injective.

The gluing construction now identifies $(\hat N_i,\hat q_i,\hat G_i)$ equivariantly with $(N_i',q_i,G_i)$ on the $R$-balls about their basepoints. In particular, $\operatorname{inj}_{\hat N_i}(\hat q_i)\ge R$.

Let $q_i'=\pi(\hat q_i)\in\hat N_i/\hat G_i$. Then $B_{R/3}(p_i)$ is isometric to $B_{R/3}(q_i')$. The local pseudodeck actions are free and have finite displacement sets on compact subballs, so the glued action is free and properly discontinuous on the relevant neighborhood. For large $i$, let $B_i$ be the connected component of $\pi^{-1}(B_{r_i}(q_i'))$ containing $\hat q_i$. The restricted projection makes $B_i$ a covering space of $B_{r_i}(q_i')$. Since $B_{r_i}(\hat q_i)\subset B_i$ and $\operatorname{inj}_{\hat N_i}(\hat q_i)\ge R$, we have $\operatorname{inj}_{B_i}(\hat q_i)\ge r_i/3$. 
The universal cover of $B_i$ is also the universal cover of $B_{r_i}(q_i')\cong B_{r_i}(p_i)$. Passing to this cover cannot decrease the injectivity radius at a lift of $\hat q_i$, contradicting the choice of $r_i$.
\end{proof}

\section{Local singular fibrations}



\subsection{Group structure of equivariant tangent cones}\label{sec:6.1}
Before constructing a local singular fibration in Theorem B, we need a stability theorem for the equivariant tangent cone. 

For any metric space $(X,p,G)$ with $G \le \mathrm{Isom}(X)$ and $r,t>0$, write $B_t^r(p)$ for the $t$-ball centered at $p$ in $rX=(X,rd)$ and $G^r$ for the rescaled action of $G$ on $rX$.
\begin{theorem}\label{theorem:stability of group}
Assume $(X,p,d,\mathcal{H}^n)$ is an $\mathrm{RCD}(-(n-1),n)$ space and satisfies the $(\rho,\delta)$ Reifenberg condition at $p$ for some $\rho>0$ and small $0<\delta<\delta(n)$. Let $G$ be a closed subgroup of $\mathrm{Isom}(X)$, and let $(Y,q,H)$ be any equivariant tangent cone of $(X,p,G)$ at $p$, obtained by rescaling along a sequence $r_i\to\infty$. Then
\[
H=\operatorname{trans}(H)\rtimes H_q,
\qquad \operatorname{trans}(H)\cong\mathbb R^k,
\qquad k=\dim(G/G_p).
\]
Moreover, after passing to a subsequence of the blow-up sequence if necessary, there is an isomorphism $\iota:G_p\to H_q$ such that, for every $g\in G_p$, the rescaled action of $g \in G^{r_i}$ is $\epsilon_i$-close to $\iota(g) \in H_q$, where $\epsilon_i$ is independent of the choice of $g$.
\end{theorem}
  
Since $Y$ is $\mathrm{RCD}(0,n)$, the splitting theorem holds. By \cite[Proposition~2.8]{DeGigli2019}, $Y$ is a metric cone. Write $Y=\mathbb{R}^l \times C(Z)$, where $C(Z)$ does not split. 

Fix $i_0$ sufficiently large that $(r_{i_0}X,p)$ is $\delta$-close to $\mathbb{R}^n$. If $\delta$ is small enough, we can identify $G_p^{r_{i_0}}$ with a subgroup of $O(n)$ via an injective homomorphism $\iota_{i_0}:G_p^{r_{i_0}} \to O(n)$. By the proof of Lemma \ref{lemma:limit R^n}, there is a $G_p^{r_{i_0}}$-equivariant $\Phi(\delta\mid n)$-almost splitting map
\[
f:B_{10}^{r_{i_0}}(p)\longrightarrow B_{10}(0^n)
\]
with $f(p)=0^n$.

We first construct a continuous injective homomorphism $\iota:G_p\to H_q$ and then prove its surjectivity in Corollary~\ref{cor:iso}.
\begin{lemma}\label{lemma:iota}
There is a continuous injective homomorphism $\iota:G_p\to H_q$
such that, after passing to a subsequence, the action of every $g\in G_p$,
viewed as an element of $G_p^{r_i}$, is $\epsilon_i$-close to $\iota(g)$
for some sequence $\epsilon_i\to0$ independent of $g$.
\end{lemma}

\begin{proof}
We first construct a $G_p$-equivariant transformation map for $f$.
For any $0<s<1/2$, we recall the construction in the transformation theorem. Consider the matrix 
\[
M_{ab}^s= \fint_{B_s^{r_{i_0}}(p)}
\langle\nabla f^a,\nabla f^b\rangle,
\]
and let $T_s$ be the lower triangular matrix with positive diagonal entries such that $T_sM^sT_s^T=\mathrm{I}$. Equivariance of $f$ implies that $M^s$ commutes with
$\iota_{i_0}(G_p)$, although $T_s$ need not.
Set $A_s=(M^s)^{-1/2}$. Then $A_sM^sA_s^T=\mathrm I$,
and $A_s$ also commutes with $\iota_{i_0}(G_p)$. Then $A_sT_s^{-1} \in O(n)$. Apply the transformation theorem also at scale $2s$. After rescaling, $s^{-1}T_{2s}\circ f$ is $\Phi(\delta\mid n)$-almost splitting on $B^{s^{-1}r_{i_0}}_2(p)$. Since $\|T_sT_{2s}^{-1}-\mathrm{I}\|\le\Phi(\delta\mid n)$ and $A_sT_s^{-1}$ is orthogonal, the same holds for $s^{-1}A_s\circ f$. Thus the canonical Reifenberg theorem applies on $B^{s^{-1}r_{i_0}}_1(p)$.

By the canonical Reifenberg theorem, there is a small
$\beta>0$, depending on $i_0$ and $\delta$, such that for all
$g_1,g_2\in G_p$ and $x\in B^{s^{-1}r_{i_0}}_1(p)$, we have
\begin{equation}\label{equation:holdfer-Lip}
\begin{aligned}
C^{-1}d_{s^{-1}r_{i_0}X}(g_1x,g_2x)^{1+\beta}
&\le |s^{-1}A_s\circ f(g_1x)-s^{-1}A_s\circ f(g_2x)|\\
&\le Cd_{s^{-1}r_{i_0}X}(g_1x,g_2x).
\end{aligned}
\end{equation}
Since $s^{-1} A_s \circ f$ is also $G_p$-equivariant with respect to $\iota_{i_0}$, 
\[
s^{-1} A_s \circ f (g_1x) - s^{-1} A_s \circ f (g_2x) = s^{-1} (\iota_{i_0}(g_1)-\iota_{i_0}(g_2)) A_s \circ f (x).
\]
Since $|s^{-1}A_s \circ f (x)| \le 2$, by \eqref{equation:holdfer-Lip} 
\begin{equation}\label{equation1}
\|\iota_{i_0}(g_1)-\iota_{i_0}(g_2)\| \ge \frac{1}{2C} d_{s^{-1}r_{i_0}X}(g_1x,g_2x)^{1+\beta}.
\end{equation}
On the other hand, we can choose a unit vector $v \in \mathbb{R}^n$ maximizing $|(\iota_{i_0}(g_1)-\iota_{i_0}(g_2))v|$ and then choose $x_0 \in B_s^{r_{i_0}}(p)$ such that $s^{-1}A_s \circ f (x_0)$ is close to $v$. By \eqref{equation:holdfer-Lip}, 
\begin{equation}\label{equation2}
\|\iota_{i_0}(g_1)-\iota_{i_0}(g_2)\| \le Cd_{s^{-1}r_{i_0}X}(g_1x_0,g_2x_0).
\end{equation}
Combining \eqref{equation1} and \eqref{equation2}, we obtain 
\begin{equation}\label{equation3}
\begin{aligned}
C^{-1}\|\iota_{i_0}(g_1)-\iota_{i_0}(g_2)\|
&\le \sup_{x\in B^{s^{-1}r_{i_0}}_1(p)}
 d_{s^{-1}r_{i_0}X}(g_1x,g_2x)\\
&\le C\|\iota_{i_0}(g_1)-\iota_{i_0}(g_2)\|^{\frac{1}{1+\beta}}.
\end{aligned}
\end{equation}

Take a countable dense subset $\{g_j\}$ of $G_p^{r_{i_0}}$. By a diagonal argument, pass to one subsequence along which $g_j\in G_p^{r_i}$ converges to $h_j\in H_q$ for every $j$. Define $\iota(g_j)=h_j$. By the right-hand side of \eqref{equation3}, $\iota$ extends to a continuous homomorphism on $G_p^{r_{i_0}}$. By the left-hand side of \eqref{equation3},  $\iota$ is injective. By the construction of $\iota$ and \eqref{equation3}, for every $g\in G_p$, the rescaled action of $g$ converges to $\iota(g)\in H_q$. The uniform convergence is due to \eqref{equation3} and the compactness of $G_p$.
\end{proof}

Next we study $f$ along the orbit $Gp$. We take a small neighborhood $U$ of $e \in G^{r_{i_0}}$ such that, for every $g \in U$, the map $x\mapsto f(g^{-1}x)$ is $\delta$-$H^{1,2}$-close to $f$ on $B_{9}^{r_{i_0}}(p)$. We can take a left Haar measure $dg$ on $G^{r_{i_0}}$ and a smooth non-negative function $\chi$ compactly supported in $U$ with $\int_{G^{r_{i_0}}} \chi(g) dg = 1$. Then the averaging map
\[
F_0(x)=\int_{G^{r_{i_0}}} \chi(g)f(g^{-1}x) dg
\] 
is $H^{1,2}$-close to $f$ on $B_{8}^{r_{i_0}}(p)$.
Since the group acts by measure-preserving isometries,
the averaging preserves the Lipschitz and Laplacian bounds
of $f$. The $H^{1,2}$-closeness controls the averaged
gradient errors. Thus $F_0$ is also a
$\Phi(\delta\mid n)$-almost splitting map.

By left invariance of Haar measure, for every $h\in U$,
\[
F_0(hx)= \int_{G^{r_{i_0}}} \chi(hg)f(g^{-1}x) dg.
\]
Thus $F_0(hx)$ is smooth in $h$. 
In particular, the induced local orbit map,
\[ 
\hat{F}_0:G^{r_{i_0}}/G_p^{r_{i_0}} \to \mathbb{R}^n,\hat{F}_0(gG_p^{r_{i_0}})=F_0(gp),
\] 
is smooth in a small neighborhood of $eG_p^{r_{i_0}}$. Since $F_0$ is almost splitting on $B_{8}^{r_{i_0}}(p)$, it is a homeomorphism. Then $\hat{F}_0$ is also injective in the small neighborhood of $eG_p^{r_{i_0}}$. Hence there exists $g_0\in G^{r_{i_0}}$, close to the identity, such that $d\hat F_0$ has full rank at $g_0G_p^{r_{i_0}}$.

Define 
\[
F:B_{7}^{r_{i_0}}(p) \to \mathbb{R}^n, F(x)=F_0(g_0x) -F_0(g_0p).
\]
Since $g_0$ is close to $e\in G^{r_{i_0}}$, $F$ is still almost splitting.
Define the induced local orbit map, 
\[
\hat{F}:G^{r_{i_0}}/G_p^{r_{i_0}} \to \mathbb{R}^n,\hat{F}(gG_p^{r_{i_0}})=F(gp).
\]
Since $d\hat{F}_0$ has rank $k$ at $g_0G_p^{r_{i_0}}$, $d\hat{F}$ has rank $k$ at $eG_p^{r_{i_0}}$.

Let $V=B_1(0^k) \subset \mathbb{R}^k$, and use $V$ as a coordinate neighborhood of $eG_p^{r_{i_0}} \in G^{r_{i_0}}/G_p^{r_{i_0}}$. Take a smooth section $a:V \to  G^{r_{i_0}}$ with $a(0^k)=e$. Since $d\hat F$ has rank $k$ at $eG_p^{r_{i_0}}$, for $v\in V$ sufficiently close to $0^k$ we have
\begin{equation}\label{equation:A}
F(a(v)p)=Av + O(|v|^2),
\end{equation}
where $A$ is an $n \times k$ matrix with rank $k$.

Since $F$ is Lipschitz, $|F(a(v)p)| =|F(a(v)p) - F(p)| \le Cd_{r_{i_0}X}(a(v)p,p)$. By \eqref{equation:A}, if $v$ is close to $0^k$, $|F(a(v)p)| \ge |Av|/2$. Since $A$ has rank $k$, $|Av| \ge c|v|$ for some $c > 0$. Then we conclude that 
\begin{equation}\label{equation:|v|}
|v| \le C' d_{r_{i_0}X}(a(v)p,p)
\end{equation} 
for some $C'>0$ and any $v \in V' \subset V$, where $V'$ is a small neighborhood of $0^k$.

We will show that every limit in $H$ of rescaled elements of $a(V')$ is a pure translation. The following estimate controls these actions across scales.
\begin{lemma}\label{lemma:control F}
There exists $C''>0$ such that for any $v \in V'$ and $x \in B_1^{r_{i_0}}(p)$, $|F(a(v)x)-F(x)-F(a(v)p)| \le C'' |v| d_{r_{i_0}X}(x,p)$.
\end{lemma}
\begin{proof}
By the definition of $F$, the left-hand side is 
\[
\begin{aligned}
&|F_0(g_0a(v)x)-F_0(g_0x)-F_0(g_0a(v)p)+F_0(g_0p)|\\
&\qquad=\left|\int_{G^{r_{i_0}}}
 (\chi(g_0a(v)g)-\chi(g_0g))
 (f(g^{-1}x)-f(g^{-1}p))\,dg\right|.
\end{aligned}
\]

By smoothness and compact support of $\chi$, after shrinking $V'$ we can choose $C'>0$ such that $|\chi(g_0 a(v) g)-\chi(g_0g)| \le C'|v|$ for all $v \in V'$ and $g\in G^{r_{i_0}}$. Since $\chi$ is supported in $U$, the measure of the set $S:=\{g\in G^{r_{i_0}}:\chi(g_0a(v)g)-\chi(g_0g)\ne0\}$ is bounded by $D$, uniformly in $v\in V'$. Then $\int_{G^{r_{i_0}}} |\chi(g_0 a(v) g)-\chi(g_0g)| dg \le DC'|v|$.

Since $f$ is Lipschitz, $|f(g^{-1}x)-f(g^{-1}p)| \le Cd_{r_{i_0}X}(x,p)$. Taking $C''=DC'C$, we finish the proof.
\end{proof}

Let $s_i=\frac{r_{i_0}}{r_i}$. Then for sufficiently large $i > i_0$, $G^{r_i}(p,1)=G^{r_{i_0}}(p,s_i) \subset a(V')G_p$. By the transformation theorem, there are normalization matrices $T_s$ such that $\|T_{s_i}\| \le Cs_i^{-\theta}$ for some $0<\theta(n,\delta)<1$ and $\|T^{-1}_{s_i}\|\le C$. For large $i$, the map $f_i=s_i^{-1}T_{s_i}F$ is almost splitting on $B_5^{r_i}(p)$. For every fixed $R\ge1$ and all sufficiently large $i$, the following estimate holds for every $v\in\mathbb R^n$:
\[
C^{-1}R^{-\theta}|v| \le |T_{s_i}T_{Rs_{i}}^{-1}v| \le CR^{\theta}|v|.
\]
After passing to a subsequence if necessary, $f_i$ converges to a map $f_Y:Y\to\mathbb R^n$ with $f_Y(q)=0$.
For every $R\ge1$ and every $x,y\in B_R(q)$, this map satisfies
\begin{equation}\label{equation:f_Y}
C^{-1}R^{-\beta-\theta}d(x,y)^{1+\beta} \le |f_Y(x)-f_Y(y)| \le CR^{\theta}d(x,y).
\end{equation} 
Thus $f_Y$ is injective. For any $x\in Y$ with $d(x,q)=R\ge1$, we have $|f_Y(x)| \ge C^{-1}d(x,q)^{1-\theta}$. Then $f_Y$ is proper.

Notice that for any $R>0$ and sufficiently large $i$, $G^{r_i}(p,R)$ is contained in $a(V')G_p$. We use $f_Y$ and Lemma~\ref{lemma:control F} to identify the limit of the rescaled actions $a(v_i)$ on $Y$.
\begin{lemma}\label{lemma:pure trans}
If there exist $v_i \in V'$ such that $a(v_i) \in G^{r_i}$ converges to $h \in H$, then $h \in \operatorname{trans}(H)$.
\end{lemma}
\begin{proof} 
Since $a(v_i) \in G^{r_i}$ converges to $h \in H$, $d_{r_{i}X}(a(v_i)p,p) \le R$ for some fixed $R>0$ and all large $i$.
Take any $y \in Y$ and $x_i \in r_iX$ converging to $y$. Then $d_{r_iX}(x_i,p) \le R'$ for some fixed $R'>0$ and all large $i$. By \eqref{equation:|v|}, $|v_i| \le C'd_{r_{i_0}X}(a(v_i)p,p) \le C'Rs_i$, $s_i=r_{i_0}/r_i$. By Lemma~\ref{lemma:control F}, 
\[
|F(a(v_i)x_i)-F(x_i)-F(a(v_i)p)| \le C'' |v_i| d_{r_{i_0}X}(x_i,p) \le Cs_i^2
\]
for all large $i$ and some $C$ independent of $i$.

Recall $f_i=s_i^{-1}T_{s_i}F$. Since $\|T_{s_i}\|\le Cs_i^{-\theta}$ for some fixed small $\theta>0$ and all sufficiently large $i$, we obtain 
\[
|f_i(a(v_i)x_i)-f_i(x_i)-f_i(a(v_i)p)| \le Cs_i^{2-1-\theta} \to 0
\]
as $i\to\infty$. It follows that 
\begin{equation}\label{equation:f_Y equ}
f_Y(hy)=f_Y(y)+f_Y(hq),
\end{equation}
for any $y \in Y$.

Assume $Y=\mathbb{R}^l \times C(Z)$ and $q=(0^l,o)$, where $C(Z)$ does not split. Then $hq=(u,o)$ for some $u \in \mathbb{R}^l$. For $R>d(hq,q)$, the triangle inequality gives $h(B_R(q))\subset B_{2R}(q)$. Applying \eqref{equation:f_Y} to $y\in B_R(q)$ and $hy$, we obtain
\[
C^{-1}R^{-\beta-\theta}d(hy,y)^{1+\beta} \le |f_Y(hy)-f_Y(y)|=|f_Y(hq)|.
\]
Thus 
\begin{equation}\label{equation:gamma}
\sup_{y \in B_R(q)}  d(hy,y) \le CR^{\gamma}, \gamma = \frac{\theta+\beta}{1+\beta} < 1,
\end{equation}
for some $C>0$ and all large $R$. In particular, $h$ has sublinear displacement growth as $\gamma<1$.

Assume $h=(A,u,g)$ where $A \in O(l)$, $u \in \mathbb{R}^l$, and $g \in \mathrm{Isom}(Z)$. If $A$ or $g$ is nontrivial, $h$ cannot have sublinear displacement growth, a contradiction. Thus $h$ is a pure translation.
\end{proof}

\begin{corollary}\label{cor:iso}
The map $\iota:G_p\to H_q$ is an isomorphism. 
\end{corollary}
\begin{proof}
Take $h\in H_q$ and choose $g_i\in G^{r_i}$ converging to $h$. For some $R>0$ and all sufficiently large $i$, we have $g_i\in G^{r_i}(p,R)=G^{r_{i_0}}(p,Rs_i)\subset a(V')G_p$. We can therefore write $g_i=a(v_i)k_i$ with $v_i\in V'$ and $k_i\in G_p$.

After passing to a subsequence, $k_i\to k_0$ in $G_p$. By the uniform convergence in Lemma~\ref{lemma:iota}, the rescaled actions of $k_i$ converge to $\iota(k_0)\in H_q$. Hence $a(v_i)\in G^{r_i}$ converges to $h\iota(k_0)^{-1}\in\operatorname{trans}(H)$ by Lemma~\ref{lemma:pure trans}. Since $h \in H_q$, $h\iota(k_0)^{-1} \in \operatorname{trans}(H) \cap H_q$. Therefore $\iota(k_0)=h$, proving surjectivity.
\end{proof}

\begin{corollary}\label{coro:gene}
$H=\operatorname{trans}(H) \rtimes H_q$. 
\end{corollary}
\begin{proof}
The decomposition used in the proof of
Corollary~\ref{cor:iso} applies to any $h\in H$,
giving $h=tg$ with $t\in\operatorname{trans}(H)$
and $g\in H_q$. Thus $H=\operatorname{trans}(H)H_q$.
Since the pure translation subgroup is normal and
intersects $H_q$ trivially, this product is semidirect.
\end{proof}

We now identify the pure translations in $H$ and determine their dimension.
For any $u \in \mathbb{R}^l$, let $T_u \in \mathrm{Isom}(Y)$ be the pure translation action along $u$. Define
\[
S:=\{ u \in \mathbb{R}^l | T_u \in H \}.
\]
Recall $q=(0^l,o) \in Y=\mathbb{R}^l \times C(Z)$. Define
\[
\psi:S \to \mathbb{R}^n, \psi(u)= f_Y(u,o)=f_Y(T_u q).
\]  
Then by the proof of Corollary~\ref{coro:gene}, any pure translation in $H$ can be approximated by some $a(v_i) \in G^{r_i}$. Therefore we can apply \eqref{equation:f_Y equ} and conclude that $\psi$ is additive. Since $f_Y$ is injective and continuous, $\psi$ is also injective and continuous.

We finish the proof of Theorem~\ref{theorem:stability of group} by the following lemma.
\begin{lemma}\label{lemma:S}
$S \cong \mathbb{R}^k$.
\end{lemma}
\begin{proof}
Fix $u\in S$. By the proof of Corollary~\ref{coro:gene},
choose $a(v_i)\in G^{r_i}$ converging to $T_u$.
Their basepoint displacements are bounded, so
\eqref{equation:|v|} gives $|v_i|\le Cs_i$ with $C$
independent of $i$. Combining \eqref{equation:A} with
$\|T_{s_i}\|\le Cs_i^{-\theta}$ yields
\[
f_i(a(v_i)p)=s_i^{-1}T_{s_i}(Av_i+O(|v_i|^2))
=T_{s_i}A(v_i/s_i)+O(s_i^{1-\theta}).
\]
Taking $i \to \infty$, we have 
\begin{equation}\label{equation:limit}
\psi(u)= f_Y(T_u q) = \lim_{i \to \infty}  T_{s_i}A (v_i/s_i).
\end{equation}

Let $V_i=\operatorname{im}(T_{s_i}A)$. Since $T_{s_i}$ is invertible and $\operatorname{rank}(A)=k$, $\dim(V_i)=k$. Passing to a subsequence, assume $V_i \to V \le \mathbb{R}^n$, $\dim(V)=k$. By \eqref{equation:limit}, $ \psi(S) \le V$.

To prove the reverse inclusion, fix $u\in V$ and choose $u_i\in V_i$ converging to $u$. The equation $T_{s_i}At_i=s_iu_i$ has a unique solution $t_i\in\mathbb R^k$. The bound $\|T_{s_i}^{-1}\|\le C$ and the rank condition on $A$ give $|t_i|\le Cs_i|u_i|$. Since $s_i\to0$ and $(u_i)$ is bounded, we have $t_i\in V'$ for all sufficiently large $i$. For every $s\in(0,1]$, \eqref{equation:A} gives $F(a(st_i)p)=Ast_i+O(|t_i|^2)$. Applying $s_i^{-1}T_{s_i}$, we obtain
\[
f_i(a(st_i)p) = ss_i^{-1}T_{s_i}At_i + O(|s_i|^{1-\theta}) = su_i+O(|s_i|^{1-\theta}).  
\]
For large $i$ and any $s \in (0,1]$, $|f_i(a(st_i)p)|=|su_i+O(|s_i|^{1-\theta})| < |u|+1$.  On the other hand, by the growth estimate for $f_Y$, there exists $R>0$ such that for any $x$ with $d_Y(x,q)=R$, $|f_Y(x)|>|u|+2$. The path $s\mapsto a(st_i)p$ remains in $B_R^{r_i}(p)$ for all sufficiently large $i$: otherwise, it would meet the $R$-sphere, contradicting the preceding estimates. Hence $a(t_i)\in G^{r_i}(p,R)$. After passing to a subsequence if necessary, $a(t_i)$ converges to $T_{u'}\in H$ for some $u'\in S\le\mathbb R^l$ with $|u'|\le R$. The preceding convergence gives $f_Y(T_{u'}q)=u$, so $u=\psi(u')\in\psi(S)$. Thus $V\subseteq\psi(S)$.

Since $f_Y$ is proper, the continuous bijection
$\psi:S\to V$ has continuous inverse.
Hence $S\cong V\cong\mathbb R^k$ as topological groups.
\end{proof}

Theorem~\ref{theorem:stability of group} follows from
Corollaries~\ref{cor:iso} and~\ref{coro:gene}, together with
Lemma~\ref{lemma:S}.

\subsection{Construction of local singular fibrations}
We prove the existence of local singular fibrations in the setting of Theorem B.
\begin{theorem}\label{singular fiber}
There exists $\delta(n)>0$ with the following property. For every $\rho>0$ and $0<\delta<\delta(n)$, let $(X_i,p_i,d_i,\mathcal H^n)$ be a sequence of $\mathrm{RCD}(-(n-1),n)$ spaces satisfying the $(\rho,\delta)$ local rewinding Reifenberg condition at $p_i$. If 
\[
(X_i,p_i,d_i,\mathcal{H}^n/\mathcal{H}^n(B_1(p_i))) \overset{\mathrm{mGH}}\longrightarrow (X,p,d,\mathfrak{m}),
\]
then there exists $0<r<\rho$ such that, for all sufficiently large $i$, there is a singular fibration $\bar F_i:B_r(p_i)\to B_r(p)$ that is $\epsilon_i$-close to the original GHA on $B_r(p_i)$.  
\end{theorem}

By Theorem~\ref{thm:RCD precompact}, we can consider relative local covers,
\begin{center}
 $\begin{CD}
		(\tilde{B}(p_i,\rho/2,\rho),\tilde{p}_i,G_i) @>{\mathrm{eGH}}>> (\tilde{X},\tilde{p},G)\\
		@VV\pi V @VV\pi V\\
		B_{\rho/2}(p_i) @>{\mathrm{GH}}>> \bar{B}_{\rho/2}(p).
	\end{CD}$
\end{center}
Although $\tilde X$ need not be globally $\mathrm{RCD}$, it is locally $\mathrm{RCD}$.
We can therefore apply Theorem~\ref{theorem:stability of group} to an equivariant tangent cone at $\tilde p$.
Denote this cone by $(Y,q,H)$, and write $Y=\mathbb R^l\times C(Z)$, where $C(Z)$ does not split. $G_0$ is nilpotent by the structure of approximate groups. $H$ satisfies property Q with $\dim(\operatorname{trans}(H))=k=\dim(G/G_{\tilde{p}})$. There is an isomorphism $\iota:G_{\tilde{p}} \to H_q$. Let $T$ be the maximal torus of $H$, and set $s=\dim(T)$. Then $H_0=\mathbb{R}^k \times T$. By Theorem~\ref{theorem:stability of group}, $\dim(G_{\tilde{p}})=\dim(H_q)=s$.

The argument of Lemma~\ref{lemma:same abstract group} also gives the following result.
\begin{lemma}\label{lemma:same limit group}
For any $R>10$, there exists $R_0>0$ such that the group defined by $G^{R_0}(\tilde{p},t)$ is the same for all $t \in [1/R,6R]$.
\end{lemma}


Fix $R>10$. The model construction in the preceding section uses only the geometry and group actions on these balls and therefore applies when the equivariant distance to the fixed cone $(Y,q,H)$ is sufficiently small. Choose the comparison error $0<\epsilon<\delta/10$ accordingly, and then choose $R_0$ large enough that $(R_0\tilde X,\tilde p,G)$ is $\epsilon/2$-close to $(Y,q,H)$ on $R$-balls about the basepoints. The comparison error $\epsilon$ may depend on the fixed cone and $R$, but is independent of $i$. Set $Y_i=\tilde{B}(p_i,\rho/2,\rho)$. Then
\[
(R_0Y_i,\tilde{p}_i,G_i) \to (R_0\tilde{X},\tilde{p},G).
\]
For the remainder of the construction, we use the rescaled metrics $R_0d_i$ and $R_0d$. In particular, the displacement sets $G_i(\tilde p_i,r)$ and $G(\tilde p,r)$ are defined with respect to these metrics.

We prove Theorem~\ref{singular fiber} by contradiction. Suppose that no singular fibrations exist on the $R$-balls about the basepoints. We shall construct such fibrations after passing to a subsequence. Our strategy is to construct nilpotent models for $(R_0\tilde{X},\tilde{p},G)$ in the same way as in the last section.

By Lemma~\ref{lemma:same limit group}, the group defined by $G(\tilde{p},t)$ is the same for all $t \in [1/R,6R]$. Denote this common group by $\hat{G}_{\infty}$. Apply the gluing construction to $\hat G_\infty$ and $R_0\tilde X$ to define
\[
\hat Y_\infty=\hat G_\infty\times_{G(\tilde p,6R)}\overline B_{3R}(\tilde p).
\]
Then $(\hat{Y}_{\infty},\hat{p}_{\infty},\hat{G}_{\infty})$ is equivariantly isometric to $(R_0\tilde{X},\tilde{p},G)$ on $R$-balls about the basepoints. 

We observe that $N_{\infty}$, the identity component of $\hat{G}_{\infty}$, is a connected nilpotent Lie group of dimension $k+s$ with a maximal torus $T_{\infty}$ of dimension $s$. Let $\pi:\tilde N_\infty\to N_\infty$ be the universal covering homomorphism, and set $D_\infty:=\ker\pi$.

By the proof of Lemma~\ref{lemma:metric} and
Corollary~\ref{cor:homo along torus}, we can equip
$\tilde N_\infty$ with a left-invariant metric and construct a
$\Phi(\delta)$-$C^\infty$-close diffeomorphism
$\tilde\varphi_\infty:\tilde N_\infty\to
\tilde H_0=\mathbb R^{k+s}$ mapping $D_\infty$ onto
$D:=\ker(\pi:\tilde H_0\to H_0)$. It descends to a $\Phi(\delta)$-$C^\infty$-close diffeomorphism $\varphi_\infty:N_\infty\to H_0$ mapping $T_\infty$ to $T$. 

Define $M_{\infty}:=N_{\infty}\times_{T_{\infty}}\mathbb{R}^{n-k}$ using $\varphi_{\infty}:T_{\infty}\to T\le O(n-k)$, and equip it with the induced connection metric. Let $p_\infty'=[e_{N_\infty},0^{n-k}]$ be the base point of $M_{\infty}$. Then $M_{\infty}$ is $\Phi(\delta)$-$C^\infty$-close to $\mathbb{R}^n$. The proof of Lemma~\ref{lem:embed G_i} similarly gives an action of $\hat G_\infty$ on $M_\infty$. After a small change of metric, this action is isometric.

With $\epsilon$ sufficiently small and $\epsilon<\delta/10$,
the argument of Lemma~\ref{lemma:same abstract group} shows that,
for all sufficiently large $i$, the groupifications of
$G_i(\tilde p_i,t)$ are naturally isomorphic for all
$t\in[1/R,6R]$. Denote this common group by $\hat G_i$. Then we construct a gluing space $(\hat{Y}_i,\hat{p}_i)$ with an isometric $\hat{G}_i$-action. The triple $(\hat{Y}_i,\hat{p}_i,\hat{G}_i)$ is equivariantly isometric to $(R_0Y_i,\tilde{p}_i,G_i)$ on $R$-balls about the basepoints. $\Gamma_i \le \hat{G}_i$ is the subgroup that is close to $H_0$.

Notice that $(\hat{Y}_i,\hat{p}_i,\hat{G}_i)$ converges to $(\hat{Y}_{\infty},\hat{p}_{\infty},\hat{G}_{\infty})$ on $R$-balls. We shall construct an orbit-preserving homeomorphism between them.
 
We apply the construction in Section \ref{sec 6} to $\hat G_i \to \hat G_\infty$. Theorem~\ref{thm:Q conv} and the proof of Corollary~\ref{cor:homo along torus} give $\epsilon_i$-$C^\infty$-close diffeomorphisms 
$\varphi_{i,\infty}:N_{i} \to N_{\infty}$ so that for any $g \in T_i$ and $h \in N_i$, $\varphi_{i,\infty}(gh)=\varphi_{i,\infty}(g)\varphi_{i,\infty}(h)$. Then $M_i$ is $\epsilon_i$-$C^\infty$-close to $M_\infty$ and therefore $\Phi(\delta)$-$C^\infty$-close to $\mathbb R^n$. Lemma~\ref{lem:G_i equi} gives $\hat G_i$-equivariant homeomorphisms
\[
f_i:B_R(\hat G_i\hat p_i)\longrightarrow B_R(\hat G_i p_i')\subset M_i.
\]

As in Lemma~\ref{lem:orbit-preserving}, there are $\epsilon_i$-$C^\infty$-close diffeomorphisms $J_i:M_i\to M_\infty$ satisfying $J_i(E_i z)=\hat G_\infty J_i(z)$ for every $z\in M_i$. The uniform local Lipschitz bounds allow us, after passing to a subsequence, to take a locally uniform limit of $J_i\circ f_i$. The convergence is also strong in $H^{1,2}_{\mathrm{loc}}$ by \cite[Theorem~4.4]{AH18}. The limit map 
\[
f_\infty:B_R(\hat G_\infty\hat p_\infty)
\longrightarrow B_R(\hat G_\infty p_\infty')\subset M_\infty,
\]
is locally $(\Phi(\delta),n)$-splitting, $\Phi(\delta)$-injective and  $\hat G_\infty$-equivariant. Hence it is a homeomorphism.

The resulting orbit-preserving map $F_i:=f_\infty^{-1}\circ J_i\circ f_i$ is
\[
F_i:B_R(\hat G_i\hat p_i)\longrightarrow B_R(\hat G_\infty\hat p_\infty).
\]
It induces a local singular fibration $\bar{F}_i:B_r(p_i)\to B_r(p)$. Lemma~\ref{lem:infranil fiber} identifies the fibers of $\bar{F}_i$ as infranilmanifolds.

To complete the proof of Theorem~\ref{singular fiber}, it is enough to show that $\bar{F}_i$ is $\epsilon_i$-close to the original GHA on $B_r(p_i)$. For any sequence $x_i\in B_R(\hat p_i)$ converging to $x\in B_R(\hat p_\infty)$, we have $J_i\circ f_i(x_i)\to f_\infty(x)$. Thus $F_i(x_i)$ is $\epsilon_i$-close to $x$. Since $F_i$ is $\epsilon_i$-close to the original eGHA, the quotient $\bar{F}_i$ is also $\epsilon_i$-close to the original GHA.
   
\section{Gluing singular fibrations}
We prove Theorem~B by an inductive gluing argument.
We first realize overlapping local singular fibrations
in a common local universal cover.
We then compare their model maps and glue the fibrations
while preserving compatibility.

\subsection{Common local universal covers}
Fix a sufficiently small dimensional constant $\eta(n)>0$, and choose $\epsilon(n)$ and $\delta(n)$ so that every $\mathrm{RCD}(-(n-1)\epsilon(n),n)$ space $(X,p,d,\mathcal{H}^n)$ satisfies $\mathcal{H}^n(B_1(p))\le(1+\eta(n))\mathrm{Vol}(B_1(0^n))$, while the $(1,\delta(n))$ Reifenberg condition at $p$ implies $\mathcal{H}^n(B_1(p))\ge(1-\eta(n))\mathrm{Vol}(B_1(0^n))$. We take $\eta(n)$ smaller than the dimensional volume loss in the overlap argument below.

\begin{lemma}\label{lemma:small radius}
Assume $(X_i,p_i,d_i,\mathcal{H}^n)$ satisfies the $(\rho,\delta)$ local rewinding Reifenberg condition at all points for some $\delta < \delta(n)$. Then it also satisfies the $(r,\delta)$ local rewinding Reifenberg condition at all points for any $r<\min (\rho/3,\epsilon(n))$.
\end{lemma}
\begin{proof}
It suffices to prove the assertion at $p_i$. Let $(\tilde{B}(p_i,r,\rho),\tilde{p}_i)$ be a connected component of $\pi^{-1}(B_r(p_i))$ in $\tilde{B}(p_i,\rho)$, with $\pi(\tilde p_i)=p_i$. By assumption, \(\tilde p_i\), viewed in the original \(\rho\)-cover, satisfies the \((r,\delta)\) Reifenberg condition. Now let $(\tilde{B}(p_i,r),\tilde{p}_i')$ be the universal cover of $B_r(p_i)$. We claim that the projection
\[
\pi: (\tilde{B}(p_i,r),\tilde{p}_i') \to (\tilde{B}(p_i,r,\rho),\tilde{p}_i)
\]
is isometric on $B_{r/3}(\tilde p_i')$. Indeed, suppose there is a nontrivial deck transformation $g$ such that $gB_{2r/3}(\tilde p_i') \cap B_{2r/3}(\tilde p_i')$ contains a point $x$.  
Then $g^{-1}B_{r/10}(x) \subset B_{r}(\tilde p_i')$. By volume comparison, 
\[
\mathrm{Vol}(B_r(\tilde p_i)) \le \mathrm{Vol}(B_r(\tilde p_i')) -\frac{1}{2}\mathrm{Vol}(B_{r/10}(x)) < (1-\eta(n))\mathrm{Vol}(B_r(0^n)), 
\]
a contradiction.
\end{proof}

By Lemma~\ref{lemma:small radius}, we may henceforth assume that $\rho<\epsilon(n)$.
An argument similar to the proof of Lemma~\ref{lemma:small radius} shows that a small change in the center of the ball used to take the universal cover preserves the relevant local information.
\begin{lemma}\label{lemma:same local uni}
Assume $(X_i,p_i,d_i,\mathcal{H}^n)$ satisfies the $(\rho,\delta)$ local rewinding Reifenberg condition at all points, with $\delta<\delta(n)$. For $x,y\in X_i$ with $d_i(x,y)<\rho/10$, let $B_1=\tilde{B}(x,\rho/3)$ and $B_2=\tilde{B}(y,\rho/3)$, with deck groups $G$ and $H$, respectively. Choose lifts $\tilde x\in B_1$ and $\tilde x'\in B_2$ of $x$. Then $(B_{\rho/200}(\tilde x),G(\tilde x,\rho/200))$ and $(B_{\rho/200}(\tilde x'),H(\tilde x',\rho/200))$ are equivariantly isometric.
\end{lemma}

Consider the convergence in Theorem B. For $x,y\in X$ with $d(x,y)<\rho/10$, choose $x_i,y_i\in X_i$ converging to $x,y$, respectively. Construct the local singular fibrations near $x$ and $y$ using the universal covers $\tilde{B}(x_i,\rho/3)$ and $\tilde{B}(y_i,\rho/3)$. These constructions use only local information in the covers; see Lemmas \ref{lemma:same abstract group} and \ref{lemma:same limit group}. By Lemma~\ref{lemma:same local uni}, the same fibration near $x$ can be constructed in $\tilde{B}(y_i,\rho/3)$. Choosing the local construction scales sufficiently small, we may therefore work in a common local universal cover whenever the fibration domains overlap.

\subsection{Comparison of overlapping local fibrations}\label{subsection:compare}
Fix $x,y\in X$ with $d(x,y)<\rho/10$, and suppose that local singular fibrations have been constructed near $x$ and $y$. Choose $x_i,y_i\in X_i$ converging to $x,y$, respectively, and let $Y_i=\tilde{B}(x_i,\rho/10,\rho/3)$ be the local relative cover. Then
\[
(Y_i,\tilde{x}_i,\tilde{y}_i,G_i) \longrightarrow (Y,\tilde{x},\tilde{y},G),
\]
where $\pi(\tilde{x}_i)=x_i$ and $\pi(\tilde{y}_i)=y_i$. Fix $R>10$. Choose $r_x<\min\{1,\rho\}/(10^4R)$ small enough that we can construct singular fibrations on $B_{Rr_x}(x)$. Choose $r_y$ similarly at $y$. 

We use $B^{r_x^{-1}}_{3R}(\tilde x_i)$ and
$G_i^{r_x^{-1}}(\tilde x_i,6R)$ to construct
$(\hat Y_{i,x},\hat G_{i,x})$;
see Subsection~\ref{subsection:same orbits}.
We similarly construct
$(\hat Y_{\infty,x},\hat G_{\infty,x})$.
The local singular fibration is obtained from the
following composition:
\[
(\hat{Y}_{i,x},\hat{x}_{i,x},\hat{G}_{i,x}) \overset{f_{i,x}}\to (M_{i,x},\hat{x}_{i,x}',\hat{G}_{i,x},E_{i,x}=\hat{G}_{i,x}N_{i,x}),
\]
\[
(M_{i,x},\hat{x}_{i,x}',\hat{G}_{i,x},E_{i,x}=\hat{G}_{i,x}N_{i,x}) \overset{J_{i,x}}\to (M_{\infty,x},\hat{x}_{\infty,x}',\hat{G}_{\infty,x}),
\]
\[ (M_{\infty,x},\hat{x}_{\infty,x}',\hat{G}_{\infty,x}) \overset{(f_{\infty,x})^{-1}}\to (\hat{Y}_{\infty,x},\hat{x}_{\infty,x},\hat{G}_{\infty,x}).
\]

The construction at $y$ is the same, with $x$ replaced by $y$. We compare the constructions when $B_{r_x}(x)\cap B_{r_y}(y)\ne\emptyset$, and may assume $r_x\ge r_y$. Since the choice of $f_{i,x}$ is not unique, the singular fibrations $\bar{F}_{i,x}$ and $\bar{F}_{i,y}$ constructed from the models at $x_i$ and $y_i$ need not coincide, even though both are close to the original GHA from $X_i$ to $X$.

We call $f_{i,x}$ the singular fiber bundle map to the nilpotent model at $\hat{x}_{i,x}$, since it determines the fiber partition of the quotient singular fibration. The limit model and $f_{\infty,x}$ determine how these fibers are assigned to points of $X$. In this subsection, we compare only $f_{i,x}$ and $f_{i,y}$. 

The difficulty in comparing $f_{i,x}$ and $f_{i,y}$ is that their nilpotent models may be very different, especially when $r_x$ is much larger than $r_y$. The key observation is that $B_{(R-2)r_y}(y_i)\subset B_{Rr_x}(x_i)$, and all the local group information used to construct the model at $y_i$ is visible in this smaller ball. We can therefore reconstruct the model at $y_i$ using the model at $x_i$.

Let $\hat{y}_{i,x} \in \hat{Y}_{i,x}$ be a lift of $y_i$ and let $f_{i,x}(\hat{y}_{i,x})=\hat{y}_{i,x}' \in M_{i,x}$. Let $s=\frac{r_y}{r_x} \le 1$. We first study $f_{i,x}$ near $\hat{y}_{i,x}$ after blowing up the metric by $s^{-1}$.

\begin{lemma}\label{lemma:group id}
We can identify $(B^{s^{-1}}_R(\hat{y}_{i,x}),\hat{y}_{i,x},\hat{G}_{i,x})$ with $(\hat{Y}_{i,y},\hat{y}_{i,y},\hat{G}_{i,y})$ on the $(R-2)$-balls about their basepoints. In particular, this induces a natural homomorphism $\theta_{i,yx}:E_{i,y} \to E_{i,x}$ which is injective on $E_{i,y}(\hat{y}_{i,y}',R-2)$. 
\end{lemma}
\begin{proof}
Recall that we can identify $(\hat{Y}_{i,x},\hat{G}_{i,x})$ with $(r_x^{-1}Y_i,\tilde{x}_i,G_i)$ on an $R$-ball and identify $(\hat{Y}_{i,y},\hat{G}_{i,y})$ with $(r_y^{-1}Y_i,\tilde{y}_i,G_i)$ on an $R$-ball. By our choice of $r_y$ and $s=r_y/r_x \le 1$, both $(B^{s^{-1}}_R(\hat{y}_{i,x}),\hat{G}_{i,x})$ and $(\hat{Y}_{i,y},\hat{G}_{i,y})$ agree with $(r_y^{-1}Y_i,\tilde{y}_i,G_i)$ on the $(R-2)$-balls about their basepoints, so we can identify them. 

We first compare the defining local groups. Let $S_{i,y}$ be the chosen defining subset of $\hat{G}_{i,y}$ for its action at $y_i$. The common covering chart identifies $S_{i,y}$ with a subset of $\hat{G}_{i,x}$ carrying the same local multiplication, and every defining product is computed in the same local covering chart. The defining property therefore gives a unique homomorphism
\[
q_{i,yx}:\hat G_{i,y}\longrightarrow\hat G_{i,x}
\]
extending this identification. 

The subgroups $\Gamma_{i,x}$ and $\Gamma_{i,y}$ converge
to the identity components of $\hat G_{\infty,x}$ and
$\hat G_{\infty,y}$, respectively.
Under the common covering identification, we claim that
\[
\Gamma_{i,y}(\hat y_{i,y},R-2)
\subset
\Gamma_{i,x}^{s^{-1}}(\hat y_{i,x},R).
\]
It suffices to check this on the short generators used
in the defining-set construction.
Let $g_i$ be such a generator, chosen to be
$\epsilon_i$-close to the identity.
We claim that it remains close to the identity in
$\hat G_{i,x}$.
Otherwise, after passing to a subsequence, it would
converge to a nontrivial action on
$B_1(\hat x_{\infty,x})\subset\hat Y_{\infty,x}$
that fixes $B_s(\hat y_{\infty,x})$ pointwise,
where $\hat y_{\infty,x}$ is the limit of $\hat y_{i,x}$.
This contradicts \cite[Lemma~3.2]{PanWang2021}.
The claimed inclusion follows.

Since $\Gamma_{i,y}(\hat{y}_{i,y}, R-2) \subset \Gamma_{i,x}^{s^{-1}}(\hat{y}_{i,x}, R)$, $q_{i,yx}$ induces an injective homomorphism, 
\[
\widetilde\beta_{i,yx}:\widetilde N_{i,y}\longrightarrow\widetilde N_{i,x}.
\]
By the proof of Lemma~\ref{lemma:same kernel}, $\widetilde\beta_{i,yx}(D_{i,y})\subset D_{i,x}$. Thus $\widetilde\beta_{i,yx}$ descends to $\beta_{i,yx}:N_{i,y}\to N_{i,x}$, agreeing with $q_{i,yx}$ on $\Gamma_{i,y}$.

For every $g\in\hat G_{i,y}(\hat y_{i,y},R-2)$,
$\beta_{i,yx}$ intertwines conjugation by $g$ with
conjugation by $q_{i,yx}(g)$.
Thus, combining $\beta_{i,yx}$ and $q_{i,yx}$
gives the required homomorphism
\[
\theta_{i,yx}:E_{i,y}\to E_{i,x}.
\]

Since all relations of $\Gamma_{i,y}(\hat{y}_{i,y}, R-2)$ can be observed in $\Gamma_{i,x}^{s^{-1}}(\hat{y}_{i,x}, R)$, $\theta_{i,yx}$ is injective on $E_{i,y}(\hat{y}_{i,y}',R-2)$. 
\end{proof}

In Theorem B, $E_{i,y}$ and $E_{i,x}$ have the same dimension, so $\theta_{i,yx}(E_{i,y})$ contains $N_{i,x}$. At scale $r_y$, both models are close to the same equivariant tangent model at $\tilde y$. After discarding the generators and period relations outside this local group, Lemma~\ref{lemma:same abstract group} identifies their local actions through $\theta_{i,yx}$.

\begin{remark}
We briefly explain which generators and relations are discarded. For simplicity, assume $G_{\tilde{y}}\le G_{\tilde{x}}$. A circle subgroup of $G_{\tilde{x}}$ that moves $\tilde{y}$ may appear as a line at scale $r_y$; its period relation is then discarded. Components of $G_{\tilde{x}}$ that move $\tilde{y}$ by a definite distance correspond to the generators that are discarded.
\end{remark}

Suppose $z_i\to z$ and a third model is constructed on $B_{Rr_z}(z_i)$, with $r_z\le r_y\le r_x$ and $B_{r_x}(x_i)\cap B_{r_y}(y_i)\cap B_{r_z}(z_i)\ne\emptyset$. The proof of Lemma \ref{lemma:group id} implies
\[
\theta_{i,zx}=\theta_{i,yx} \circ \theta_{i,zy}.
\]

\begin{defn}\label{defn:compat}
Let $S\subset B_{r_x}(x_i)\cap B_{r_y}(y_i)$.
We say that $f_{i,x}$ is compatible with $f_{i,y}$ on $S$
if, after identifying the covering neighborhoods as in
Lemma~\ref{lemma:group id}, the following implication
holds for all lifts $z_1,z_2\in\hat Y_{i,y}$ lying over
$S$ and all $g\in E_{i,y}$:
\[
f_{i,y}(z_1)=g f_{i,y}(z_2)
\quad\Longrightarrow\quad
f_{i,x}(z_1)=\theta_{i,yx}(g)f_{i,x}(z_2).
\]
\end{defn}
Because $E_{i,y}$ and $E_{i,x}$ have the same dimension, the identification of local actions above shows that compatible singular fiber bundle maps induce the same fiber partition on the overlap.


The choices of almost splitting maps in the constructions of $f_{i,x}$ and $f_{i,y}$ are not unique, so these maps need not be compatible. The next lemma constructs a map $f_{i,x,y}$ in the model at $y_i$ that is close to $f_{i,y}$ and with which $f_{i,x}$ is compatible.

\begin{lemma}\label{lemma:change scale}
There exists a $\hat{G}_{i,y}$-equivariant almost splitting map $f_{i,x,y}$, $\Phi(\delta)$-close to $f_{i,y}$,
\[
(\hat{Y}_{i,y},\hat{y}_{i,y},\hat{G}_{i,y}) \overset{f_{i,x,y}}\to (M_{i,y},\hat{y}_{i,y}',\hat{G}_{i,y},E_{i,y}),
\]
such that $f_{i,x}$ is compatible with $f_{i,x,y}$ on $B_{r_y}(y_i)$.
\end{lemma}
\begin{proof}
Apply the transformation theorem to $f_{i,x}$ in normal coordinates at $\hat y_{i,x}'$ with $s=r_y/r_x$. Write $A_{i,s}$ for its linear transformation matrix and $T_{i,s}$ for the corresponding coordinate map on $M_{i,x}$. Then $s^{-1}T_{i,s}\circ f_{i,x}$ remains $(\Phi(\delta)+\epsilon_i)$-almost splitting on the $R$-ball. Choose these matrices consistently with the $H^{1,2}$ convergence to the limit map, and let $T_{\infty,s}$ be the resulting transformation on $M_{\infty,x}$. Transfer this coordinate change back to \(M_{i,x}\) by setting
\[
T_{i,s}'=J_{i,x}^{-1}\circ T_{\infty,s}\circ J_{i,x}.
\]
For the fixed ratio $s=r_y/r_x$ and each derivative order $m$, $T_{i,s}'$ differs from $T_{i,s}$ by $\epsilon_i$ in the coordinate $C^m$ norm. Hence $f_{i,x,y}=s^{-1}T_{i,s}'\circ f_{i,x}$ is still $(\Phi(\delta)+\epsilon_i)$-almost splitting. Here $s^{-1}$ denotes rescaling of the target metric.

Conjugate the original $\hat G_{i,x}$- and $E_{i,x}$-actions by $T_{i,s}'$. With these conjugated actions, $f_{i,x,y}$ is $\hat G_{i,x}$-equivariant. Moreover,
\[
f_{i,x}(z_1)=g f_{i,x}(z_2)
\quad\Longleftrightarrow\quad
f_{i,x,y}(z_1)=
\bigl(T_{i,s}'g(T_{i,s}')^{-1}\bigr)f_{i,x,y}(z_2),
\]
so the change preserves compatibility with $f_{i,x}$.

The map $T_{i,s}'$ need not be an isometry, so the conjugated $E_{i,x}$-action need not preserve the original metric on $s^{-1}M_{i,x}$. We therefore average the metric while changing it only slightly.

We observe that $f_{i,x,y}:B^{s^{-1}}_R(\hat{y}_{i,x}) \to (s^{-1}M_{i,x},\hat{y}_{i,x}')$ is a $\hat{G}_{i,x}$-equivariant almost splitting map and $\hat{G}_{i,x}$ acts isometrically on $B^{s^{-1}}_R(\hat{y}_{i,x})$. Thus the conjugated $\hat{G}_{i,x}$-action is almost isometric on $(s^{-1}M_{i,x},\hat{y}_{i,x}')$. Since $\hat{G}_{i,x}$ is $\epsilon_i$-dense in $E_{i,x}$, the conjugated $E_{i,x}$-action is almost isometric on $(s^{-1}M_{i,x},\hat{y}_{i,x}')$. 

We now obtain smooth control of these conjugated actions. Use the linear transformation matrix $A_{i,s}$ introduced above. The transformation estimates give
\[
\|A_{i,s}\|\le Cs^{-\alpha},\qquad
\|A_{i,s}^{-1}\|\le C,\qquad 0<\alpha<1.
\]
For any isometry $g_i \in E_{i,x}^{s^{-1}}(\hat{y}_{i,x}',R)$ on $s^{-1}M_{i,x}$,  we have $\|D^jg_i\|\le C_js^{j-1}$. Since $A_{i,s}$ is linear,
\[
D^j(A_{i,s}g_iA_{i,s}^{-1})(z)
=A_{i,s}(D^jg_i)(A_{i,s}^{-1}z)
 [A_{i,s}^{-1},\ldots,A_{i,s}^{-1}],
\]
and therefore
\[
\|D^j(A_{i,s}g_iA_{i,s}^{-1})\|
\le C_j's^{j-1-\alpha}\qquad(j\ge2).
\]
These bounds are uniform. The $C^\infty$ closeness of the metric to the flat metric controls the coordinate and metric corrections on these comparison balls. Replacing the coordinate map by $T_{i,s}'$ introduces an error
of size $\epsilon_i$ for each fixed $s$ and derivative order. Together with the preceding $C^0$ closeness, interpolation
on smaller balls gives $(\Phi_m(\delta)+\epsilon_i)$-$C^m$
closeness to Euclidean isometries in these coordinates,
for every fixed $m$, where $\Phi_m(\delta)\to0$
as $\delta\to0$.

Apply the averaging construction of Theorem~\ref{thm:inv} to the noncompact virtually nilpotent group $E_{i,x}$. This gives a new metric on $s^{-1}M_{i,x}$ that is $(\Phi(\delta)+\epsilon_i)$-$C^\infty$-close to the original metric. The conjugated $E_{i,x}$-action is isometric for the new metric. This change of metric does not affect the compatibility of $f_{i,x}$ and $f_{i,x,y}$.

It remains to compare this rescaled model with the model at $y_i$. Through $\theta_{i,yx}$, both carry the same local $E_{i,y}$-action. We seek an $E_{i,y}$-equivariant diffeomorphism, $(\Phi(\delta)+\epsilon_i)$-$C^\infty$-close to the coordinate comparison,
\[
g_i:(s^{-1}M_{i,x},\hat y_{i,x}',E_{i,y})\longrightarrow(M_{i,y},\hat y_{i,y}',E_{i,y}),
\] 
such that $g_i\circ f_{i,x,y}$ is $(\Phi(\delta)+\epsilon_i)$-close to $f_{i,y}$.

Since both $s^{-1}M_{i,x}$ and $M_{i,y}$ are $C^{\infty}$-close to $\mathbb{R}^n$, there is a diffeomorphism (on an $R$-ball) $g_i':(s^{-1}M_{i,x},\hat{y}_{i,x}') \to (M_{i,y},\hat{y}_{i,y}')$ so that $g_i' \circ f_{i,x,y}$ is close to $f_{i,y}$ and $g_i'$ is $C^{\infty}$-close to an orthogonal linear map on $\mathbb{R}^n$.

The coordinate comparison $g_i'$ may not be equivariant. Since $f_{i,y}$ and $f_{i,x,y}$ are $\hat G_{i,y}$-equivariant and $g_i'\circ f_{i,x,y}$ is close to $f_{i,y}$, $g_i'$ is almost $\hat G_{i,y}$-equivariant. Density in $E_{i,y}$ gives almost $E_{i,y}$-equivariance on local charts of the action. We extend this comparison through the local action charts
and apply the same averaging construction to obtain an
$E_{i,y}$-equivariant comparison $g_i$.
The smooth estimates above keep $g_i$ $C^{\infty}$-close to
$g_i'$, so it remains a diffeomorphism on smaller balls.
Replacing $f_{i,x,y}$ by $g_i\circ f_{i,x,y}$ completes
the proof.
\end{proof}

\begin{remark}
The construction above transfers the coordinate change $T_{\infty,s}$ to $M_{i,x}$ through $J_{i,x}$. In Lemma~\ref{lemma:glue fibers}, we apply the same coordinate change to $f_{\infty,x}$, so $J_{i,x}$ remains an orbit-preserving diffeomorphism realizing $C^{\infty}$ convergence. In the proof of Theorem C, where no limit space is fixed, we instead change coordinates using the transformation map of $f_{i,x}$ itself.
\end{remark}

The preceding lemmas compare fibers on $B_{r_x}(x)$ and $B_{r_y}(y)$. To glue them in Theorem B, we work on a ball $B_h(z)\subset B_{r_x}(x)\cap B_{r_y}(y)$ and represent $f_{i,x}$ and $f_{i,y}$ in a common model.

Assume $z_i \in X_i$ converges to $z \in X$ and choose a lift $\hat{z}_{i,y}\in\hat Y_{i,y}$ of $z_i$. We describe the map $f_{i,y}$ on the covering neighborhood of $B_h(z_i)$. Let $s_y= \frac{h}{r_y} \le 1$ and $\hat{z}_{i,y}'=f_{i,y}(\hat{z}_{i,y})$.

Applying the transformation theorem and averaging argument as in Lemma~\ref{lemma:change scale} gives the following result.
\begin{lemma}\label{lemma:fiber at small scale}
After a small change of the metric on $(s^{-1}_yM_{i,y},\hat{z}_{i,y}')$ if necessary, there exists a $\hat{G}_{i,y}$-equivariant almost splitting map 
\[
f_{i,y,z}: (s^{-1}_y\hat{Y}_{i,y},\hat{z}_{i,y},\hat{G}_{i,y}) \to (s^{-1}_yM_{i,y},\hat{z}_{i,y}',\hat{G}_{i,y},E_{i,y}),
\]
such that $f_{i,y}$ is compatible with $f_{i,y,z}$.
\end{lemma}

Combining Lemmas \ref{lemma:change scale} and \ref{lemma:fiber at small scale}, we can represent the singular fiber bundle map $f_{i,x}$ on $B_h(z_i)$ in the nilpotent model at $y_i$.
\begin{lemma}\label{lemma:x-fiber at small scale}
With the metric supplied by Lemma~\ref{lemma:fiber at small scale}, there is a $\hat{G}_{i,y}$-equivariant almost splitting map 
\[
f_{i,x,y,z}: (s^{-1}_y\hat{Y}_{i,y},\hat{z}_{i,y},\hat{G}_{i,y}) \to (s^{-1}_yM_{i,y},\hat{z}_{i,y}',\hat{G}_{i,y},E_{i,y}),
\]
such that $f_{i,x}$ is compatible with $f_{i,x,y,z}$. Moreover, $f_{i,x,y,z}$ is $(\Phi(\delta)+\epsilon_i)$-close to $f_{i,y,z}$ in Lemma~\ref{lemma:fiber at small scale}.
\end{lemma}

\subsection{Gluing two local singular fibrations}
Retain the notation of the preceding subsection. Suppose that, for some $h>0$, $B_{2h}(z_1),B_{2h}(z_2)\subset B_{r_x}(x)\cap B_{r_y}(y)$, and choose $z_{i,1},z_{i,2}\in X_i$ converging to $z_1,z_2\in X$, respectively.
Let $\bar{F}_{i,x}:B_h(z_{i,1}) \to B_{h}(z_1)$ be the singular fibration induced by $f_{i,x}$ through the model construction at $x_i\to x$. Similarly, let $\bar{F}_{i,y}:B_h(z_{i,2}) \to B_{h}(z_2)$ be the singular fibration induced by $f_{i,y}$ through the model construction at $y_i\to y$. Recall that both $\bar{F}_{i,x}$ and $\bar{F}_{i,y}$ are $\epsilon_i$-close to the original GHA from $X_i$ to $X$. 

Our goal in this subsection is to glue $\bar{F}_{i,x}$ and $\bar{F}_{i,y}$.
\begin{lemma}\label{lemma:glue fibers}
There exists a singular fibration
\[
\bar{F}_{i,xy}:B_{0.9h}(z_{i,1})\cup B_{0.9h}(z_{i,2})
\longrightarrow B_{0.9h}(z_1)\cup B_{0.9h}(z_2)
\]
induced by a singular fiber bundle map $f_{i,xy}$, with the following properties:
\begin{itemize}
\item  On $B_{0.9h}(z_{i,1})$, we have
$\bar F_{i,xy}=\bar F_{i,x}$, and $f_{i,x}$ is
compatible with $f_{i,xy}$.
\item On $B_{0.9h}(z_{i,2})\setminus B_{1.1h}(z_{i,1})$, $f_{i,y}$ is compatible with $f_{i,xy}$.
\item $\bar{F}_{i,xy}$ is $\epsilon_i$-close to the original GHA on the same domain.
\end{itemize}
\end{lemma}

We may assume $d(z_1,z_2)<2h$, since otherwise there is no need to glue them. We retain the notation of the preceding subsection and set $s_y=h/r_y$. Let $\hat{z}_{i,y,1},\hat{z}_{i,y,2} \in s^{-1}_y\hat{Y}_{i,y}$ be lifts of $z_{i,1}$ and $z_{i,2}$. Let $\hat{z}_{i,y,1}'=f_{i,y}(\hat{z}_{i,y,1})$ and $\hat{z}_{i,y,2}'=f_{i,y}(\hat{z}_{i,y,2})$.

We apply Lemmas \ref{lemma:x-fiber at small scale} and \ref{lemma:fiber at small scale} for $x,y$ and $z=z_1$. On
\[
(s^{-1}_y\hat{Y}_{i,y},\hat{z}_{i,y,1},\hat{z}_{i,y,2},\hat{G}_{i,y}),
\]
we have two $\hat{G}_{i,y}$-equivariant almost splitting maps $f_{i,x,y,z}$ and $f_{i,y,z}$ such that $f_{i,x}$ is compatible with $f_{i,x,y,z}$ and $f_{i,y}$ is compatible with $f_{i,y,z}$. They are $(\Phi(\delta)+\epsilon_i)$-close to each other. We first glue these two maps.

\begin{lemma}\label{lemma:glue almost splitting}
There exists a $\hat{G}_{i,y}$-equivariant almost splitting map 
\[
\begin{aligned}
f_{i,xy}:{}&\bigl(B_{0.95}(\hat{G}_{i,y}\hat{z}_{i,y,1})
 \cup B_{0.95}(\hat{G}_{i,y}\hat{z}_{i,y,2}),
 \hat{z}_{i,y,1},\hat{z}_{i,y,2},\hat{G}_{i,y}\bigr)\\
&\longrightarrow(s_y^{-1}M_{i,y},\hat{z}_{i,y,1}',
 \hat{z}_{i,y,2}',\hat{G}_{i,y},E_{i,y}),
\end{aligned}
\]
which is equal to $f_{i,x,y,z}$ on $B_{0.95}(\hat{G}_{i,y}\hat{z}_{i,y,1})$ and equal to $f_{i,y,z}$ on $B_{0.95}(\hat{G}_{i,y}\hat{z}_{i,y,2})\setminus B_{1.05}(\hat{G}_{i,y}\hat{z}_{i,y,1})$. Moreover, $f_{i,xy}$ is $(\Phi(\delta)+\epsilon_i)$-close to $f_{i,y,z}$ and $f_{i,x,y,z}$.
\end{lemma}

We first construct invariant cutoff functions on $s^{-1}_yM_{i,y}$ with uniform derivative bounds. 
\begin{lemma}\label{lemma:cut-off}
For $\delta$ small enough, there is an $E_{i,y}$-invariant smooth non-negative function $0 \le \chi_1 \le 1$ on $(s^{-1}_yM_{i,y},\hat{z}_{i,y,1}')$ such that $\operatorname{supp}(\chi_1) \subset B_{1.03}(E_{i,y}\hat{z}_{i,y,1}')$ and $\chi_1=1$ on $B_{0.97}(E_{i,y}\hat{z}_{i,y,1}')$. Moreover, for any $j$, $|D^j \chi_1| \le C_j(n,\delta)$ where $C_j$ is independent of $i$.
\end{lemma}

\begin{proof}
Choose smooth functions $0 \le \psi_1 \le \psi_2 \le 1$ on $s^{-1}_yM_{i,y}$ with $C^j$ norms uniformly bounded by $C_j' >0$, so that $\psi_1=1$ on $B_{1.02}(\hat{z}_{i,y,1}')$, $\operatorname{supp}(\psi_1) \subset B_{1.03}(\hat{z}_{i,y,1}')$, $\psi_2=1$ on $B_{1.1}(\hat{z}_{i,y,1}')$, and $\operatorname{supp}(\psi_2) \subset B_{1.15}(\hat{z}_{i,y,1}')$. This can be done since $s^{-1}_yM_{i,y}$ is $(\Phi(\delta)+\epsilon_i)$-$C^{\infty}$-close to $\mathbb{R}^n$.

The group $E_{i,y}$ is a finite extension of the nilpotent group $N_{i,y}$ and is therefore unimodular. We use its Haar measure $dg$, which is both left- and right-invariant. 
For $z\in s_y^{-1}M_{i,y}$, define
\[
A(z)=\int_{E_{i,y}}\psi_1(gz)\,dg,
\qquad
B(z)=\int_{E_{i,y}}\psi_2(gz)\,dg.
\]
These integrals are finite because the action is proper
and the functions have compact support.
Set
\[
\psi_3(z)=
\begin{cases}
A(z)/B(z),&B(z)>0,\\
0,&B(z)=0.
\end{cases}
\]
Since $\psi_2=1$ on a neighborhood of
$\operatorname{supp}(\psi_1)$, this zero extension is smooth.
Haar invariance makes $\psi_3$ $E_{i,y}$-invariant.
Moreover, $0\le\psi_3\le1$ and
\[
\operatorname{supp}(\psi_3)
\subset B_{1.03}(E_{i,y}\hat z_{i,y,1}').
\]

Let $\mu_i$ denote the Haar measure $dg$ fixed above. Let $D_1(z)=\mu_i(E_{i,y}(z,0.05))$ and $D_2(z)=\mu_i(E_{i,y}(z,2.3))$. We claim that $D_2(z)\le ND_1(z)$ for a uniform constant $N$. By volume comparison, cover $\bar B_{2.3}(z)$ by at most $N$ balls $B_j$ of diameter $0.05$. Define \[ S_j:=\{g\in E_{i,y}(z,2.3):gz\in B_j\}. \] For each nonempty $S_j$, choose $h_j\in S_j$. If $g\in S_j$, then $h_j^{-1}g\in E_{i,y}(z,0.05)$. Thus $S_j\subset h_jE_{i,y}(z,0.05)$, and left invariance gives $\mu_i(S_j)\le D_1(z)$. Since the sets $S_j$ cover $E_{i,y}(z,2.3)$, summing these inequalities proves the claim.

Notice that for any $z \in B_{0.97}(\hat{z}_{i,y,1}')$, $A(z) \ge D_1(z)$ and $B(z) \le D_2(z)$. Thus $\psi_3(z) \ge 1/N$.

We next establish uniform $C^j$ bounds for $\psi_3$.
By $E_{i,y}$-invariance, it suffices to consider
$z\in B_{1.03}(\hat z_{i,y,1}')$.
\[
D \psi_3 = \frac{B\,DA-A\,DB}{B^2}.
\]
At $z$, we have $|DA|,|DB|\le C_1'D_2$, $|A|,|B|\le D_2$, and $B\ge D_1$. Hence $|D\psi_3(z)|\le 2C_1'N^2$. The higher derivative bounds follow in the same way.

Choose a smooth nondecreasing function
$\psi_4:[0,1]\to[0,1]$ satisfying $\psi_4(0)=0$ and
$\psi_4(t)=1$ for $t\ge1/N$, and set
$\chi_1=\psi_4\circ\psi_3$.
The invariance, support, and derivative bounds follow
from those of $\psi_3$. The lower bound
$\psi_3\ge1/N$ gives $\chi_1=1$ on the required
inner orbit neighborhood.
\end{proof}

\begin{proof}[Proof of Lemma~\ref{lemma:glue almost splitting}]
Choose the invariant cutoff $\chi_1$ from Lemma~\ref{lemma:cut-off} and define the weights
\[
\lambda_1:=\chi_1\circ f_{i,x,y,z},
\qquad
\lambda_2:=1-\lambda_1.
\]
These weights are $\hat G_{i,y}$-invariant. For large $i$, we have $\lambda_1=1$ on $B_{0.95}(\hat G_{i,y}\hat z_{i,y,1})$ and $\lambda_1=0$ outside $B_{1.05}(\hat G_{i,y}\hat z_{i,y,1})$.

Take the weighted center of mass of $f_{i,x,y,z}$ and $f_{i,y,z}$ with weights $\lambda_1$ and $\lambda_2$, respectively, and denote the resulting map by $f_{i,xy}$. The center-of-mass estimates preserve the almost splitting property and the closeness to both maps. The invariant weights and equivariance of the two maps make $f_{i,xy}$ $\hat G_{i,y}$-equivariant. Where $\lambda_1=1$, it equals $f_{i,x,y,z}$; where $\lambda_1=0$, it equals $f_{i,y,z}$. This gives the claimed agreement in the lemma.
\end{proof}

To obtain a quotient fibration that is $\epsilon_i$-close
to the original GHA, we use the limit $f_{\infty,xy}$
of the glued maps to label their fibers.
Using $f_{\infty,y}$ would leave an error controlled
only by $\Phi(\delta)+\epsilon_i$.
We then adjust the base map to recover agreement with
$\bar F_{i,x}$ on $B_{0.9h}(z_{i,1})$.

\begin{proof}[Proof of Lemma~\ref{lemma:glue fibers}]
Write $\hat z_{i,y,1}'\to
\hat z_{\infty,y,1}'$.
Apply the coordinate changes used in constructing
$f_{i,y,z}$ to both the approximating and limit models,
and conjugate their group actions accordingly.
Average the metrics over the conjugated $E_{i,y}$-
and $\hat G_{\infty,y}$-actions, respectively.
These constructions are compatible with the convergence.
Thus
\[
J_{i,y}:
(s_y^{-1}M_{i,y},\hat z_{i,y,1}')
\longrightarrow
(s_y^{-1}M_{\infty,y},\hat z_{\infty,y,1}')
\]
remains an $\epsilon_i$-$C^\infty$-close diffeomorphism
for the new metrics.
It also remains $\epsilon_i$-almost equivariant with
respect to the group comparison and preserves the
corresponding orbits.

After passing to a subsequence, $J_{i,y}\circ f_{i,xy}$ converges in $H^{1,2}$ to a $\hat{G}_{\infty,y}$-equivariant almost splitting map $f_{\infty,xy}$. We use $f_{\infty,xy}$, instead of $f_{\infty,y}$, to assign the fibers determined by $f_{i,xy}$ to points of the base. This gives an orbit-preserving homeomorphism $F_{i,xy}'$ and its quotient singular fibration
\[
\bar{F}_{i,xy}':B_{0.94h}(z_{i,1})\cup B_{0.94h}(z_{i,2})
\longrightarrow B_{0.94h}(z_1)\cup B_{0.94h}(z_2)
\]
that is $\epsilon_i$-close to the original GHA.

We next deform $\bar{F}_{i,xy}'$ so that it agrees with $\bar{F}_{i,x}$ on $B_{0.9h}(z_{i,1})$. Set $S=B_{0.93h}(z_1)$.
For sufficiently large $i$, every fiber of
$\bar F_{i,xy}'$ over $S$ lies in
$B_{0.94h}(z_{i,1})$.
On this region, $\bar F_{i,xy}'$ and $\bar F_{i,x}$
have the same fiber partition. Hence they determine
an open embedding $g_i:S\to B_{0.95h}(z_1)$ satisfying
\[
(\bar F_{i,xy}')^{-1}(z)
=\bar F_{i,x}^{-1}(g_i(z)).
\]
The $\epsilon_i$-closeness of both maps to the original
GHA implies that $g_i$ converges to the identity.

The equivariant homeomorphism $f_{\infty,y}$ identifies
the quotient locally with the quotient of the smooth model.
Choose an open neighborhood $U$ with
$B_{3h}(z_2)\subset U\subset B_{10h}(z_2)$ and equip it
with the induced orbit-type CS stratification.
For each relevant closed isotropy subgroup $H$, let
$C_H\subset U$ consist of points having a lift whose
isotropy contains a conjugate of $H$.
The slice structure shows that $C_H$ is a relatively
closed union of orbit-type stratum components, hence
a sub-CS set. Compatibility of the model maps gives
\[
g_i(S\cap C_H)=g_i(S)\cap C_H.
\]
Thus $g_i$ thoroughly respects these subsets.

By Siebenmann's deformation theorem
\cite[Theorem~2.3 and Complement~4.3]{Siebenmann1972},
for all sufficiently large $i$ there is an ambient isotopy
$g_{i,t}:U\to U$, $t\in[0,1]$, preserving the specified
sub-CS sets, with $g_{i,0}=\mathrm{Id}$, $g_{i,1}=g_i$ on
$K=\overline B_{0.92h}(z_1)\subset S$, and support in
$B_{0.94h}(z_1)$. After choosing another $\epsilon_i \to 0$ if necessary, these homeomorphisms are $\epsilon_i$-close to the identity. Since $\bar{F}_{i,xy}'(B_{0.9h}(z_{i,1}))\subset K$ for large $i$, the map $\bar{F}_{i,xy}=g_{i,1}\circ \bar{F}_{i,xy}'$ agrees with $\bar{F}_{i,x}$ on $B_{0.9h}(z_{i,1})$. 

This deformation changes only the assignment to points of the base, leaving the maps to the nilpotent models unchanged. Lemma~\ref{lemma:glue almost splitting} gives the required compatibility.
\end{proof}
\begin{remark}
The deformed map is still the quotient of a local
orbit-preserving homeomorphism.
Let $\tilde U$ denote the preimage of $U$ in the
limiting local covering space.
Since $g_{i,t}$ preserves orbit types, the results of
\cite{Palais1960,Palais1961} give an equivariant lift
$\tilde g_{i,t}$ on $\tilde U$, with
$\tilde g_{i,0}=\mathrm{Id}$.
The lifts may be chosen to converge to the identity
as $i\to\infty$.
Consequently, $\bar F_{i,xy}$ is the quotient of
$\tilde g_{i,1}\circ F_{i,xy}'$.
\end{remark}

\subsection{Global gluing}
In this subsection, we prove Theorem B by gluing the local singular fibrations. Fix $R>10$ and choose finitely many points $T=\{y_j\mid 1\le j\le J\}\subset X$ and radii $r_j>0$ such that singular fibrations can be constructed on $B_{Rr_j}(y_j)$ and $\bigcup_j B_{r_j}(y_j)=X$. Choose $T_i=\{y_{i,j}\mid 1\le j\le J\}\subset X_i$ with $y_{i,j}\to y_j$, and construct the corresponding nilpotent models. We fix these models and construct a global singular fibration by gluing their local maps.

Volume comparison gives a constant $C(n)$ with the following property.
For every $r\le1$, any $r$-separated subset of $X$ can be partitioned into at most $C(n)$ families, each of which is $10r$-separated.
Since the balls $B_{r_j}(y_j)$ cover $X$, choose $0<h\le1$ such that every $10h$-ball lies in one of these balls.
Let $S=\{x_s\mid1\le s\le s_0\}\subset X$ be a maximal $h$-separated set.
Maximality gives $\bigcup_{x\in S}B_h(x)=X$. Divide $S$ into disjoint subsets $S_a$, $a \le C(n)$, so that each $S_a$ is $10h$-separated. For each $s$, choose $x_{i,s} \in X_i$ converging to $x_s \in X$. Then define $S_i= \{ x_{i,s} | 1 \le s \le s_0\}$ and $S_{i,a}=\{ x_{i,s} | x_s \in S_a\}$. 

For each $1\le s\le s_0$, choose $j_s$ such that $B_{10h}(x_s)\subset B_{r_{j_s}}(y_{j_s})$. We call the nilpotent model at $y_{i,j_s}$ the initial model at $x_{i,s}$ and the corresponding singular fiber bundle map $f_{i,j_s}$ the initial map at $x_{i,s}$. 

We glue the local singular fibrations by induction on
$a\in\{1,\dots,C(n)\}$.
Set $S(a)=\bigcup_{j=1}^a S_j$ and
$S_i(a)=\bigcup_{j=1}^a S_{i,j}$.
At step $a$, we construct a single singular fibration
on $B_{4.5h-2ah/C(n)}(S_i(a))$ that is
$\epsilon_i$-close to the original
GHA from $X_i$ to $X$. For each $x_s\in S(a)$, assign a model based at
$y_{i,\kappa_s}\in T_i$ and a singular fiber bundle
map $f_{i,s}$ to this model.
The map $f_{i,s}$ induces the constructed singular
fibration on $B_{4.5h-2ah/C(n)}(x_{i,s})$.
The assigned maps are compatible wherever these
balls intersect.

For $a=1$, use the initial models and Lemma~\ref{lemma:change scale} to construct local singular fibrations on $B_{4.5h}(x_{i,s})$ for $x_{i,s}\in S_{i,1}$, and set $\kappa_s=j_s$, taking the resulting models and maps as the assigned ones. These balls are disjoint, so no gluing is needed.

Assume the construction is complete at step $a-1$. We describe the $a$-th step.
For sufficiently large $i$, $S_{i,a}$ is $9h$-separated,
so the constructions at its distinct points can be
performed independently.
Fix $x_{i,s}\in S_{i,a}$.
Let $K$ index the points $x_{i,s_k}\in S_i(a-1)$ such that
\[
B_{4.5h-2(a-1)h/C(n)}(x_{i,s_k})
\cap
B_{4.5h-2(a-1)h/C(n)}(x_{i,s})
\ne\emptyset.
\]
Following Lemma~\ref{lemma:glue fibers}, we extend the
previously constructed fibration to the chosen ball $B_{4.5h-2ah/C(n)}(x_{i,s})$ while preserving it on
$B_{4.5h-2ah/C(n)}(x_{i,s_k})$ for every $k\in K$.

Take the initial nilpotent model at $y_{i,j_s}$ for $x_{i,s}$ and the assigned models at $y_{i,\kappa_{s_k}}$ for $x_{i,s_k}$. Choose $j_0\in\{j_s\}\cup\{\kappa_{s_k}:k\in K\}$ such that $r_{j_0}$ is minimal. Apply Lemmas~\ref{lemma:fiber at small scale} and~\ref{lemma:x-fiber at small scale} to represent the initial map $f_{i,j_s}$ in the nilpotent model at $y_{i,j_0}$ on $B_{4.5h}(x_{i,s})$.
Retain the notation $f_{i,j_s}$ for this representation.
For each $k\in K$, likewise represent the assigned map $f_{i,s_k}$ by a compatible singular fiber bundle map in the same model.
Each resulting map is $(\Phi(\delta)+\epsilon_i)$-close to $f_{i,j_s}$.
Retain the notation $f_{i,s_k}$ for these representations.

We next glue $f_{i,s_k}$ and $f_{i,j_s}$ in the model $(M_{i,j_0},\hat{G}_{i,j_0},E_{i,j_0})$. We choose lifts $\hat{x}_{i,s_k},\hat{x}_{i,s} \in \hat{Y}_{i,j_0}$ of $x_{i,s_k}$ and $x_{i,s}$, respectively. Then let $\hat{x}_{i,s_k}'=f_{i,s_k}(\hat{x}_{i,s_k})$ and $\hat{x}_{i,s}'= f_{i,j_s}(\hat{x}_{i,s})$, both in $M_{i,j_0}$. Rescale the model metric to the \(h\)-scale, as in Lemma~\ref{lemma:fiber at small scale}. 

Using Lemma~\ref{lemma:cut-off}, for each $k \in K$, let $\chi_k$ be an $E_{i,j_0}$-invariant cut-off function which is $1$ on the $4.5-(2a-1)/C(n)$-ball of $E_{i,j_0}\hat{x}_{i,s_k}'$, and $0$ outside the $4.5-(2a-1.5)/C(n)$-ball of $E_{i,j_0}\hat{x}_{i,s_k}'$. Let $\psi$ be a smooth non-negative function so that $\psi(t)=1$ if $t \le 1/4$ and $\psi(t)=0$ if $t \ge 1/2$. 

For each $k\in K$, define $\lambda_k:=\chi_k\circ f_{i,s_k}$ on the domain of the previously assigned map and extend it by zero. The cutoff supports lie in the interiors of these domains by the local comparison estimates. Use the weight functions
\[
\begin{aligned}
\frac{\lambda_k}{\psi(\sum_{k\in K}\lambda_k)+\sum_{k\in K}\lambda_k}
&\quad\text{on }f_{i,s_k},\\
\frac{\psi(\sum_{k\in K}\lambda_k)}{\psi(\sum_{k\in K}\lambda_k)+\sum_{k\in K}\lambda_k}
&\quad\text{on }f_{i,j_s}.
\end{aligned}
\]
Taking the pointwise weighted center of mass of these maps gives a $\hat{G}_{i,j_0}$-equivariant almost splitting map $f_{i,s}$. Set $\kappa_s=j_0$ and take $f_{i,s}$ as the assigned map at $x_{i,s}$ in this model.

On the preimage under $f_{i,s_k}$ of the $4.5-(2a-1)/C(n)$-ball about $E_{i,j_0}\hat{x}_{i,s_k}'$, we have $\lambda_k=1$ and hence $\psi(\sum_{k\in K}\lambda_k)=0$, so only previously assigned maps enter the average. These maps are compatible on their overlaps. Since the model cutoffs are $E_{i,j_0}$-invariant, the weights $\lambda_k$ are constant on the elements of the common fiber partition determined by the previously assigned maps. Their center of mass therefore remains compatible with these maps, preserving the induction hypothesis.

The construction in the proof of Lemma~\ref{lemma:glue fibers} yields a singular fibration induced by $f_{i,s}$ over
\[
B_{4.5h-(2a-\tfrac12)h/C(n)}(x_s)
\cup\bigcup_{k\in K}B_{4.5h-(2a-\tfrac12)h/C(n)}(x_{s_k})
\]
that is $\epsilon_i$-close to the original GHA.
As in Lemma~\ref{lemma:glue fibers}, deform this local
singular fibration so that it agrees with the previous
one on $B_{4.5h-2ah/C(n)}(x_{i,s_k})$ for every $k\in K$.
Together with the unchanged fibration outside these
neighborhoods, this completes the induction step.

Each gluing step adds an error of order $\Phi(\delta)$
to the regularity estimates for the models and almost
splitting maps.
Since there are at most $C(n)$ steps, the total error
remains controlled.
After the final step, the domain contains
$B_{2.5h}(S_i)=X_i$ for all sufficiently large $i$,
because the $h$-balls centered at $S$ cover $X$.
The resulting singular fibrations
$\bar F_i:X_i\to X$ prove Theorem~B.

\section{Nilpotent collapsing structure}
We now formulate Theorem~C using the nilpotent models constructed in Section~\ref{sec 6}.
\begin{theorem}\label{Theorem C}
Assume $(X,p,d,\mathcal{H}^n)$ is an $\mathrm{RCD}(-(n-1),n)$ space satisfying the $(\rho,\delta)$ local rewinding Reifenberg condition for some $\rho>0$ and $0<\delta<\delta(n)$. Then there exist positive constants $\epsilon=\epsilon(n,\rho)$
and $h=h(n,\rho)$ such that the $\epsilon$-collapsed region of $X$
admits a cover by balls $\{B_h(x_s):x_s\in X\}$.
Each ball admits a singular fiber bundle map from a covering
neighborhood to a nilpotent model, and these maps are compatible
on overlaps.
\end{theorem}

Unlike Theorem B, Theorem C has no fixed limit space. We first prove a quantitative existence result for nilpotent models. Their fiber partitions can then be made compatible on overlaps by gluing, so that on each overlap, one fiber partition refines the other.

\subsection{Quantitative existence of nilpotent models}

Recall that if $H$ is a closed subgroup of $\mathrm{Isom}(\mathbb{R}^n)$ satisfying property~Q from Section~\ref{sec 6}, then $\dim(H) \le n$. Let $S_n$ be the set of closed subgroups $H$ of $\mathrm{Isom}(\mathbb{R}^n)$ satisfying property~Q and $\dim(H)=n$.
\begin{lemma}\label{lemma:dim n}
Fix $R>10$. Then there exists $\delta_n>0$ such that for every $0<\delta<\delta_n$ the following holds. Assume an $\mathrm{RCD}(-(n-1)\delta, n)$ space $(X,p,G,\mathcal{H}^n)$ is $\delta$-close to $(\mathbb{R}^n,q=0^n,H,\mathcal{H}^n)$ with $H \in S_n$. Suppose $G$ is a discrete group acting freely by isometries. Then $(B_R(p),p,G(p,R))$ is equivariantly isomorphic, on an
$R$-ball, to some $(\hat X,\hat p,\hat G)$.
There is a nilpotent model $(M,N,\hat G,E=N\hat G)$ with
basepoint $\hat p'$ such that $\dim N=n$, $M$ admits a
$\Phi(\delta\mid n,R)$-$C^\infty$-close diffeomorphism to
$\mathbb R^n$, and $\hat G$ is $\Phi(\delta\mid n,R)$-dense in $E$.
There is also a $\hat G$-equivariant almost splitting map
$f:(\hat X,\hat p)\to(M,\hat p')$ defined on $B_R(\hat G\hat p)$. Moreover, the quotient map $\bar{f}:B_{R}(\bar{p}) \to M/E$, where $\bar{p} = \pi(p) \in X/G$, is a singular fibration with infranilmanifold fibers.
\end{lemma} 
\begin{proof}
Unlike in Theorem~\ref{thm:R^n singular fiber bundle}, the comparison group \(H\) is allowed to vary.

Suppose, to the contrary, that the conclusion fails along a sequence with $\epsilon_i\to0$. The spaces $(X_i,p_i,G_i,\mathcal H^n)$ are $\epsilon_i$-close to $(\mathbb R^n,q=0^n,H_i,\mathcal H^n)$, with $H_i\in S_n$. After passing to a subsequence, both $(X_i,p_i,G_i,\mathcal H^n)$ and $(\mathbb R^n,0^n,H_i,\mathcal H^n)$ converge to $(\mathbb R^n,q,H,\mathcal H^n)$. Since $H_i\in S_n$, $H_i=\mathbb{R}^n \rtimes (H_i)_{q}$ where $(H_i)_{q}$ is finite. Then $H$ also contains $\mathbb{R}^n$. By the structure of approximate groups, $H_0$ is nilpotent, so $H=\mathbb{R}^n \rtimes H_q$ where $H_q$ is finite.
The compact isotropy groups $(H_i)_q$ converge to $H_q$; Montgomery--Zippin's theorem identifies them for large $i$, up to orthogonal conjugacy \cite{MontgomeryZippin1942}. Applying the proof of Theorem~\ref{thm:R^n singular fiber bundle} to
\[
(X_i,p_i,G_i,\mathcal{H}^n) \to (\mathbb{R}^n,q=0^n,H,\mathcal{H}^n),
\] 
gives a contradiction.
\end{proof}

We next consider groups of dimension less than $n$.
For $j=n-1,\dots,0$, define $S_j$ inductively as follows.
Consider the closed subgroups $H\le\mathrm{Isom}(\mathbb R^n)$ satisfying property~Q and having dimension $j$.
Include $H$ in $S_j$ precisely when its distance from every closed subgroup of larger dimension satisfying property~Q is at least $\delta_{j+1}/10$.
Under the remaining hypotheses of Lemma~\ref{lemma:dim n},
the same proof gives a constant $0<\delta_j\le\delta_{j+1}/10$
such that, whenever $(X,p,G,\mathcal H^n)$ is $\delta$-close
to $(\mathbb R^n,0^n,H,\mathcal H^n)$ with $H\in S_j$
and $0<\delta<\delta_j$, there is a nilpotent model
$(M,N,\hat G,E=N\hat G)$ for $(X,p,G)$, with basepoint
$\hat p'$ and $\dim N=j$.

The gluing argument in the next subsection has at most $C(n)$ steps. Choose $\delta_n,\ldots,\delta_0$ sufficiently small that the preceding constructions remain valid after the errors from these steps have accumulated.

We can now construct local singular fibrations in the setting of Theorem C.
\begin{lemma}\label{lemma:small szie closeness}
Fix $\rho>0$ and $R>10$. Let $R_0=\min\left\{\sqrt{\delta_0},\frac{\rho}{100R}\right\}$. There exists $0<r_0<R_0$ such that the following holds for sufficiently small $\delta$, depending on $\delta_0$ and $n$. Suppose that $(X,p,G,\mathcal{H}^n)$ is an $\mathrm{RCD}(-(n-1),n)$ space satisfying the $(\rho,\delta)$ Reifenberg condition at $p$ and $G$ is a discrete group acting freely by isometries. Then there exist $0 \le j \le n$ and $r \in [r_0,R_0]$ such that $(r^{-1}X,p,G,\mathcal{H}^n)$ is $\delta_j/2$-close to $(\mathbb{R}^n,q=0^n,H,\mathcal{H}^n)$ with $H \in S_{j}$.
\end{lemma} 
\begin{proof}
Suppose, to the contrary, that there are $r_i\to0$ and spaces $(X_i,p_i,G_i,\mathcal{H}^n)$ for which the conclusion fails with $r_0=r_i$. Passing to a subsequence, we have
\[
(X_i,p_i,G_i,\mathcal{H}^n) \to (X,p,G,\mathcal{H}^n).
\] 
We may take $R_i \to \infty$ such that $R_i \le 1/ \sqrt{r_i}$ and $(R_iX_i,p_i,G_i,\mathcal{H}^n)$ converges to a tangent cone of $(X,p,G,\mathcal{H}^n)$, denoted by $(Y,q,H,\mathcal{H}^n)$.

By Theorem~\ref{theorem:stability of group}, $H$ satisfies property Q. For small $\delta$, $(Y,q,H,\mathcal H^n)$ is $\delta_0/10$-close to some $(\mathbb R^n,0^n,H',\mathcal H^n)$ with $H'$ satisfying property Q; write $j=\dim(H')$. If $H'\in S_j$, convergence gives the required contradiction with $r=R_i^{-1}$ for large $i$. Otherwise, the definition of $S_j$ gives $s>j$ and  $d(H',H'')<\delta_{j+1}/10$ for some $H'' \le \mathrm{Isom}(\mathbb{R}^n)$ with property Q and dimension $s$. If $H''\in S_s$, the same convergence argument gives a contradiction. Otherwise, repeat the comparison with a group of larger dimension satisfying property~Q. Since the dimension strictly increases and is bounded, the process terminates after finitely many steps at some $H''' \in S_{j'}$ where $j' \ge j$. The total error is at most $\sum_{l=0}^{j'} \delta_l/10 <\delta_{j'}/2$, a contradiction.
\end{proof}

Combining Lemma~\ref{lemma:small szie closeness} with the preceding discussion, if $\delta< \delta(n)$ is small enough, we can construct nilpotent models on local universal covers. The index $j$ in Lemma~\ref{lemma:small szie closeness} need not
be unique; we always choose the largest index for which the
conclusion of the lemma holds. Notice that if $\mathcal{H}^n(B_1(\bar{p}))$ is small enough, we always have $j \ge 1$.

In Theorem B, the singular fibration must be $\epsilon_i$-close to the original GHA. Theorem C does not have a fixed limit space, so our construction only needs to control the error $\Phi(\delta)$.

Consider the example $X_i\to[0,\pi/2]$ in Section~\ref{example B}. Fix a sufficiently small $t>0$ and
choose $x_i\in X_i$ converging to $t$. At a lift of $t$, the equivariant tangent group is
$\mathbb R^5$. Our pointwise tangent-model construction therefore gives a regular fibration chart. This regular chart cannot contain the singular fiber over $0$, so its radius must tend to zero as $t\to0$.

For Theorem C, we may instead use a $(t+\epsilon_i)$-GHA that sends $x_i$ to $0$, and construct
the nilpotent model using a fixed endpoint neighborhood.
For sufficiently small $t$ and sufficiently large $i$, the model is then defined near $x_i$ on a ball whose radius is bounded below independently of $t$ and $i$. For fixed $t>0$, the small error bound $\Phi(t+\epsilon_i)$ need not tend to zero as $i\to\infty$. 
\subsection{Gluing nilpotent models}
Let $(X,p,d,\mathcal H^n)$ satisfy the hypotheses of Theorem~C, and let $x,y\in X$ satisfy $d(x,y)\le\rho/10$. Let $\tilde{X}=\tilde{B}(x,\rho/2,\rho)$ be the relative covering space with deck group $G$. By the argument in the last subsection, there are scales $R_0 \ge r_x,r_y\ge r_0$ at which we can construct nilpotent models $(M_x,\hat{x}'_x,N_x,\hat{G}_x,E_x)$ for $B_{Rr_x}(x)$ and $(M_y,\hat{y}'_y,N_y,\hat{G}_y,E_y)$ for $B_{Rr_y}(y)$.

Assume $B_{r_x}(x)\cap B_{r_y}(y)\ne\emptyset$ and $r_x\ge r_y$; let $s:=r_y/r_x\in[r_0/R_0,1]$. The proof of Lemma~\ref{lemma:group id} gives
\[
\theta_{yx}:E_y\longrightarrow E_x,
\]
injective on the local group $E_y(\hat y_y',R-2)$. The inclusion
$\Gamma_y(\hat y_y,R-2)\subset
\Gamma_x^{s^{-1}}(\hat y_x,R)$
follows for sufficiently small $\delta$ by the argument
of Lemma~\ref{lemma:group id}. More precisely, under the common covering identification,
if $g\in\Gamma_y(\hat y_y,R-2)$ is $\epsilon$-close to the identity,
then the corresponding element of $\hat G_x(\hat y_x,R)$ is
$\Phi(\epsilon+\delta)$-close to the identity in $\hat G_x$. 
This convergence is uniform over the admissible models and scale ratios when $n$, $R$, the threshold hierarchy, and the positive lower bound $r_0/R_0$ are fixed.
Suppose this uniform conclusion fails.
We could then choose counterexamples with $\epsilon_i\to0$ and $\delta_i\to0$.
Their scale ratios lie in $[r_0/R_0,1]$.
The argument of Lemma~\ref{lemma:group id} would produce a limit action $h$ that acts nontrivially on $B_4(0^n)$.
For some $t$ with $r_0/R_0\le t\le1$, this action would fix every point of $B_t(q)\subset B_3(0^n)$.
This is impossible.

After retaining only the generators and relations visible at scale $r_y$, we can identify a defining subset of $E_y$ with a local subgroup of $E_x$ by $\theta_{yx}$.
Some elements of the short basis of $\hat{G}_x$ at scale $r_x$ may lie far from the identity at scale $r_y$. The corresponding directions in $E_x$ are therefore discarded when it is viewed at the smaller scale. Using this identification, the construction in the preceding section still applies.

Using $\theta_{yx}$, we define compatibility as in Definition~\ref{defn:compat}. Unlike in Theorem B, we may have $\dim E_y<\dim E_x$. Compatibility means that the fiber partition determined by $f_y$ refines that determined by $f_x$, with the group relations transported by $\theta_{yx}$.

We now use the gluing argument of Theorem B to construct nilpotent models on a cover of the $\epsilon(n,\rho)$-collapsed region of $X$ with compatible fiber partitions. Since no fibers need to be assigned to points of a prescribed limit space, the final deformation is unnecessary.

\begin{proof}[Proof of Theorem C]
Fix $R>10$. Using the constant $C(n)$ from the proof of Theorem~B, we proceed in at most $C(n)$ steps.

Choose $\epsilon(n,\rho)$ such that Lemma~\ref{lemma:small szie closeness} holds for some $j \ge 1$ if $\mathcal{H}^n(B_1(p)) \le \epsilon(n,\rho)$. Let $h=r_0/(10C(n))$ and let $S=\{x_{s}\}$ be a maximal $h$-separated subset of the $\epsilon(n,\rho)$-collapsed region of $X$. Then we can divide $S$ into at most $C(n)$ disjoint subsets $S_a$ such that in each $S_a$, all points are $10h$-separated. By the construction in Section~\ref{sec 6}, each point $x \in S$ has a nilpotent model
$(M_x,N_x,\hat G_x,E_x)$ with basepoint $\hat x_x'$ and a
$\hat G_x$-equivariant almost splitting map
\[ 
f_x:(\hat{X}_x,\hat{x},\hat{G}_x) \to (M_x,\hat{x}'_x,\hat{G}_x,E_x)
\]
which induces a singular fibration on $B_{Rr_x}(x)$, where $r_x\ge r_0$. We use these as the initial models and maps.

Apply the same induction and gluing argument as in Theorem B. This gives assigned models and maps describing the singular fibration on $B_h(x_s)$ for each $x_s \in S$, and these maps are compatible on overlaps. This completes the proof. 
\end{proof}  

\appendix
\section{Nilprogressions at large scales}
This appendix proves Theorem~\ref{thm:large nil}. For our applications, it suffices to consider the case in which $\Gamma$ is bi-Lipschitz isomorphic to $\mathbb{R}^{k}\times\mathbb{Z}^{l}$ equipped with the standard product metric. We outline the argument in this case.

\begin{proof}[Proof sketch in the special case]
Via the given isomorphism, we may work with the standard product metric, changing only the constants in the conclusion.

For each fixed scale $r$, consider the symmetric sets
\[
A_i(r):=
\bar\varphi_i^{-1}(B_r(e_\Gamma))
\cap
\bigl(\bar\varphi_i^{-1}(B_r(e_\Gamma))\bigr)^{-1}.
\]
The uniform covering estimates in $\mathbb{R}^{k}\times\mathbb{Z}^{l}$ give an approximate-group constant $C$ independent of $r$. For each fixed $r$, these sets are strong $C$-approximate groups for all sufficiently large $i$. 

Consequently, the escape norm estimates in
Theorem~\ref{thm:escape} hold with constants independent
of the scale. The argument in \cite{BGT}
gives the required generators and estimates.
\end{proof}

In the general setting, large balls need not be strong $C$-approximate groups with $C$ independent of the radius. The preceding argument therefore does not directly yield uniform escape norm estimates on large scales. 

\subsection{The connected case}
We first treat the case in which $\Gamma$ is connected, simply connected, and nilpotent, with $\dim\Gamma=k$. We show that the nilprogressions in Theorem~\ref{thm:nil} extend to arbitrarily large scales.

For any nilprogression $P=P(u_1,\dots,u_k;C_1,\dots,C_k)$ and $r>0$, we define $rP:=P(u_1,\dots,u_k;rC_1,\dots,rC_k)$. In this subsection, $\Lambda:[1,\infty)\to(0,\infty)$ denotes a function tending to infinity as $r\to\infty$; its choice may change between intermediate estimates.

\begin{theorem}\label{thm:con nil}
In the setting of Theorem~\ref{thm:nil}, let
\[
P_i=P(u_{i,1},\dots,u_{i,k};C_{i,1},\dots,C_{i,k}),
\] 
and assume in addition that $\Gamma$ is a connected, simply connected nilpotent Lie group. Then there is a function $\Lambda$ such that, for every fixed $r\ge1$ and all sufficiently large $i$ depending on $r$, the following hold:
\begin{itemize}
\item $\bar{\varphi}_i(rP_i) \subset B_{(r+1)C}(e_{\Gamma})$, $\bar{\varphi}_i^{-1}(B_{\Lambda(r)}(e_{\Gamma})) \subset rP_i$.

\item For integer exponents satisfying $\lvert n_{i,j}\rvert,\lvert n'_{i,j}\rvert\le rC_{i,j}$ for $1\le j\le k$, the equality $u_{i,1}^{n_{i,1}}\dots u_{i,k}^{n_{i,k}} = u_{i,1}^{n'_{i,1}}\dots u_{i,k}^{n'_{i,k}}$ implies that $n_{i,j}=n_{i,j}'$ for each $j$.
\end{itemize}
\end{theorem}

We first prove the outer containment. Since $\bar{\varphi}_i(P_i)\subset B_C(e_\Gamma)$ and $\bar{\varphi}_i$ is almost homomorphic, we have $\bar{\varphi}_i(rP_i)\subset B_{(r+1)C}(e_\Gamma)$ for all sufficiently large $i$, after increasing $C$ if necessary.

\begin{lemma}\label{lem:contain a large ball}
There is a function $\Lambda$ such that, for every fixed $r\ge1$, the inclusion $\bar{\varphi}_i^{-1}(B_{\Lambda(r)}(e_\Gamma))\subset rP_i$ holds for all sufficiently large $i$ depending on $r$.
\end{lemma}

\begin{proof}
Let $(v_1,\dots,v_k)$ be the limiting strong Malcev basis introduced before Theorem~\ref{thm:Q conv}. By the short-basis construction underlying Theorem~\ref{thm:nil}, there are $r_j>0$ such that $\bar\varphi_i(\{u_{i,j}^m:|m|\le C_{i,j}\})$ converges to $\{\exp_\Gamma(tv_j):|t|\le r_j\}$.

For each fixed $r\ge1$, the sets
$\bar\varphi_i(\{u_{i,j}^m:|m|\le rC_{i,j}\})$
converge to $\{\exp_\Gamma(tv_j):|t|\le rr_j\}$. The associated ordered exponential chart is a diffeomorphism, so, for all sufficiently large $i$,
$\bar{\varphi}_i(rP_i) \cap B_{\Lambda(r)}(e_{\Gamma})$ is $\epsilon_i$-dense in $B_{\Lambda(r)}(e_{\Gamma})$.

For any $g\in\Gamma_i/H_i$ such that $\bar\varphi_i(g)\in B_{\Lambda(r)}(e_\Gamma)$, there exists $h \in rP_i$ so that $d(\bar{\varphi}_i(h)^{-1}\bar{\varphi}_i(g) ,e_{\Gamma}) \le \epsilon_i$.
Then $d(\bar{\varphi}_i(h^{-1}g),e_{\Gamma})\le 2\epsilon_i$. Since $\bar{\varphi}_i^{-1}(B_{\delta}(e_{\Gamma})) \subset P_i$ and $2\epsilon_i<\delta$ for sufficiently large $i$, $h^{-1}g\in P_i$. 

By Theorem~\ref{thm:Q conv}, the inclusions $h\in rP_i$ and $h^{-1}g\in P_i$ imply that there exists $R(r)>0$, depending on $r$ and the Lie algebra structure expressed in the basis $(v_1,\dots,v_k)$, such that $g=h(h^{-1}g)\in R(r)P_i$. Therefore $\bar{\varphi}_i^{-1}(B_{\Lambda(r)}(e_{\Gamma}))\subset R(r)P_i$. Since $R(r)$ is independent of $i$, reparametrizing the progression scale gives the claimed inclusion.
\end{proof}

We next prove the second conclusion of Theorem~\ref{thm:con nil}.
\begin{lemma}\label{lem:proper}
For every fixed $r>0$ and all sufficiently large $i$, the ordered
representation $u_{i,1}^{n_{i,1}}\cdots u_{i,k}^{n_{i,k}}$ is unique
for integer exponents satisfying $|n_{i,j}|\le rC_{i,j}$ for
$1\le j\le k$.
\end{lemma}
\begin{proof}
Assume $u_{i,1}^{n_{i,1}}u_{i,2}^{n_{i,2}}\dots u_{i,k}^{n_{i,k}} = u_{i,1}^{n_{i,1}'}u_{i,2}^{n_{i,2}'}\dots u_{i,k}^{n_{i,k}'}$ for integer exponents satisfying $-rC_{i,j}\le n_{i,j},n_{i,j}' \le rC_{i,j}$.
After passing to a subsequence, we may assume that $\bar{\varphi}_i(u_{i,j}^{n_{i,j}}) \to g_j$ and $\bar{\varphi}_i(u_{i,j}^{n_{i,j}'}) \to g_j'$. Then $g_1g_2\dots g_k = g_1'g_2'\dots g_k'$. Uniqueness in the ordered exponential chart of this strong Malcev basis gives $g_j=g_j'$ for every $j$. Therefore for each $j$, $\frac{n_{i,j}-n_{i,j}'}{C_{i,j}} \to 0$ as $i \to \infty$. Using $Q_i$ in Theorem~\ref{thm:Q conv}, we have
\[
(u_{i,1}^{n_{i,1}}u_{i,2}^{n_{i,2}}\dots u_{i,k}^{n_{i,k}})(u_{i,1}^{n_{i,1}'}u_{i,2}^{n_{i,2}'}\dots u_{i,k}^{n_{i,k}'})^{-1} = u_{i,1}^{n_{i,1}''}u_{i,2}^{n_{i,2}''}\dots u_{i,k}^{n_{i,k}''}.
\]
Since $\frac{n_{i,j}-n_{i,j}'}{C_{i,j}}\to0$ for every $j$, Lie algebra convergence gives $\frac{n_{i,j}''}{C_{i,j}}\to0$ for every $j$. Therefore $u_{i,1}^{n_{i,1}''}u_{i,2}^{n_{i,2}''}\dots u_{i,k}^{n_{i,k}''} \in P_i$. Since this product equals the identity, uniqueness in Theorem~\ref{thm:nil} implies that $n_{i,j}''=0$ for every $j$.

Observe that the $u_{i,1}$ component of the product on the left-hand side of the preceding display is $u_{i,1}^{n_{i,1}-n_{i,1}'}$. Therefore $n_{i,1}=n_{i,1}'$ and $u_{i,2}^{n_{i,2}}\dots u_{i,k}^{n_{i,k}} = u_{i,2}^{n_{i,2}'}\dots u_{i,k}^{n_{i,k}'}$. An induction argument shows that $n_{i,j}=n_{i,j}'$ for each $j$.
\end{proof}

\subsection{The disconnected case}
In this subsection, we prove Theorem~\ref{thm:large nil}.
Let $\Gamma_i':=\{ g \in \Gamma_i/H_i | \bar{\varphi}_i(g) \in \Gamma_0\}$. Choose $\epsilon>0$ sufficiently small that $B_{\epsilon}(e_{\Gamma})\subset\Gamma_0$. Since $\bar{\varphi}_i$ is an $\epsilon_i$-almost homomorphism and $\epsilon_i<\epsilon$ for all sufficiently large $i$, the restrictions $\bar{\varphi}_i:\Gamma_i'\to\Gamma_0$ form a discrete approximation. Then we can apply Theorem~\ref{thm:con nil} to $\bar{\varphi}_i:\Gamma_i' \to \Gamma_0$ to obtain $P(u_{i,1},\dots,u_{i,k};rC_{i,1},\dots,rC_{i,k})$.

On a sufficiently large ball, the local quotient $(\Gamma_i/H_i)/\Gamma_i'$ is naturally identified with a local subgroup of $\Gamma/\Gamma_0$. Thus $\Gamma_i'$ approximates $\Gamma_0$, while the component representatives describe the discrete quotient.

The progression in Theorem~\ref{thm:large nil} formed from the original generators $w_{i,s}$ and $u_{i,j}$ satisfies the required properties except possibly the mixed commutator condition \eqref{e0}; for instance, $[w_{i,1},u_{i,k}]=u_{i,k-1}$ may occur. To obtain the mixed commutator condition, we replace \(u_{i,1},\ldots,u_{i,k}\) with a basis adapted to the component generators. 

We recall the short-basis construction in Theorem~\ref{thm:nil}. Use the strong approximate groups $A_i$ chosen in Theorem~\ref{thm:strong}, with nested working neighborhoods chosen as in \cite[proof of Lemma~9.4]{BGT}. Since the small subgroups have been removed, choose $u_{i,k}\in A_i\setminus\{e\}$ with minimal escape norm. The commutator estimate in Theorem~\ref{thm:escape} shows that $u_{i,k}$ commutes with every sufficiently short generator, and hence with $A_i$.

Write $S_{i,k}:=\min\{m\ge1:u_{i,k}^m\notin A_i\}$ and $\mu_i:=\|u_{i,k}\|_{A_i}=S_{i,k}^{-1}$.
The normalized progression $P(u_{i,k};S_{i,k})$ tends to a nontrivial local one-parameter subgroup. Taking the quotient by this central direction and comparing escape norms yield the usual dimension induction; see \cite[Lemmas~9.4--9.5]{BGT} and \cite[Theorem~7.1 and Lemma~7.2]{Zamora2020}.
We first construct $u_{i,k}'$ and then proceed by induction on the dimension in the same way.

We first observe that for each $1 \le j \le k-1$ and large $i$, $\|u_{i,j+1}\|_{A_i} \le C \|u_{i,j}\|_{A_i}$. Moreover, there exists a sufficiently small constant $c>0$, independent of $i$ and $j$, such that every $g\in A_i$ satisfying $\|g\|_{A_i}<c\|u_{i,j}\|_{A_i}$ belongs to the subgroup generated by $u_{i,j+1},\dots,u_{i,k}$. After passing to a subsequence, we may assume that the limits
$\lim_{i\to\infty}
\frac{\|u_{i,j+1}\|_{A_i}}{\|u_{i,j}\|_{A_i}}$
exist for every $j$ and are at most $C$.
Choose $0=s_0<s_1<\dots<s_b=k$ so that
$\{s_l\mid1\le l\le b-1\}$ consists precisely of the indices $j$
for which
$\lim_{i\to\infty}
\frac{\|u_{i,j+1}\|_{A_i}}{\|u_{i,j}\|_{A_i}}=0$.

Since commutators have much smaller escape norms, for any $1\le l_1\le l_2\le b$ and indices $j_1,j_2$ satisfying $s_{l_1-1}+1\le j_1\le s_{l_1}$ and $s_{l_2-1}+1\le j_2\le s_{l_2}$, we have
\[
[u_{i,j_1},u_{i,j_2}] \in P(u_{i,s_{l_2}+1},\dots,u_{i,k};\frac{CC_{i,s_{l_2}+1}}{C_{i,j_1}C_{i,j_2}},\dots,\frac{CC_{i,k}}{C_{i,j_1}C_{i,j_2}}).
\]
In particular, every element of $P_{i,1}=P(u_{i,s_{b-1}+1},\dots,u_{i,k};C_{i,s_{b-1}+1},\dots,C_{i,k})$ commutes with every element of $A_i$.

We next show that $P_{i,1}$ is preserved up to a fixed dilation under conjugation by $w_{i,s}$ for each $1\le s\le l$. 

\begin{lemma}\label{lemma:conj}
For all sufficiently large $i$ and each $1 \le s \le l$,
we have
\[
\|g\|_{A_i}/C
\le \|w_{i,s}gw_{i,s}^{-1}\|_{A_i}
\le C\|g\|_{A_i}
\]
for every $g \in A_i$.
\end{lemma}
\begin{proof}
Fix \(s\). We prove the left-hand inequality; the right-hand inequality follows similarly. Suppose that there are nontrivial $h_i\in A_i$ with
\[
\frac{\|w_{i,s}h_iw_{i,s}^{-1}\|_{A_i}}{\|h_i\|_{A_i}}\longrightarrow0.
\]
Let $T_i$ be the first escape time of $h_i$ from $A_i$. After passing to a subsequence, we may assume that $\bar\varphi_i(h_i^{T_i})\longrightarrow h\ne e_\Gamma$.
On the other hand, the product estimate for the escape norm gives $\|(w_{i,s}h_iw_{i,s}^{-1})^{T_i}\|_{A_i}\to0$. Hence $\bar\varphi_i((w_{i,s}h_iw_{i,s}^{-1})^{T_i})\to e_\Gamma$. Almost multiplicativity and $\bar\varphi_i(w_{i,s})\to g_s$ also imply that this sequence converges to $g_shg_s^{-1}\ne e_\Gamma$. This is a contradiction.
\end{proof}

\begin{corollary}\label{cor:inv}
For large $i$ and each $1 \le s \le l$, $w_{i,s}P_{i,1}w_{i,s}^{-1} \subset CP_{i,1}$.
\end{corollary}

The images of $P_{i,1}$ converge to a central local subgroup $G_1\le\Gamma_0$ of dimension $k-s_{b-1}$. Its Lie algebra is spanned by the limiting tail vectors $v_{s_{b-1}+1},\dots,v_k$. By Corollary~\ref{cor:inv}, this subspace is invariant under every $\operatorname{Ad}(g_s)$.

By Theorem \ref{thm:escape}, we can use the same escape time $S_{i,k}$ to normalize every element whose escape norm is comparable to $\mu_i$.
\begin{lemma}\label{lem:power}
For every fixed $C\ge1$, there exist constants $c_C,R_C>0$
with the following property. For every sequence
$h_i\in P_{i,1}\setminus\{e\}$ satisfying $\|h_i\|_{A_i}\le C\mu_i$,
we have
\[
c_C\le d_\Gamma\bigl(\bar\varphi_i(h_i^{S_{i,k}}),e_\Gamma\bigr)
\le R_C
\]
for all sufficiently large $i$.
\end{lemma}

Next we find a nontrivial element of controlled escape norm in the last central block that commutes with every component generator.
\begin{lemma}\label{lemma:commute}
There is a constant $C>0$ such that, for all sufficiently large $i$,
one can choose a nontrivial $h_i\in P_{i,1}$ that commutes with every
$w_{i,s}$, $1\le s\le l$, and satisfies
$\|h_i\|_{A_i}\le C\|u_{i,k}\|_{A_i}$.
\end{lemma}
\begin{proof}
We start with $u_{i,k}$. If $w_{i,s}u_{i,k}=u_{i,k}w_{i,s}$ for all large $i$ and $1 \le s \le l$, we simply take $h_i=u_{i,k}$.

Otherwise, after passing to a subsequence, choose $\sigma_1$ such that $h_{i,1}:=[u_{i,k},w_{i,\sigma_1}]$ is nontrivial. Lemma~\ref{lemma:conj} and the product estimate give $\|h_{i,1}\|_{A_i}\le C\mu_i$. The short-basis ordering therefore places $h_{i,1}$ and the relevant conjugates in the last central block. Consequently, $h_{i,1}^{S_{i,k}}=[u_{i,k}^{S_{i,k}},w_{i,\sigma_1}]$. Passing to a subsequence, write $\bar\varphi_i(u_{i,k}^{S_{i,k}})\to u_\infty$. Lemma~\ref{lem:power} shows that the limit $[u_\infty,g_{\sigma_1}]$ of $\bar\varphi_i(h_{i,1}^{S_{i,k}})$ is nontrivial.

If $h_{i,1}$ commutes with every $w_{i,s}$, take $h_i=h_{i,1}$. Otherwise repeat the construction with a nontrivial $h_{i,2}:=[h_{i,1},w_{i,\sigma_2}]$. Its escape norm is at most $C_2\mu_i$, it remains in the central block, and its $S_{i,k}$-th power has the nontrivial limit $[[u_\infty,g_{\sigma_1}],g_{\sigma_2}]$. At every stage the same argument applies, with a fixed larger constant. The process must stop within $k+l$ steps, since $\Gamma$ is nilpotent of step at most $k+l$. The resulting nontrivial element has the required properties.
\end{proof}

To construct the adapted basis, first choose its central
generator $u'_{i,k}$.
Let
\[
A_i':=\{g\in A_i:gw_{i,s}=w_{i,s}g
                    \text{ for every }1\le s\le l\}.
\]
Choose $u'_{i,k}$ minimizing $\|g\|_{A_i}$ over
$A_i'\setminus\{e\}$.
Lemma~\ref{lemma:commute} gives
$\|u'_{i,k}\|_{A_i}\le C\mu_i$.

By \cite{BGT}, after passing to a subsequence if necessary, the image of $P(u_{i,k}';S_{i,k})$ under $\bar\varphi_i$ converges to a $1$-dimensional local subgroup.
Choose $r_i\to\infty$ sufficiently slowly, for example with $r_i\ll R_i$ and $r_i\ll\epsilon_i^{-1}$.
The image of $P(u_{i,k}';r_iS_{i,k})$ under $\bar\varphi_i$ then converges to a $1$-dimensional subgroup $L$ of $\Gamma$.
The bound $\|u'_{i,k}\|_{A_i}\le C\mu_i$ and the commutator estimate imply that $u'_{i,k}$ commutes with every short generator of $A_i$.
By the definition of $A_i'$, it also commutes with every $w_{i,s}$.
These facts imply that the limiting one-parameter subgroup $L$ is central in $\Gamma$.

Take the quotient local group $\overline{\Gamma}_i:=(\Gamma_i/H_i)/P(u_{i,k}';r_iS_{i,k})$, and let $\varphi_i^L:\overline{\Gamma}_i\to\Gamma/L$ be the induced quotient approximation. Let $\overline{A}_i\subset\overline{\Gamma}_i$ be the image of $A_i$ under the projection. As in \cite{BGT}, we show that the escape norm in $\overline{A}_i$ is comparable with the norm in $A_i$.

\begin{lemma}\label{lemma:comp}
There exists $C>0$ such that, for all sufficiently large $i$, every $\bar h\in\overline A_i$ has a lift $h\in A_i$ satisfying $\|h\|_{A_i}\le C\|\bar h\|_{\overline A_i}$.
\end{lemma}
For every lift \(h\in A_i\) of \(\bar h\), the reverse inequality $\|\bar{h}\|_{\overline{A}_i} \le \|h\|_{A_i}$ follows directly from the definition.
\begin{proof}
Choose a sequence of nontrivial elements
$\bar h_i\in\overline A_i$ and, for each $i$, a lift
$h_i\in A_i$ minimizing $\|h_i\|_{A_i}$ among all lifts
of $\bar h_i$ in $A_i$.
Assume that
$\|\bar h_i\|_{\overline A_i}/\|h_i\|_{A_i}\to0$.

Let $T_i$ be the first escape time of $h_i$ from $A_i$. The argument of \cite[Lemma~9.5]{BGT} shows that the limiting one-parameter subgroup obtained from the normalized powers of \(h_i\) lies in $L$. If $\|h_i\|_{A_i}\ge\|u'_{i,k}\|_{A_i}$, the same lifting argument provides an integer $m_i$, with $|m_i|\le CS_{i,k}$, such that $\|h_i(u'_{i,k})^{m_i}\|_{A_i}<\|h_i\|_{A_i}$. This element is another lift of $\bar h_i$, contradicting minimality.

It remains to consider $\|h_i\|_{A_i}<\|u'_{i,k}\|_{A_i}\le C\mu_i$. Then $h_i\notin A_i'$ and $h_i$ lies in the last central block. After passing to a subsequence, $\bar\varphi_i(h_i^{S_{i,k}})\longrightarrow g\in L$. Choose $s$ such that $[h_i,w_{i,s}]$ is nontrivial. Its escape norm is at most a fixed multiple of $\mu_i$. As in Lemma~\ref{lemma:commute}, the central power identity and Lemma~\ref{lem:power} show that $[g,g_s]\ne e_\Gamma$. This contradicts the centrality of $L$.
\end{proof}

Lemma~\ref{lemma:comp} also implies that the quotient approximation has no small subgroups, so the induction requires no further small-subgroup quotient.

\begin{proof}[Proof of Theorem~\ref{thm:large nil}]
If $k=0$, then $\bar{\varphi}_i$ is a local group isomorphism,
and the theorem follows. Assume $k\ge1$ and suppose that
Theorem~\ref{thm:large nil} holds with $k$ replaced by $k-1$. The preceding construction gives $u'_{i,k}$ central in $A_i$ and commuting with every $w_{i,s}$. Apply the induction hypothesis to the local quotient approximation $\varphi_i^L:\overline{\Gamma}_i\to\Gamma/L$ introduced above. We obtain $\bar u'_{i,1},\dots,\bar u'_{i,k-1}\in\overline A_i$ satisfying the conclusions of Theorem~\ref{thm:nil}
and condition~\eqref{e0}. For each $1 \le j \le k-1$, choose $u_{i,j}' \in \Gamma_i/H_i$ as a preimage of $\bar{u}_{i,j}' \in \overline{\Gamma}_i$ with a minimal escape norm. Then $u_{i,1}',\dots,u_{i,k}'$ satisfy \eqref{e0}. Lemma~\ref{lemma:comp} and the proof of \cite[Lemma 9.4]{BGT} show that $u_{i,1}',\dots,u_{i,k}'$ satisfy the properties in Theorem~\ref{thm:nil} as well. 
The same arguments as in Lemmas~\ref{lem:contain a large ball}
and~\ref{lem:proper} show that
$w_{i,1},\dots,w_{i,l},u_{i,1}',\dots,u_{i,k}'$
satisfy the conclusions of Theorem~\ref{thm:large nil}.
\end{proof}
\bibliographystyle{plain}
\bibliography{bib}

@article{PanWang2021,
  author        = {Pan, Jiayin and Wang, Jikang},
  title         = {Some topological results of {Ricci} limit spaces},
  journal       = {Transactions of the American Mathematical Society},
  year          = {2022},
  volume        = {375},
  number        = {12},
  pages         = {8445--8464},
  doi           = {10.1090/tran/8549},
  archiveprefix = {arXiv},
  eprint        = {2103.11344},
  primaryclass  = {math.DG},
  url           = {https://doi.org/10.1090/tran/8549},
}

@article{Palais1961,
  author        = {Palais, Richard S.},
  title         = {On the existence of slices for actions of non-compact {Lie} groups},
  journal       = {Annals of Mathematics},
  year          = {1961},
  volume        = {73},
  number        = {2},
  pages         = {295--323},
  doi           = {10.2307/1970335},
  issn          = {0003-486X},
  url           = {https://doi.org/10.2307/1970335},
}

@article{CheegerColding1997,
  author        = {Cheeger, Jeff and Colding, Tobias H.},
  title         = {On the structure of spaces with {Ricci} curvature bounded below. {I}},
  journal       = {Journal of Differential Geometry},
  year          = {1997},
  volume        = {46},
  number        = {3},
  pages         = {406--480},
  doi           = {10.4310/jdg/1214459974},
  url           = {https://doi.org/10.4310/jdg/1214459974},
}

@article{SormaniWei2004,
  author        = {Sormani, Christina and Wei, Guofang},
  title         = {Universal covers for {Hausdorff} limits of noncompact spaces},
  journal       = {Transactions of the American Mathematical Society},
  year          = {2004},
  volume        = {356},
  number        = {3},
  pages         = {1233--1270},
  doi           = {10.1090/S0002-9947-03-03412-3},
  issn          = {0002-9947},
  url           = {https://doi.org/10.1090/S0002-9947-03-03412-3},
}

@article{CheegerColding2000a,
  author        = {Cheeger, Jeff and Colding, Tobias H.},
  title         = {On the structure of spaces with {Ricci} curvature bounded below. {II}},
  journal       = {Journal of Differential Geometry},
  year          = {2000},
  volume        = {54},
  number        = {1},
  pages         = {13--35},
  doi           = {10.4310/jdg/1214342145},
  url           = {https://doi.org/10.4310/jdg/1214342145},
}

@article{CheegerNaber2015,
  author        = {Cheeger, Jeff and Naber, Aaron},
  title         = {Regularity of {Einstein} manifolds and the codimension 4 conjecture},
  journal       = {Annals of Mathematics},
  year          = {2015},
  volume        = {182},
  number        = {3},
  pages         = {1093--1165},
  doi           = {10.4007/annals.2015.182.3.5},
  issn          = {0003-486X},
  url           = {https://doi.org/10.4007/annals.2015.182.3.5},
}

@misc{KapovitchWilking2011,
  author        = {Kapovitch, Vitali and Wilking, Burkhard},
  title         = {Structure of fundamental groups of manifolds with {Ricci} curvature bounded below},
  year          = {2011},
  howpublished  = {arXiv:1105.5955},
  archiveprefix = {arXiv},
  eprint        = {1105.5955},
  primaryclass  = {math.DG},
  url           = {https://arxiv.org/abs/1105.5955},
}

@article{Fukaya1986,
  author        = {Fukaya, Kenji},
  title         = {Theory of convergence for {Riemannian} orbifolds},
  journal       = {Japanese Journal of Mathematics. New Series},
  year          = {1986},
  volume        = {12},
  number        = {1},
  pages         = {121--160},
  doi           = {10.4099/math1924.12.121},
  issn          = {0289-2316},
  url           = {https://doi.org/10.4099/math1924.12.121},
}

@article{FukayaYamaguchi1992,
  author        = {Fukaya, Kenji and Yamaguchi, Takao},
  title         = {The fundamental groups of almost nonnegatively curved manifolds},
  journal       = {Annals of Mathematics},
  year          = {1992},
  volume        = {136},
  number        = {2},
  pages         = {253--333},
  doi           = {10.2307/2946606},
  issn          = {0003-486X},
  url           = {https://doi.org/10.2307/2946606},
}

@article{HKRX2020,
  author        = {Huang, Hongzhi and Kong, Lingling and Rong, Xiaochun and Xu, Shicheng},
  title         = {Collapsed manifolds with {Ricci} bounded covering geometry},
  journal       = {Transactions of the American Mathematical Society},
  year          = {2020},
  volume        = {373},
  number        = {11},
  pages         = {8039--8057},
  doi           = {10.1090/tran/8177},
  issn          = {0002-9947},
  url           = {https://doi.org/10.1090/tran/8177},
}

@article{Huang2020,
  author        = {Huang, Hongzhi},
  title         = {Fibrations and stability for compact group actions on manifolds with local bounded {Ricci} covering geometry},
  journal       = {Frontiers of Mathematics in China},
  year          = {2020},
  volume        = {15},
  number        = {1},
  pages         = {69--89},
  doi           = {10.1007/s11464-020-0824-2},
  issn          = {1673-3452},
  url           = {https://doi.org/10.1007/s11464-020-0824-2},
}

@article{Wang2021,
  author        = {Wang, Jikang},
  title         = {{Ricci} limit spaces are semi-locally simply connected},
  journal       = {Journal of Differential Geometry},
  year          = {2024},
  volume        = {128},
  number        = {3},
  pages         = {1301--1314},
  doi           = {10.4310/jdg/1729092461},
  url           = {https://doi.org/10.4310/jdg/1729092461},
}

@article{Fukaya1987,
  author        = {Fukaya, Kenji},
  title         = {Collapsing {Riemannian} manifolds to ones of lower dimensions},
  journal       = {Journal of Differential Geometry},
  year          = {1987},
  volume        = {25},
  number        = {1},
  pages         = {139--156},
  doi           = {10.4310/jdg/1214440728},
  url           = {https://doi.org/10.4310/jdg/1214440728},
}

@article{Fukaya1988,
  author        = {Fukaya, Kenji},
  title         = {A boundary of the set of the {Riemannian} manifolds with bounded curvatures and diameters},
  journal       = {Journal of Differential Geometry},
  year          = {1988},
  volume        = {28},
  number        = {1},
  pages         = {1--21},
  doi           = {10.4310/jdg/1214442157},
  url           = {https://doi.org/10.4310/jdg/1214442157},
}

@article{Fukaya1989,
  author        = {Fukaya, Kenji},
  title         = {Collapsing {Riemannian} manifolds to ones with lower dimension {II}},
  journal       = {Journal of the Mathematical Society of Japan},
  year          = {1989},
  volume        = {41},
  number        = {2},
  pages         = {333--356},
  doi           = {10.2969/jmsj/04120333},
  url           = {https://doi.org/10.2969/jmsj/04120333},
}

@article{Anderson1992,
  author        = {Anderson, Michael T.},
  title         = {{Hausdorff} perturbations of {Ricci}-flat manifolds and the splitting theorem},
  journal       = {Duke Mathematical Journal},
  year          = {1992},
  volume        = {68},
  number        = {1},
  pages         = {67--82},
  doi           = {10.1215/S0012-7094-92-06803-7},
  url           = {https://doi.org/10.1215/S0012-7094-92-06803-7},
}

@misc{Rong2019,
  author        = {Rong, Xiaochun},
  title         = {A new proof of the {Gromov}'s theorem on almost flat manifolds},
  year          = {2019},
  howpublished  = {arXiv:1906.03377},
  archiveprefix = {arXiv},
  eprint        = {1906.03377},
  primaryclass  = {math.DG},
  url           = {https://arxiv.org/abs/1906.03377},
}

@article{CFG,
  author        = {Cheeger, Jeff and Fukaya, Kenji and Gromov, Mikhael},
  title         = {Nilpotent structures and invariant metrics on collapsed manifolds},
  journal       = {Journal of the American Mathematical Society},
  year          = {1992},
  volume        = {5},
  number        = {2},
  pages         = {327--372},
  doi           = {10.1090/S0894-0347-1992-1126118-X},
  issn          = {0894-0347},
  url           = {https://doi.org/10.1090/S0894-0347-1992-1126118-X},
  language      = {English},
}

@article{Gromov1978,
  author        = {Gromov, Mikhael},
  title         = {Almost flat manifolds},
  journal       = {Journal of Differential Geometry},
  year          = {1978},
  volume        = {13},
  number        = {2},
  pages         = {231--241},
  doi           = {10.4310/jdg/1214434488},
  url           = {https://doi.org/10.4310/jdg/1214434488},
}

@article{Ruh1982,
  author        = {Ruh, Ernst A.},
  title         = {Almost flat manifolds},
  journal       = {Journal of Differential Geometry},
  year          = {1982},
  volume        = {17},
  number        = {1},
  pages         = {1--14},
  doi           = {10.4310/jdg/1214436698},
  url           = {https://doi.org/10.4310/jdg/1214436698},
}

@article{GuSan2019,
  author        = {Guijarro, Luis and Santos-Rodr{\'i}guez, Jaime},
  title         = {On the isometry group of {$\mathrm{RCD}^{*}(K,N)$}-spaces},
  journal       = {Manuscripta Mathematica},
  year          = {2019},
  volume        = {158},
  number        = {3--4},
  pages         = {441--461},
  doi           = {10.1007/s00229-018-1010-7},
  url           = {https://doi.org/10.1007/s00229-018-1010-7},
}

@article{Sosa2018,
  author        = {Sosa, Gerardo},
  title         = {The isometry group of an {$\mathrm{RCD}^{*}$} space is {Lie}},
  journal       = {Potential Analysis},
  year          = {2018},
  volume        = {49},
  number        = {2},
  pages         = {267--286},
  doi           = {10.1007/s11118-017-9656-4},
  url           = {https://doi.org/10.1007/s11118-017-9656-4},
}

@article{MondinoWei2019,
  author        = {Mondino, Andrea and Wei, Guofang},
  title         = {On the universal cover and the fundamental group of an {$\mathrm{RCD}^{*}(K,N)$}-space},
  journal       = {Journal f{\"u}r die reine und angewandte Mathematik},
  year          = {2019},
  volume        = {753},
  pages         = {211--237},
  doi           = {10.1515/crelle-2016-0068},
  url           = {https://doi.org/10.1515/crelle-2016-0068},
}

@article{Gigli2013,
  author        = {Gigli, Nicola},
  title         = {The splitting theorem in non-smooth context},
  journal       = {Memoirs of the American Mathematical Society},
  year          = {2026},
  volume        = {317},
  number        = {1609},
  doi           = {10.1090/memo/1609},
  archiveprefix = {arXiv},
  eprint        = {1302.5555},
  primaryclass  = {math.MG},
  url           = {https://doi.org/10.1090/memo/1609},
}

@article{SantosZamora,
  author        = {Santos-Rodr{\'i}guez, Jaime and Zamora-Barrera, Sergio},
  title         = {On fundamental groups of {RCD} spaces},
  journal       = {Journal f{\"u}r die reine und angewandte Mathematik},
  year          = {2023},
  volume        = {799},
  pages         = {249--286},
  doi           = {10.1515/crelle-2023-0027},
  url           = {https://doi.org/10.1515/crelle-2023-0027},
}

@article{Wang2022,
  author        = {Wang, Jikang},
  title         = {{$\mathrm{RCD}^{*}(K,N)$} spaces are semi-locally simply connected},
  journal       = {Journal f{\"u}r die reine und angewandte Mathematik},
  year          = {2024},
  volume        = {806},
  pages         = {1--7},
  doi           = {10.1515/crelle-2023-0058},
  url           = {https://doi.org/10.1515/crelle-2023-0058},
}

@article{Zamora2020,
  author        = {Zamora, Sergio},
  title         = {Limits of almost homogeneous spaces and their fundamental groups},
  journal       = {Groups, Geometry, and Dynamics},
  year          = {2024},
  volume        = {18},
  number        = {3},
  pages         = {761--798},
  doi           = {10.4171/GGD/792},
  url           = {https://doi.org/10.4171/GGD/792},
}

@article{Rong2022,
  author        = {Rong, Xiaochun},
  title         = {Collapsed manifolds with local {Ricci} bounded covering geometry},
  journal       = {Science China Mathematics},
  year          = {2026},
  volume        = {69},
  number        = {6},
  pages         = {1585--1614},
  doi           = {10.1007/s11425-025-2473-y},
  url           = {https://doi.org/10.1007/s11425-025-2473-y},
}

@article{BGT,
  author        = {Breuillard, Emmanuel and Green, Ben and Tao, Terence},
  title         = {The structure of approximate groups},
  journal       = {Publications Math{\'e}matiques de l'Institut des Hautes {\'E}tudes Scientifiques},
  year          = {2012},
  volume        = {116},
  pages         = {115--221},
  doi           = {10.1007/s10240-012-0043-9},
  url           = {https://doi.org/10.1007/s10240-012-0043-9},
}

@book{Thurston,
  author        = {Thurston, William P.},
  title         = {Three-dimensional geometry and topology, volume 1},
  year          = {1997},
  series        = {Princeton Mathematical Series},
  volume        = {35},
  publisher     = {Princeton University Press},
  isbn          = {978-0-691-08304-9},
  note          = {Edited by Silvio Levy},
}

@article{CheegerNaberJiang2021,
  author        = {Cheeger, Jeff and Jiang, Wenshuai and Naber, Aaron},
  title         = {Rectifiability of singular sets of noncollapsed limit spaces with {Ricci} curvature bounded below},
  journal       = {Annals of Mathematics},
  year          = {2021},
  volume        = {193},
  number        = {2},
  pages         = {407--538},
  doi           = {10.4007/annals.2021.193.2.2},
  url           = {https://doi.org/10.4007/annals.2021.193.2.2},
}

@article{CheegerGromoll1971,
  author        = {Cheeger, Jeff and Gromoll, Detlef},
  title         = {The splitting theorem for manifolds of nonnegative {Ricci} curvature},
  journal       = {Journal of Differential Geometry},
  year          = {1971},
  volume        = {6},
  number        = {1},
  pages         = {119--128},
  doi           = {10.4310/jdg/1214430220},
  url           = {https://doi.org/10.4310/jdg/1214430220},
}

@article{CG2,
  author        = {Cheeger, Jeff and Gromov, Mikhael},
  title         = {Collapsing {Riemannian} manifolds while keeping their curvature bounded. {II}},
  journal       = {Journal of Differential Geometry},
  year          = {1990},
  volume        = {32},
  number        = {1},
  pages         = {269--298},
  doi           = {10.4310/jdg/1214445047},
  url           = {https://doi.org/10.4310/jdg/1214445047},
}

@article{CG1,
  author        = {Cheeger, Jeff and Gromov, Mikhael},
  title         = {Collapsing {Riemannian} manifolds while keeping their curvature bounded. {I}},
  journal       = {Journal of Differential Geometry},
  year          = {1986},
  volume        = {23},
  number        = {3},
  pages         = {309--346},
  doi           = {10.4310/jdg/1214440117},
  url           = {https://doi.org/10.4310/jdg/1214440117},
}

@misc{Wang2023,
  author        = {Wang, Jikang},
  title         = {On the limit of simply connected manifolds with discrete isometric cocompact group actions},
  year          = {2023},
  howpublished  = {arXiv:2307.07658},
  archiveprefix = {arXiv},
  eprint        = {2307.07658},
  primaryclass  = {math.DG},
  url           = {https://arxiv.org/abs/2307.07658},
}

@book{BuserKarcher,
  author        = {Buser, Peter and Karcher, Hermann},
  title         = {{Gromov}'s almost flat manifolds},
  year          = {1981},
  series        = {Ast{\'e}risque},
  volume        = {81},
  publisher     = {Soci{\'e}t{\'e} math{\'e}matique de France},
  address       = {Paris},
  url           = {https://www.numdam.org/item/AST_1981__81__1_0/},
}

@article{Xu2023,
  author        = {Xu, Shicheng},
  title         = {Precompactness of domains with lower {Ricci} curvature bound under {Gromov}-{Hausdorff} topology},
  journal       = {Journal f{\"u}r die reine und angewandte Mathematik},
  year          = {2025},
  volume        = {828},
  pages         = {149--173},
  doi           = {10.1515/crelle-2025-0060},
  url           = {https://doi.org/10.1515/crelle-2025-0060},
}

@article{ZamoraZhu2024,
  author        = {Zamora, Sergio and Zhu, Xingyu},
  title         = {Topological rigidity of small {$\mathrm{RCD}(K,N)$} spaces with maximal rank},
  journal       = {Nonlinear Analysis},
  year          = {2026},
  volume        = {267},
  eid           = {114046},
  doi           = {10.1016/j.na.2025.114046},
  note          = {Article 114046},
  url           = {https://doi.org/10.1016/j.na.2025.114046},
}

@article{HondaPeng2024,
  author        = {Honda, Shouhei and Peng, Yuanlin},
  title         = {A note on the topological stability theorem from {RCD} spaces to {Riemannian} manifolds},
  journal       = {Manuscripta Mathematica},
  year          = {2023},
  volume        = {172},
  number        = {3--4},
  pages         = {971--1007},
  doi           = {10.1007/s00229-022-01418-7},
  url           = {https://doi.org/10.1007/s00229-022-01418-7},
}

@Article{BS2023,
  author  = {Bru{\`e}, Elia and Semola, Daniele},
  journal = {Communications on Pure and Applied Mathematics},
  title   = {Constancy of the dimension for {RCD(K,N)} spaces via regularity of {Lagrangian} flows},
  year    = {2020},
  number  = {6},
  pages   = {1141--1204},
  volume  = {73},
  doi     = {10.1002/cpa.21849},
}

@article{BNS,
  author        = {Bru{\`e}, Elia and Naber, Aaron and Semola, Daniele},
  title         = {Boundary regularity and stability for spaces with {Ricci} bounded below},
  journal       = {Inventiones Mathematicae},
  year          = {2022},
  volume        = {228},
  number        = {2},
  pages         = {777--891},
  doi           = {10.1007/s00222-021-01092-8},
  url           = {https://doi.org/10.1007/s00222-021-01092-8},
}

@article{DeGigli2019,
  author        = {De Philippis, Guido and Gigli, Nicola},
  title         = {Non-collapsed spaces with {Ricci} curvature bounded from below},
  journal       = {Journal de l'{\'E}cole polytechnique --- Math{\'e}matiques},
  year          = {2018},
  volume        = {5},
  pages         = {613--650},
  doi           = {10.5802/jep.80},
  url           = {https://doi.org/10.5802/jep.80},
}

@article{DSZZ2023,
  author        = {Deng, Qin and Santos-Rodr{\'i}guez, Jaime and Zamora, Sergio and Zhao, Xinrui},
  title         = {{Margulis} lemma on {$\mathrm{RCD}(K,N)$} spaces},
  journal       = {Analysis \& PDE},
  year          = {2025},
  volume        = {18},
  number        = {10},
  pages         = {2481--2536},
  doi           = {10.2140/apde.2025.18.2481},
  url           = {https://doi.org/10.2140/apde.2025.18.2481},
}

@article{AH18,
  author        = {Ambrosio, Luigi and Honda, Shouhei},
  title         = {Local spectral convergence in {$\mathrm{RCD}^{*}(K,N)$} spaces},
  journal       = {Nonlinear Analysis},
  year          = {2018},
  volume        = {177},
  number        = {A},
  pages         = {1--23},
  doi           = {10.1016/j.na.2017.04.003},
  url           = {https://doi.org/10.1016/j.na.2017.04.003},
}

@book{KirbySiebenmann,
  author        = {Kirby, Robion C. and Siebenmann, Laurence C.},
  title         = {Foundational essays on topological manifolds, smoothings, and triangulations},
  year          = {1977},
  series        = {Annals of Mathematics Studies},
  volume        = {88},
  publisher     = {Princeton University Press},
  url           = {https://www.jstor.org/stable/j.ctt1b9s024},
}

@misc{HHWZ,
  author        = {Huang, Hongzhi and Huang, Xian-Tao and Wang, Jikang and Zhu, Xingyu},
  title         = {Fibrations, the first {Betti} number, and almost nonnegative {Ricci} curvature},
  year          = {2026},
  howpublished  = {arXiv:2605.24380},
  archiveprefix = {arXiv},
  eprint        = {2605.24380},
  primaryclass  = {math.DG},
  url           = {https://arxiv.org/abs/2605.24380},
}

@Article{Wang2024,
  author  = {Wang, Jikang},
  journal = {Advances in Mathematics},
  title   = {Non-collapsed {RCD} spaces with locally bounded covering geometry},
  year    = {2026},
  pages   = {111373},
  volume  = {507},
  doi     = {10.1016/j.aim.2026.111373},
}

@article{Wilking2000,
  author        = {Wilking, Burkhard},
  title         = {Rigidity of group actions on solvable {Lie} groups},
  journal       = {Mathematische Annalen},
  year          = {2000},
  volume        = {317},
  number        = {2},
  pages         = {195--237},
  doi           = {10.1007/s002089900091},
  url           = {https://doi.org/10.1007/s002089900091},
}

@article{HuangHuang2024,
  author        = {Huang, Hongzhi and Huang, Xian-Tao},
  title         = {Almost splitting maps, transformation theorems and smooth fibration theorems},
  journal       = {Advances in Mathematics},
  year          = {2024},
  volume        = {457},
  eid           = {109914},
  doi           = {10.1016/j.aim.2024.109914},
  note          = {Article 109914},
  url           = {https://doi.org/10.1016/j.aim.2024.109914},
}

@article{Harvey2016,
  author        = {Harvey, John},
  title         = {Convergence of isometries, with semicontinuity of symmetry of {Alexandrov} spaces},
  journal       = {Proceedings of the American Mathematical Society},
  year          = {2016},
  volume        = {144},
  number        = {8},
  pages         = {3507--3515},
  doi           = {10.1090/proc/12994},
  url           = {https://doi.org/10.1090/proc/12994},
}

@book{Dekimpe,
  author        = {Dekimpe, Karel},
  title         = {Almost-{Bieberbach} groups: affine and polynomial structures},
  year          = {1996},
  series        = {Lecture Notes in Mathematics},
  volume        = {1639},
  publisher     = {Springer},
  address       = {Berlin, Heidelberg},
  doi           = {10.1007/BFb0094472},
  isbn          = {978-3-540-61899-7},
  url           = {https://doi.org/10.1007/BFb0094472},
}

@article{MazurRongWang,
  author        = {Mazur, Marcin and Rong, Xiaochun and Wang, Yusheng},
  title         = {{Margulis} lemma for compact {Lie} groups},
  journal       = {Mathematische Zeitschrift},
  year          = {2008},
  volume        = {258},
  number        = {2},
  pages         = {395--406},
  doi           = {10.1007/s00209-007-0178-4},
  url           = {https://doi.org/10.1007/s00209-007-0178-4},
}

@Article{Siebenmann1972,
  author  = {Siebenmann, L. C.},
  journal = {Commentarii Mathematici Helvetici},
  title   = {Deformation of homeomorphisms on stratified sets},
  year    = {1972},
  pages   = {123--163},
  volume  = {47},
}

@Article{MontgomeryZippin1942,
  author  = {Montgomery, Deane and Zippin, Leo},
  journal = {Bulletin of the American Mathematical Society},
  title   = {A theorem on {Lie} groups},
  year    = {1942},
  number  = {6},
  pages   = {448--452},
  volume  = {48},
  doi     = {10.1090/S0002-9904-1942-07699-3},
}

@Article{Breuillard2014,
  author  = {Breuillard, Emmanuel},
  journal = {Groups, Geometry, and Dynamics},
  title   = {Geometry of locally compact groups of polynomial growth and shape of large balls},
  year    = {2014},
  number  = {3},
  pages   = {669--732},
  volume  = {8},
  doi     = {10.4171/GGD/244},
}

@Article{KPT2010,
  author  = {Kapovitch, Vitali and Petrunin, Anton and Tuschmann, Wilderich},
  journal = {Annals of Mathematics},
  title   = {Nilpotency, almost nonnegative curvature, and the gradient flow on {Alexandrov} spaces},
  year    = {2010},
  number  = {1},
  pages   = {343--373},
  volume  = {171},
  doi     = {10.4007/annals.2010.171.343},
  url     = {https://annals.math.princeton.edu/2010/171-1/p06},
}

@Book{Palais1960,
  author    = {Palais, Richard S.},
  publisher = {American Mathematical Society},
  title     = {The classification of {$G$}-spaces},
  year      = {1960},
  address   = {Providence, RI},
  number    = {36},
  series    = {Memoirs of the American Mathematical Society},
  doi       = {10.1090/memo/0036},
}

@Article{Huang2024,
  author  = {Huang, Hongzhi},
  journal = {International Mathematics Research Notices},
  title   = {A finite topological type theorem for open manifolds with non-negative {Ricci} curvature and almost maximal local rewinding volume},
  year    = {2024},
  number  = {10},
  pages   = {8568--8591},
  volume  = {2024},
  doi     = {10.1093/imrn/rnad311},
}
\end{document}